\documentclass[a4paper,fleqn]{cas-sc}

\usepackage[numbers,sort&compress]{natbib}

\def\tsc#1{\csdef{#1}{\textsc{\lowercase{#1}}\xspace}}
\tsc{WGM}
\tsc{QE}

\usepackage{algorithm}
\usepackage{algpseudocode}

\usepackage{bm}
\usepackage{mathtools}
\usepackage{placeins}

\usepackage{cleveref}

\hypersetup{
    hidelinks,
    pdftitle={Efficient second-order physics-preserving schemes  for incompressible flows with variable density and viscosity},
    pdfauthor={Baolin Kuang, Lanxin Sun, Rihui Lan, and Hongfei Fu}
}

\newcommand{\maA}{\mathcal{A}}
\newcommand{\maB}{\mathcal{B}}
\newcommand{\maC}{\mathcal{C}}

\newcommand{\maO}{\mathcal{O}}
\newcommand{\maN}{\mathcal{N}}

\newcommand{\maL}{\mathcal{L}}
\newcommand{\maV}{\mathcal{V}}
\newcommand{\mrs}{\mathrm{s}}

\newtheorem{theorem}{Theorem}[section]

\newtheorem{proposition}[theorem]{Proposition}
\newproof{proof}{Proof}

\newdefinition{example}{Example}[section]

\newdefinition{remark}{Remark}[section]

\begin{document}

\let\WriteBookmarks\relax
\def\floatpagepagefraction{1}
\def\textpagefraction{.001}

\shorttitle{Efficient physics-preserving schemes for incompressible flows}
\shortauthors{B. Kuang et al.}

\title[mode=title]{Efficient second-order physics-preserving schemes for incompressible flows with variable density and viscosity}



%


\author[1]{Baolin Kuang}
\fnmark[1]
\ead{kuang_bl124@163.com}
\credit{Writing – review \& editing, Writing – original draft, Software,  Methodology, Formal analysis, Conceptualization}
\affiliation[1]{
    organization={School of Science, Qingdao University of Technology},
    city={Qingdao},
    postcode={266520},
    state={Shandong},
    country={China}
}

\author[2]{Lanxin Sun}%
\fnmark[1]
\ead{sunlanxin@stu.ouc.edu.cn}
\credit{Writing – original draft, Software,  Visualization, Methodology}

\author[2,3]{Rihui Lan}
\ead{lanrihui@ouc.edu.cn}
\credit{Writing – review \& editing, Methodology, Funding  acquisition,  Conceptualization}

\author[2,3]{Hongfei Fu}
\cormark[1]
\ead{fhf@ouc.edu.cn}
\credit{Writing – review \& editing,  Methodology, Funding  acquisition, Conceptualization}
\affiliation[2]{
    organization={School of Mathematical Sciences, Ocean University of China},
    city={Qingdao},
    postcode={266100},
    state={Shandong},
    country={China}
}
\affiliation[3]{
    organization={Laboratory of Marine Mathematics, Ocean University of China},
    city={Qingdao},
    postcode={266100},
    state={Shandong},
    country={China}
}

 \cortext[1]{Corresponding author}

\fntext[1]{The two authors contribute equally to the article.}


\begin{abstract}
This paper proposes a novel second-order prediction--filter--projection framework for incompressible flows with variable density and viscosity. In this framework, a constrained $L^2$-projection is enforced to unconditionally preserve the prescribed density bounds and discrete mass conservation, while a generalized-Laplacian reformulation combined with the D-CSAV approach is used to discretize the variable-viscosity momentum equation. In particular, a density-weighted second-order time filter is proposed, which eliminates the additional uniform-in-time bound restriction imposed on the time difference quotient of the discrete density in existing stability analyses. The resulting scheme is shown to satisfy an unconditional modified-energy dissipation law, without the viscosity-related restrictions imposed in earlier work. Then, an additional scalar auxiliary variable and its compatible temporal discretization are introduced to construct a variable-step prediction--filter--projection scheme. This scheme is also bound-preserving and mass-conserving, and satisfies an unconditional modified-energy dissipation law. Moreover, the modified energy is much closer to the ``original'' kinetic energy, differing from it only through the auxiliary-variable corrections. Based on the variable-step scheme, a reliable energy-based adaptive time-stepping strategy is developed. Finally, numerical experiments demonstrate the accuracy, physics-preserving properties, and efficiency of the proposed schemes.
\end{abstract}




\maketitle












\section{Introduction}\label{sec:intro}
Incompressible flows are ubiquitous in various scientific, engineering, and biological applications \cite{kwak2005computational}. As a classical model of viscous incompressible flow, the constant-density, constant-viscosity incompressible Navier--Stokes equations have been extensively studied from both the mathematical and numerical perspectives \cite{temam2001,constantin1988navier,doanjcp2025,liShen2020SAVMAC,wangmathcomp2000,linjcp2019}. Nevertheless, many flows of practical interest are inherently variable-density, and in multi-phase settings the constituent fluids may even exhibit different viscosities  \cite{shapiroDrikakis2005,dong2012,elAmrani2024}. For such applications, constant-density, constant-viscosity models are insufficient to capture the essential physics. In this paper, we focus on the numerical simulation of the following variable-density, variable-viscosity incompressible Navier--Stokes system \cite{espinozaM2AN2025,lundgrenNazarov2024}:
\begin{align}
   &  \rho_t+\nabla\cdot(\rho\bm{u})=0, && \text{in } \Omega\times(0,T], \label{model:rho} \\
   & \rho\bigl(\partial_t\bm{u}+(\bm{u}\cdot\nabla)\bm{u}\bigr)
      -\nabla\cdot\bigl(2\nu\nabla^{\mrs}\bm{u}\bigr)+\nabla p
    =\bm{f}, && \text{in } \Omega\times(0,T],  \label{model:u} \\
   &  \nabla\cdot\bm{u}=0, && \text{in } \Omega\times(0,T].\label{model:div}
\end{align}
Here, \(\Omega\subset\mathbb{R}^2\) denotes a bounded Lipschitz domain with boundary $\partial\Omega$; $\rho>0$, \(\bm u\), and \(p\) represent the density, velocity, and pressure, respectively. The function $\bm{f}=\bm{f}(\bm{x},t)$ is a prescribed body force per unit volume. The dynamic viscosity $\nu=\nu(\bm{x},t)>0$ may depend on the density and other flow variables through a prescribed constitutive relation. In general, the viscosity is assumed to be uniformly bounded; that is,
\[
       0<\nu_{\min} \leq \nu(\bm{x},t) \leq \nu_{\max} <\infty.
\]
The symmetric gradient of $\bm{u}$ is denoted by
$
   \nabla^{\mrs}\bm{u}:=\frac{1}{2}\bigl(\nabla\bm{u} + \nabla^\top \bm{u}\bigr).
$
The system is supplemented with the initial conditions
 \begin{equation}\label{model:ic}      
  \rho(\bm{x},0)=\rho^0(\bm{x}),\qquad \bm{u}(\bm{x},0)=\bm{u}^0(\bm{x}), \qquad \text{in}~ \Omega, 
 \end{equation}
 and the homogeneous no-slip boundary condition
 \begin{equation}\label{model:bc}
\bm u=\bm 0, \qquad \text{on}~ \partial\Omega. 
 \end{equation}
 
Under the boundary condition \eqref{model:bc}, the system preserves the
initial density bounds \cite{guermond2000}: if
\[
  0<\rho_{\min} \leq \rho^0(\bm{x}) \leq \rho_{\max}, \qquad \text{in}~\Omega,
\]
then the bounds 
\begin{equation} \label{model:bounds}
    \rho_{\min}\leq \rho(\bm{x},t)\leq\rho_{\max}
\end{equation}
hold for all $(\bm{x},t)\in\Omega\times[0,T]$. It also conserves the total mass, i. e.,
\begin{equation}\label{model:MassConserv}
  \frac{d}{dt}M[\rho(t)]=0, \qquad
   M[\rho(t)]:= \int_{\Omega}\rho(\bm{x},t)\,d\bm{x},
\end{equation}
and moreover, for $\bm{f}=\bm 0$, it satisfies the kinetic-energy dissipation law
\begin{equation}\label{model:Energy}
   \frac{d}{dt}E(t) = -2  \left\|\sqrt{\nu}\,\nabla^{\mrs}\bm u \right\|^2 \leq 0, \qquad
    E(t)=E[\rho(t),\bm u(t)]:= \frac12 \left\|\sqrt{\rho}\,\bm u\right\|^2.
\end{equation}
Therefore, a desirable numerical method should preserve these properties simultaneously at the discrete level. Otherwise, discretizations of the variable-density incompressible Navier--Stokes system may produce nonphysical solutions and may even render the discretization of the momentum equation \eqref{model:u}--\eqref{model:div} ill-posed.

Despite the importance of density bounds, relatively few studies have specifically addressed density-bound preservation for the variable-density incompressible Navier--Stokes equations.  Liu and Walkington \cite{liu2007convergence} employed a piecewise-constant upwind discontinuous Galerkin (DG) discretization for the density transport equation, whose monotonicity preserves the initial density bounds. The resulting density approximation is, however, only formally first-order accurate in space. Li et al. \cite{liJCP2021} developed a high-order DG method equipped with a bound-preserving limiter, for which the prescribed density bounds are ensured under a suitable CFL-type condition. To guarantee the strict lower density bound, Li et al. \cite{liJCP2024} introduced a transformed density variable
$$ 
    \rho=\widetilde{\sigma}^{\,2}+\rho_{\min},
 $$
and developed two efficient energy-stable schemes based on the scalar auxiliary variable (SAV) approach \cite{shenXuYang2018}. However, the preservation of the upper density bound in \eqref{model:bounds} was not considered. A natural approach to enforcing physical constraints is to project the solution onto the constrained manifold~\cite{tong2024pnp}. The underlying idea is to project a predicted solution, which typically fails to preserve positivity and mass conservation, onto an admissible set characterized by these two properties. Since the projection step is a convex optimization problem, it can be characterized by the well-known Karush–Kuhn–Tucker (KKT) conditions, which are also closely related to the Lagrange multiplier approach \cite{cheng2022I,cheng2022II}. To the best of our knowledge, no existing work has considered a density-projection approach for the variable-density incompressible Navier--Stokes system.

Energy-stable discretizations have been extensively studied for the classical (constant density and viscosity) Navier–Stokes equations and related coupled incompressible flow systems \cite{linjcp2019, liShen2020SAVMAC, linNiYangDong2020, yang2021IJNME, yangSISC2021, kuang2025unconditional, kuangFuLi2026, zhang2024unified}, particularly through the SAV approach combined with the zero-energy-contribution (ZEC) idea \cite{yang2021}. However, a spatially varying viscosity introduces an additional numerical difficulty in the momentum equation: the viscous stress couples the velocity components and, even under the incompressibility constraint, cannot be simplified to the standard Laplacian form. An efficient numerical treatment must therefore not only decouple the variable-viscosity contribution but also preserve the underlying viscous dissipation structure. A straightforward way to decouple the velocity components is to split the stress-divergence term into two parts, namely
\begin{equation*}\label{nu_SD}
      \nabla\cdot\bigl(2\nu\nabla^{\mrs}\bm u\bigr) 
       = \nabla\cdot(\nu\nabla\bm u) + \nabla\cdot\bigl(\nu\nabla^\top\bm u\bigr),   
\end{equation*}
and then treat the transpose-gradient term explicitly and the first term semi-implicitly (with $\nu$ evaluated explicitly). However, such a treatment may suffer from poor stability for strongly varying viscosity fields. To improve stability of such implicit-explicit (IMEX) treatments, Barrenechea et al. \cite{barrenecheaSISC2024} systematically studied the generalized-Laplacian (GL) reformulation
\begin{equation}\label{nu_GL}
      \nabla\cdot\bigl(2\nu\nabla^{\mrs}\bm u\bigr) = \nabla\cdot(\nu\nabla\bm u) + (\nabla^\top\bm u)\nabla\nu, 
\end{equation}
which allows the componentwise diffusion term to be treated semi-implicitly and the coupling term explicitly. Espinoza-Contreras et al. \cite{espinozaM2AN2025} subsequently extended this strategy to the variable-density, variable-viscosity incompressible flows, and developed some constant-step first- and second-order linearized IMEX schemes with provable unconditional stability in time. However, their energy estimates for the IMEX schemes retain an explicit dependence on $\|\nabla\nu\|_{L^\infty(\Omega)}$, owing to the explicit treatment of the viscosity-gradient coupling term. This motivates us to develop an energy-stability-compatible treatment of the coupling term that removes such dependence. 
 
Another challenge for the variable-density, variable-viscosity incompressible flows is the construction of an efficient variable-step second-order physics-preserving scheme, which provides a reliable foundation for developing an efficient adaptive time-stepping strategy for stable long-term simulations.  The time filter technique \cite{guzel2018time,decaria2020time,songCMAME2026} improves the temporal accuracy of backward Euler schemes to second order at low cost, while remaining A-stable and readily compatible with time adaptivity. Recently, Li et al. \cite{li2023filtered} developed some variable-step and adaptive time-filtered algorithms for the variable-density, constant-viscosity incompressible Navier--Stokes equations. However, the corresponding energy estimate therein relies on assumptions of uniform lower and upper bounds on the numerical density, as well as a uniform-in-time estimate of its discrete time-difference quotient (cf. \eqref{assump:e1}). These assumptions are imposed in the numerical analysis rather than guaranteed by the scheme itself. This motivates a variable-step second-order formulation that enforces the density bounds numerically while avoiding the additional density-variation assumption in the energy analysis.

This paper aims to present a second-order physics-preserving prediction--filter--projection framework for the variable-density, variable-viscosity incompressible flows. First, a bound-preserving and mass-conserving $L^2$-projection strategy is developed in the density update. Second, a novel density-weighted time filter, together with a recently proposed decoupling constant scalar auxiliary variable (D-CSAV) approach \cite{zhangjcp2025} and the GL reformulation \eqref{nu_GL}, is adopted to treat the momentum equation. The main algorithmic and theoretical contributions of the proposed framework are summarized as follows:
\begin{itemize}
\item  It unconditionally preserves the density bounds and discrete mass conservation, thereby yielding a physically meaningful density solution and eliminating the previous density-bound assumption in the stability analysis.

\item It yields an unconditional modified-energy dissipation law without the explicit dependence on $\|\nabla\nu\|_{L^\infty(\Omega)}$ appearing in previous analyses. Moreover, an additional SAV and its compatible temporal discretization in the variable-step scheme enable a closed stability analysis on nonuniform grids and yield a modified energy much closer to the ``original'' energy. 

\item It is highly efficient: the density, velocity--pressure, and SAV are solved sequentially and are mutually decoupled. In particular, the auxiliary variables in both the constant-step and variable-step schemes are solved explicitly, and, as discussed in \eqref{proj-clamp}, the projection step admits an  explicit formula, except for an efficient scalar nonlinear iteration for the scalar Lagrange multiplier $\xi^{n+1}$ in \eqref{eq:F-xi-cont}. Moreover, the variable-step-based adaptive strategy substantially reduces the computational cost. 

\item Within this framework, a novel second-order density-weighted time filter is proposed for the momentum equation; it eliminates the additional uniform-in-time bound restriction imposed on the time difference quotient of the discrete density required in previous velocity-filtered analyses, while retaining the same low-cost structure as the classical time filter.
\end{itemize}

The remainder of this paper is organized as follows. Section~\ref{sec:2} presents the reformulated system and its modified-energy structure. Section~\ref{sec:scheme} develops the constant-step prediction--filter--projection scheme and establishes its physics-preserving properties. Section~\ref{sec:variable-scheme} extends the framework to variable time steps, derives the corresponding modified-energy law using discrete orthogonal convolution (DOC) kernels, and presents the energy-based adaptive time-stepping strategy. Section~\ref{sec:test} presents numerical results, and Section~\ref{sec:conclusion} concludes the paper.

\section{Model reformulation}\label{sec:2} 
For constant-density incompressible flow with $\rho\equiv\rho_0>0$, the convective term makes no contribution to the kinetic-energy balance, since
\begin{equation*}
      \int_{\Omega}\rho_0(\bm u\cdot\nabla)\bm v\cdot\bm v\,d\bm{x}=0.
\end{equation*}
For variable-density flow, however, this cancellation no longer holds in general, which makes it difficult to establish an energy estimate at the discrete level.

To restore an analogous energy structure, following Guermond and Quartapelle \cite{guermond2000}, we introduce $\sigma:=\sqrt{\rho}$ and, using the density equation \eqref{model:rho}, rewrite the time-derivative term in the momentum equation as
\begin{equation}\label{model:re_u}
     \rho\bm{u}_t = \sigma(\sigma\bm{u})_t -\frac{1}{2}\rho_t\bm{u}  
         = \sigma(\sigma\bm{u})_t +\frac{1}{2} \nabla\cdot(\rho\bm{u}) \bm{u}.
\end{equation}
For convenience, we introduce the following simplified differential operators 
\[
       N(\rho,\bm{u},\bm{v}):= (\rho\bm{u}\cdot\nabla)\bm{v}
           +\frac{1}{2}\nabla\cdot(\rho\bm{u})\,\bm{v},\qquad
     L(\nu,\bm{u},\bm{v}):= \nabla\cdot(\nu\nabla\bm{u})
          +(\nabla^\top \bm{v} )\nabla\nu,
\]
for $\bm{u},\bm{v} \in H_0^1(\Omega)^d$, and for given $\rho$ and $\nu$, define 
the associated forms:
\[
    \maN(\rho,\bm{u},\bm{v};\bm{w})
       :=\bigl(N(\rho,\bm{u},\bm{v}),\bm{w}\bigr),
       \]
       \[
    \maL(\nu,\bm{u},\bm{v};\bm{w})
      :=\bigl(L(\nu,\bm{u},\bm{v}),\bm{w}\bigr)
       = -(\nu\nabla\bm{u},\nabla\bm{w})
          +\bigl((\nabla^\top \bm{v} )\nabla\nu,\bm{w}\bigr),
\]
for $\bm{u},\bm{v},\bm{w} \in H_0^1(\Omega)^d$. Note that integration by parts directly gives
\begin{equation}\label{eq0:tri}
  \maN(\rho,\bm{u},\bm v;\bm v)
    = \int_{\Omega} (\rho\bm{u}\cdot\nabla)\bm v\cdot\bm v\,d\bm{x}
        +\frac{1}{2} \int_{\Omega}\nabla\cdot(\rho\bm{u})\,|\bm v|^2\,d\bm{x}    =0.
\end{equation}
Thus, the modified convective term makes no contribution to the kinetic-energy balance, a property commonly referred to as the ``zero-energy-contribution'' property \cite{yang2021,yang2021IJNME}.

Next, following the D-CSAV framework proposed in \cite{zhangjcp2025}, we introduce a viscous D-CSAV equation:
\begin{equation}\label{sav:M}
   \alpha_t = Q \left[ 
    -2\bigl(\nu\nabla^{\mrs}\bm{u}, \nabla^{\mrs}\bm{u}\bigr)
     -\alpha\,\maL(\nu,\bm{u},\bm{u};\bm{u})\right], \quad 0<t\le T;
   \qquad \alpha(0)= \alpha_0,
\end{equation}
with $Q>0$ a user-prescribed stabilization parameter. 
By the identity \eqref{nu_GL} and the homogeneous boundary condition \eqref{model:bc}, we have
\begin{equation*}
   \maL(\nu,\bm u,\bm u;\bm u)
   =-2\bigl(\nu\nabla^{\mrs}\bm u,\nabla^{\mrs}\bm u\bigr)
   =:\maV(t),
\end{equation*}
and hence, \eqref{sav:M} reduces to
\begin{equation}\label{sav:alpha-linear}
   \alpha_t=Q(1-\alpha)\maV(t),
\end{equation}
whose unique solution with $\alpha(0)=1$ is
$ 
   \alpha(t)\equiv 1.
$ 

Now, the original system \eqref{model:rho}--\eqref{model:div} can be equivalently reformulated as
\begin{align}
  & \rho_t+\nabla\cdot(\rho\bm{u})=0, \label{sav:eq:rho}\\
  & \sigma(\sigma\bm{u})_t +N(\rho,\bm{u},\bm{u})
    -\alpha\,L(\nu,\bm{u},\bm{u}) +\nabla p 
    =\bm{f},  \label{sav:eq:sigma}\\
  & \nabla\cdot\bm{u}=0, \label{sav:eq:div}\\
  & \alpha_t = Q \bigl[ 
    -2\bigl(\nu\nabla^{\mrs}\bm{u}, \nabla^{\mrs}\bm{u}\bigr)
     -\alpha\,\maL(\nu,\bm{u},\bm{u};\bm{u})\bigr],\quad \alpha(0)=1.\label{sav:eq}
\end{align}

The following theorem establishes an equivalent modified kinetic-energy dissipation law for the reformulated system.
\begin{theorem} \label{thm:mEnergy} Let $\bm{f}=\bm{0}$. Then the solutions of the reformulated system \eqref{sav:eq:rho}--\eqref{sav:eq} satisfy the following modified-energy dissipation law:
\begin{equation}\label{mEnergy:e1}
  \frac{d}{dt}E_{\mathrm{m}}[\rho,\bm{u},\alpha]
  = -2\bigl\|\sqrt{\nu}\nabla^{\mrs}\bm{u}\bigr\|^2 \leq 0,
  \qquad  
  E_{\mathrm{m}}[\rho,\bm{u},\alpha]
  :=\frac{1}{2}\|\sigma\bm{u}\|^2+\frac{\alpha}{Q}.
\end{equation}
\end{theorem}
\begin{proof} Taking the $L^2$ inner product of the momentum equation \eqref{sav:eq:sigma} with $\bm{u}$, eliminating the pressure term via \eqref{sav:eq:div} and the convective term via \eqref{eq0:tri}, we arrive at
\begin{equation} \label{mth:1}
    \frac{1}{2}\frac{d}{dt} \|\sigma \bm{u}\|^2 
      -  \alpha \,\maL(\nu, \bm{u}, \bm{u}; \bm{u})   = 0.
\end{equation}
Finally, dividing the auxiliary equation \eqref{sav:eq} by $Q$, and adding the resulting equation to \eqref{mth:1} yields the desired energy identity.  \qed
\end{proof}

\begin{remark} At the continuous level, $E_{\mathrm m}[\rho,\bm u,\alpha]=E[\rho,\bm u]+1/Q$; hence the modified energy differs from the original one only by the constant $1/Q$ and obeys the same dissipation law. In the next section, we shall construct an efficient second-order-in-time numerical scheme based on \eqref{sav:eq:rho}--\eqref{sav:eq} that unconditionally preserves a discrete counterpart of the modified energy \eqref{mEnergy:e1}. Meanwhile, the scheme is also designed to enforce the initial density bounds \eqref{model:bounds} and to preserve mass conservation \eqref{model:MassConserv} in the discrete setting. 
\end{remark}

\section{A constant-step prediction--filter--projection scheme}\label{sec:scheme}
Let $\tau=T/N$ and $t^n=n\tau$, $0\le n\le N$. Let $\psi^n$
denote an approximation to $\psi(t^n)$. Define
\[
  \delta\psi^{n} := \psi^n-\psi^{n-1}, \quad \delta^2\psi^{n} := \psi^{n+1}-2\psi^n+\psi^{n-1}, \quad \psi^{*,n+1}:= 2\psi^n-\psi^{n-1},\quad \text{for}~ n\geq 1.
\]
Moreover, define $g(w):=(\rho_{\max}-w)(w-\rho_{\min})$,
so that $g(w)\ge0$ if and only if $\rho_{\min}\leq w\leq\rho_{\max}$.

\subsection{Second-order numerical scheme}
Firstly, the density is updated by first‑order prediction, second‑order time filtering, and projection onto the admissible set 
\begin{equation} \label{adm:set}
    \maC:=\bigl\{v \in L^2(\Omega) \mid g(v) \ge 0\quad \text{a.e. in }\Omega; \, (v,1)=(\rho^0,1)\bigr\}.
\end{equation}
The three steps are specified as follows.

\medskip
\noindent\textbf{Step 1 (Density prediction).}
Compute the predicted density $\overline{\rho}^{n+1}$ from
\begin{equation}\label{density:predict}
    \frac{\overline{\rho}^{n+1}-\rho^n}{\tau}
     + \nabla\cdot\bigl(\overline{\rho}^{n+1}\bm{u}^{*,n+1}\bigr) =0.
\end{equation}
 
\noindent\textbf{Step 2 (Time filtering).}
Apply the time filter to $\overline{\rho}^{n+1}$ to obtain the second-order approximation $\widehat{\rho}^{n+1}$:
\begin{equation}\label{density:filter}
   \widehat{\rho}^{n+1}  = \overline{\rho}^{n+1} 
       - \frac{1}{3} (\overline{\rho}^{n+1} - 2 \rho^n + \rho^{n-1}).
\end{equation}

\noindent\textbf{Step 3 (Bound-preserving and mass-conserving projection).}
Project $\widehat{\rho}^{n+1}$ onto the admissible set $\maC$:
\begin{equation}\label{density:projection}
   \rho^{n+1} = \mathbb{P}_{\maC}\widehat{\rho}^{n+1}:= \operatorname*{arg\,min}_{v\in \maC} \frac{1}{2}\bigl\|v-\widehat{\rho}^{n+1}\bigr\|^2.
\end{equation}

In fact, the projection \eqref{density:projection} is a convex optimization problem, and can be characterized by the following KKT conditions:
\begin{align}
    & \rho^{n+1} =
  \widehat{\rho}^{n+1}  +\lambda^{n+1}g'\bigl(\rho^{n+1}\bigr)
  -\xi^{n+1},  \label{2th:rhohat}\\
   & \lambda^{n+1}g\bigl(\rho^{n+1}\bigr)  =0,
  \quad  g\bigl(\rho^{n+1}\bigr)\geq 0,
  \quad  \lambda^{n+1}\geq 0,  \label{2th:comple}\\
   & (\rho^{n+1},1)=(\rho^{0},1).  \label{2th:mass}
\end{align}
Here $\lambda^{n+1}=\lambda^{n+1}(\bm{x})$ and
$\xi^{n+1}\in\mathbb R$ are the Lagrange multipliers associated with the
bound and mass constraints, respectively. Define $z^{n+1}:=\widehat{\rho}^{n+1}-\xi^{n+1}$.
Then, the complementary condition \eqref{2th:comple} yields an explicit formula for $\rho^{n+1}$ \cite{cheng2022II,tong2024pnp}:
\begin{equation}\label{proj-clamp}
  \rho^{n+1} = \max\left\{ \rho_{\min}, \min\left\{ \rho_{\max}, z^{n+1} \right\} \right\}.
\end{equation}
Substituting it into \eqref{2th:mass} yields that $ \xi^{n+1}\in\mathbb{R}$ is the solution of the scalar equation
\begin{equation}\label{eq:F-xi-cont}
  F(\xi)
  := \int_{\Omega} \max\bigl\{\rho_{\min}, \min\bigl\{ \rho_{\max}, \widehat{\rho}^{n+1}-\xi \bigr\}\bigr\}\,d\bm{x} - \int_{\Omega}\rho^0\,d\bm{x}=0.
\end{equation}
Note that the admissible set \eqref{adm:set} is nonempty, since $v=\rho^0$ is always feasible. Moreover, the function $F(\cdot)$ is continuous and nonincreasing, which can be solved by the semismooth Newton method \cite{facchinei2003vol2} or the secant method \cite{cheng2022II,KUANG2025} in only a few iterations.

Secondly, the velocity and pressure pairs are advanced in a two‑step procedure: first‑order prediction and second‑order time filtering.

\medskip
\noindent\textbf{Step 4 (Velocity--pressure prediction).} Set
$\sigma^{n+1}=\sqrt{\rho^{n+1}}$. Compute the predicted velocity $\widehat{\bm{u}}^{n+1}$ and pressure $p^{n+1}$ from the first-order velocity--pressure system: 
\begin{align}
 & \sigma^{n+1}  \frac{\sigma^{n+1}\widehat{\bm{u}}^{n+1}-\sigma^n\bm{u}^n }{\tau}
  + N\bigl( \rho^{n+1}, \bm{u}^{*,n+1}, \widehat{\bm{u}}^{n+1} \bigr)
    -\alpha^{*,n+1} L\bigl(\nu^{n+1},  \widehat{\bm{u}}^{n+1},\bm{u}^{*,n+1}\bigr)\nonumber\\
   & \qquad \qquad\qquad\qquad\quad\quad+\nabla p^{n+1}  = \bm{f}^{n+1},  \label{2th:up}\\
 & \nabla\cdot\widehat{\bm{u}}^{n+1}=0.
  \label{2th:div}
\end{align}

\medskip
\noindent\textbf{Step 5 (Weighted-velocity time filtering).}  Apply the time filter to the weighted velocity $\bm{w}=\sigma\bm{u}$ by
\begin{equation}\label{2th:TF:u}
   \sigma^{n+1}\bm{u}^{n+1} 
   = \sigma^{n+1}\widehat{\bm{u}}^{n+1}
      -\frac{1}{3}\left(\sigma^{n+1}\widehat{\bm{u}}^{n+1} -2\sigma^n\bm{u}^n +\sigma^{n-1}\bm{u}^{n-1}\right).
\end{equation}

Finally, the auxiliary variable $\alpha$ is updated by the following explicit second‑order prediction--then--filter procedure.

\medskip
\noindent\textbf{Step 6 (D-CSAV prediction).}
Once $\rho^{n+1}$ and $\bm{u}^{n+1}$ are available, we compute the predicted auxiliary variable $\widehat{\alpha}^{n+1}$ explicitly from  
\begin{equation}\label{2th:Mhat}
 \frac{\widehat{\alpha}^{n+1}-\alpha^n}{\tau}
   = Q\bigl[ -2\bigl(\nu^{n+1}\nabla^{\mrs}\bm{u}^{n+1},\nabla^{\mrs}\bm{u}^{n+1}\bigr) 
   -\alpha^{*,n+1}\maL\bigl(\nu^{n+1}, \widehat{\bm{u}}^{n+1},\bm{u}^{*,n+1};
      \widehat{\bm{u}}^{n+1}\bigr) \bigr].
\end{equation}

\noindent\textbf{Step 7 (D-CSAV time filtering).}
Apply the time filter to $\widehat{\alpha}^{n+1}$ to obtain $\alpha^{n+1}$:
\begin{equation}\label{2th:TF:M}
  \alpha^{n+1} = \widehat{\alpha}^{n+1}-\frac{1}{3}\left(\widehat{\alpha}^{n+1}-2\alpha^n+\alpha^{n-1}\right).
\end{equation}

\begin{remark} The auxiliary variable $\alpha$ serves to compensate for the explicit treatment of the viscosity-gradient coupling in \eqref{2th:up}, while preserving the viscous dissipation in the discrete energy law established below. Although $\alpha(t)\equiv1$ at the continuous level, its discrete counterpart may deviate from unity due to temporal discretization. Indeed, \eqref{sav:alpha-linear} shows that a smaller value of the parameter $Q>0$ drives the filtered value $\alpha^{n+1}$ toward one. This behavior is confirmed by the numerical tests in Section~\ref{sec:test}. Unless otherwise specified, $Q=10^{-3}$ is adopted in the following experiments.
\end{remark}
 
To initialize the two-step scheme, we compute $\rho^1$, $\bm u^1$, $p^1$, and $\alpha^1$
by the following first-order start-up procedure:
\begin{align}
  & \frac{\overline{\rho}^{1}-\rho^0}{\tau}
    +\nabla\cdot\bigl(\overline{\rho}^{1}\bm u^0\bigr)=0,  \label{startup:rho}\\
  & \rho^1=\mathbb{P}_{\maC}\, \overline{\rho}^{1},  
  \label{startup:projection}\\
  & \sigma^1\frac{\sigma^1\bm u^1-\sigma^0\bm u^0}{\tau}   +N(\rho^1,\bm u^0,\bm u^1)
    -L(\nu^1,\bm u^1,\bm u^0)  +\nabla p^1  =\bm{f}^1,  \label{startup:up}\\
  & \nabla\cdot\bm u^1=0,   \label{startup:div}\\
  & \frac{\alpha^1-\alpha^0}{\tau}  =
   Q\bigl[ -2\bigl(\nu^1\nabla^{\mrs}\bm u^1,\nabla^{\mrs}\bm u^1\bigr)
         -\maL(\nu^1,\bm u^1,\bm u^0;\bm u^1) \bigr],\quad \alpha^0=1.
  \label{startup:alpha}
\end{align}

\begin{remark}\label{rem:mu:1} 
In practical simulations, the dynamic viscosity $\nu$ may depend nonlinearly on the density and other flow variables through a prescribed constitutive relation. For instance, the viscosity depends on the density and the strain-rate magnitude in the form \cite{espinozaM2AN2025}
$$
\nu:=\nu\!\left(\rho, |\nabla^{\mathrm{s}}\bm{u}|\right).
$$
Then, prior to the velocity--pressure step, we evaluate
\[
\nu^{n+1}
=
\nu\bigl(\rho^{n+1},
|\nabla^{\mrs}\bm u^{*,n+1}|\bigr).
\]
 In particular, if the viscosity depends solely on the density (see, e.g., Example \ref{subsec:e2}), then we can simply take $\nu^{n+1}=\nu(\rho^{n+1})$. 
\end{remark}

\subsection{Physics-preserving analysis}  \label{subsec:analysis}
We now establish the density bounds, discrete mass conservation, and an unconditional modified-energy dissipation law of the proposed scheme. For $n\geq1$, define
\begin{equation*}
\maA(v^{n+1}):=\frac{3}{2}v^{n+1}-2v^n+\frac{1}{2}v^{n-1},\qquad \maB(v^{n+1}):=\frac{3}{2}v^{n+1}-v^n+\frac{1}{2}v^{n-1}.
\end{equation*}
Besides, we also introduce a weighted velocity $\bm w^n:=\sigma^n\bm{u}^n$. Then, the momentum equation \eqref{2th:up}  and the density-weighted filter \eqref{2th:TF:u} can be recast as
\begin{align}
  & \sigma^{n+1}\frac{\maA(\bm w^{n+1})}{\tau}+N\bigl(\rho^{n+1},\bm{u}^{*,n+1},\widehat{\bm{u}}^{n+1}\bigr) 
   -\alpha^{*,n+1}L\bigl(\nu^{n+1},\widehat{\bm{u}}^{n+1},\bm{u}^{*,n+1}\bigr) \notag\\
    &\qquad  \qquad  \qquad ~  +\nabla p^{n+1}=\bm{f}^{n+1}, \label{eq:TF:u}\\
    &  \sigma^{n+1}\widehat{\bm{u}}^{n+1}=\frac{3}{2}\bm w^{n+1}-\bm w^n+\frac{1}{2}\bm w^{n-1}=\maB(\bm w^{n+1}). \label{eq:wB} 
\end{align}

The main physics-preserving properties of the scheme are summarized in the following theorem.
\begin{theorem} \label{thm:pme-conservation} 
Let $\{\rho^n,\bm{u}^n,p^n\}_{n\ge 1}$ be solutions of the constant-step prediction--filter--projection scheme \eqref{density:predict}--\eqref{2th:TF:M}. Then the following properties hold. 
\begin{itemize}
    \item[(i)] Bound preservation: If $\rho_{\min} \leq \rho^0(\bm{x}) \leq \rho_{\max}$, then
    \begin{equation*} \label{th:bound}
	  \rho_{\min} \leq \rho^{n} \leq \rho_{\max}, \quad 1 \le n \le N.
     \end{equation*}

   \item[(ii)] Mass conservation: For any $1 \le n \le N$,  
      \begin{equation*} \label{th:mass}
		 \int_{\Omega}  \rho^{n} \, d \bm{x} = \int_{\Omega}  \rho^{0} \, d \bm{x}.
	\end{equation*}	

   \item[(iii)] Unconditional modified-energy dissipation: If $\bm{f}=\bm{0}$, then
    \begin{equation*}\label{th:energy}
        \frac{\widehat{E}_{\mathrm{m}}^{n+1}-\widehat{E}_{\mathrm{m}}^n}{\tau}
        \leq -2 \bigl\|\sqrt{\nu^{n+1}}\nabla^{\mrs}\bm{u}^{n+1}\bigr\|^2,\quad 1 \le n \le N-1,
    \end{equation*}
where $\widehat{E}_{\mathrm{m}}^{n} =\widehat{E}_{\mathrm{m},1}^{n}+\widehat{E}_{\mathrm{m},2}^{n}$ with 
$$\widehat{E}_{\mathrm{m},1}^{n}
    :=\frac14\bigl(\|\bm w^{n}\|^2 + \|2\bm w^{n}-\bm w^{n-1}\|^2 
        + \|\delta\bm w^{n}\|^2\bigr),~\widehat{E}_{\mathrm{m},2}^{n}:=\frac{1}{2Q} \bigl(3\alpha^{n}-\alpha^{n-1}\bigr).
$$        
\end{itemize}
\end{theorem}
\begin{proof} The first two assertions follow directly from the constraints imposed in \eqref{density:projection}. It remains to prove the unconditional modified-energy dissipation law.

Note that for arbitrary $a,b,c \in \mathbb{R}$, the following algebraic identity holds:
\begin{align*}
   & (3a-4b+c)(3a-2b+c)\\
   &=\bigl[a^2+ (2a-b)^2+ (a-b)^2\bigr]-\bigl[b^2+ (2b-c)^2+ (b-c)^2\bigr]  +3 (a-2b+c)^2.
\end{align*}
Therefore, it follows from \eqref{eq:wB} that
\begin{equation}\label{th:au:sigmu}
\begin{aligned}
   \bigl(\maA(\bm w^{n+1}),\sigma^{n+1}\widehat{\bm{u}}^{n+1}\bigr)
   =\bigl(\maA(\bm w^{n+1}),\maB(\bm w^{n+1})\bigr)
   = \widehat{E}_{\mathrm{m},1}^{n+1}-\widehat{E}_{\mathrm{m},1}^{n}
     +\frac34\|\delta^2\bm w^{n}\|^2.
\end{aligned}
\end{equation}
Now, taking the inner product of \eqref{eq:TF:u} with $\tau\widehat{\bm{u}}^{n+1}$, and using \eqref{2th:div}, the boundary condition $\widehat{\bm{u}}^{n+1}\cdot\bm n=0$ on $\Gamma$, together with \eqref{eq0:tri} and \eqref{th:au:sigmu}, we obtain, for $\bm{f}^{n+1}=\bm 0$:
\begin{equation}\label{th:momentum-energy}
\begin{aligned}
   \widehat{E}_{\mathrm{m},1}^{n+1}-\widehat{E}_{\mathrm{m},1}^{n}
     +\frac34\|\delta^2\bm w^{n}\|^2
     =\tau\alpha^{*,n+1}\maL\bigl(\nu^{n+1},\widehat{\bm{u}}^{n+1},\bm{u}^{*,n+1};\widehat{\bm{u}}^{n+1}\bigr).
\end{aligned}
\end{equation}
On the other hand, multiplying \eqref{2th:Mhat} by $\tau/Q$ and eliminating $\widehat\alpha^{n+1}$ through the time filter relation \eqref{2th:TF:M} yield
\begin{equation}\label{th:KM}
\begin{aligned}
   \widehat{E}_{\mathrm{m},2}^{n+1}-\widehat{E}_{\mathrm{m},2}^{n} 
   &=-2\tau \bigl\|\sqrt{\nu^{n+1}}\nabla^{\mrs}\bm{u}^{n+1}\bigr\|^2 -\tau\alpha^{*,n+1} \maL \bigl(\nu^{n+1},\widehat{\bm{u}}^{n+1},\bm{u}^{*,n+1};\widehat{\bm{u}}^{n+1}\bigr).
\end{aligned}
\end{equation}

Therefore, adding \eqref{th:momentum-energy} and \eqref{th:KM} together, the terms involving $\alpha^{*,n+1}$ are canceled, and we obtain
\begin{equation}\label{th:energy:eq}
    \widehat{E}_{\mathrm{m}}^{n+1}-\widehat{E}_{\mathrm{m}}^n 
    +\frac{3}{4}\|\delta^2\bm w^{n}\|^2
    =-2\tau \bigl\|\sqrt{\nu^{n+1}}\nabla^{\mrs}\bm{u}^{n+1}\bigr\|^2,
\end{equation}
which yields the desired energy inequality, a discrete counterpart of Theorem \ref{thm:mEnergy}. \qed
\end{proof}

Unlike the predicted velocity $\widehat{\bm{u}}^{n+1}$, the nonlinear recovery $\bm{u}^{n+1}=\bm w^{n+1}/\sigma^{n+1}$ is not, in general, divergence-free. Nevertheless, the energy identity \eqref{th:energy:eq} yields a quantitative bound on the divergence approximation error $\nabla\cdot (\bm{u}^{n+1}-u(t^{n+1}))=\nabla\cdot\bm{u}^{n+1}$, as stated in the following proposition.
\begin{proposition}\label{prop:divergence-err}
Let $\bm{f}=\bm{0}$. The divergence error of the time-filtered velocity satisfies
\begin{equation*}
  \sum_{n=1}^{m} \|\nabla\cdot\bm{u}^{n+1} \|_{-1}^2
   \leq\frac{\widehat{E}_{\mathrm{m}}^1-\widehat{E}_{\mathrm{m}}^{m+1}}{3\sigma_{\min}^2},\qquad \sigma_{\min}:=\sqrt{\rho_{\min}},
\end{equation*}
where the negative norm is defined as
$
  \|q\|_{-1}
  :=\sup_{0\neq\varphi\in H_0^1(\Omega)}
  \frac{\langle q,\varphi\rangle}{\|\nabla\varphi\|}.
$
\end{proposition}
\begin{proof} We rewrite \eqref{eq:wB} as
$
\widehat{\bm{u}}^{n+1}=\bm{u}^{n+1}+\frac{1}{2\sigma^{n+1}}\delta^2\bm w^{n}.
$
Taking the divergence of both sides and invoking the incompressibility condition $\nabla\cdot\widehat{\bm{u}}^{n+1}=0$, we obtain
\begin{equation*} \label{cor:divu}
\nabla\cdot\bm{u}^{n+1}=-\frac{1}{2}\nabla\cdot\left(\frac{\delta^2\bm w^{n}}{\sigma^{n+1}}\right),
\end{equation*}
from which and the bound $\sigma^{n+1}\geq\sigma_{\min}$, we further have
\begin{equation*}
  \|\nabla\cdot\bm{u}^{n+1}\|_{-1}\leq\frac{1}{2\sigma_{\min}}\|\delta^2\bm w^{n}\|.
\end{equation*}
Recall that the energy identity \eqref{th:energy:eq} implies that
\begin{equation*}
  \frac{3}{4}\sum_{n=1}^{m}\|\delta^2\bm w^{n}\|^2\leq\widehat{E}_{\mathrm{m}}^1-\widehat{E}_{\mathrm{m}}^{m+1}.
\end{equation*}
Then, we have
\begin{equation*}
    \sum_{n=1}^{m} \|\nabla\cdot\bm{u}^{n+1} \|_{-1}^2\leq\frac{1}{4\sigma_{\min}^2}\sum_{n=1}^{m}\|\delta^2\bm w^{n}\|^2\leq\frac{\widehat{E}_{\mathrm{m}}^1-\widehat{E}_{\mathrm{m}}^{m+1}}{3\sigma_{\min}^2},
\end{equation*}
which proves the desired conclusion. \qed
\end{proof}

\begin{remark} The density-weighted filter is motivated by the discrete energy structure.
In \cite{li2023filtered}, the time filter is applied directly to the velocity $\bm u$, and
the stability analysis additionally assumes
\begin{equation}\label{assump:e1}
  \max_{0\leq k\leq N-1} \Bigl\|\frac{\rho^{k+1}-\rho^k}{\tau}\Bigr\|\leq C^*\rho_{\min},
\end{equation}
 with $C^*>0$ independent of $\tau$. 
However, this assumption is restrictive for flows exhibiting rapid temporal variations in density. In contrast, by filtering $\bm w=\sigma\bm{u}$, the energy stability proved in Theorem~\ref{thm:pme-conservation} is obtained without the additional uniform-in-time bound assumption \eqref{assump:e1}. 
To the best of our knowledge, this is the first density-weighted time-filtered scheme proposed for the variable-density, variable-viscosity Navier--Stokes equations. Most importantly, the proposed scheme satisfies an unconditional modified-energy dissipation law without imposing such a stringent boundedness assumption on the discrete density variation.
\end{remark}
\begin{remark}
The backward Euler prediction combined with the classical time filter yields a second-order temporal approximation \cite{decaria2020time,guzel2018time}. Although a complete convergence analysis of the proposed scheme is beyond the scope of this paper, we note that the projection step does not reduce the temporal accuracy.
Assume that
$
\|\widehat{\rho}^{n+1}-\rho(t^{n+1})\|=\maO(\tau^2).$
Since $\maC$ is closed and convex and
$\rho(t^{n+1})\in\maC$, the $L^2$-projection $\mathbb{P}_{\maC}$ is non-expansive;
hence,
\[
\|\rho^{n+1}-\rho(t^{n+1})\|
  = \|\mathbb{P}_{\maC}\widehat{\rho}^{n+1}  -\mathbb{P}_{\maC}\rho(t^{n+1})\|  
  \le \|\widehat{\rho}^{n+1}-\rho(t^{n+1})\|=\maO(\tau^2).
\]
Thus, the projection preserves the second-order temporal accuracy of the
filtered density. Since $\rho\ge\rho_{\min}>0$, it also follows that
\[
\|\sigma^{n+1}-\sigma(t^{n+1})\|=\maO(\tau^2).
\]
If, in addition,
$\|\bm w^{n+1}-\bm w(t^{n+1})\|=\maO(\tau^2)$ and
$\bm u(t^{n+1})\in L^\infty(\Omega)^d$, then
\[
\begin{aligned}
&\|\bm u^{n+1}-\bm u(t^{n+1})\| \\
&\le
\frac{1}{\sqrt{\rho_{\min}}}
\big(
\|\bm w^{n+1}-\bm w(t^{n+1})\|
+
\|\bm u(t^{n+1})\|_{L^\infty(\Omega)}
\|\sigma^{n+1}-\sigma(t^{n+1})\|
\big)
=\maO(\tau^2).
\end{aligned}
\]
Hence, the recovered velocity $\bm u^{n+1}=\bm w^{n+1}/\sigma^{n+1}$ is formally second-order accurate.
\end{remark}

\section{A variable-step prediction--filter--projection scheme}\label{sec:variable-scheme}
We extend the constant-step scheme to the nonuniform grid
\begin{equation*}
  0=t_0<t_1<\cdots<t_N=T,\quad
  \tau_n:=t_n-t_{n-1},\quad
  r_{n}:=\frac{\tau_{n+1}}{\tau_{n}},\quad n\ge 1,
\end{equation*}
and devise an efficient adaptive time-stepping strategy for long-term simulations.

For a sequence $\{v^n\}_{n\geq 0}$, set 
$ \partial_\tau v^n:={(v^n-v^{n-1})}/{\tau_n}$, and introduce the variable-step second-order backward differentiation formula (BDF2) and the second-order extrapolation formula by
\begin{equation*}\label{vs:D2}
  D_2v^n
  :=\frac{1+2r_{n-1}}{1+r_{n-1}}\partial_\tau v^n
  -\frac{r_{n-1}}{1+r_{n-1}}\partial_\tau v^{n-1}, ~ v^{*,n}:=(1+r_{n-1})v^{n-1}-r_{n-1}v^{n-2},  ~ n\geq 2.
\end{equation*}
Given a predicted value $\widehat v^n$, the variable-step time-filtered counterpart of \eqref{density:filter} is defined by \cite{guzel2018time,li2023filtered}:
\begin{equation}\label{vs:timefilter}
  v^n  =\widehat v^n -\frac{r_{n-1}}{1+2r_{n-1}}
  \left[\widehat v^n-(1+r_{n-1})v^{n-1}+r_{n-1}v^{n-2}
  \right],  \qquad n \geq 2.
\end{equation}

For the constant-step scheme, Theorem~\ref{thm:pme-conservation} establishes an unconditional energy stability in terms of the modified kinetic energy $\widehat E_{\mathrm m,1}^{n+1}$, which differs from the ``original” physical kinetic energy $ E^{n+1}:=\frac12\|\sigma^{n+1}\bm u^{n+1}\|^2$. On nonuniform temporal grids, an additional difficulty arises: the temporal term generated by the weighted-velocity filter,
$(D_2\bm w^n,\widehat{\bm w}^{\,n})$, does not generally coincide with $\frac12D_2\|\bm w^n\|^2$. To recover a $D_2$-compatible physical kinetic energy, we introduce an additional SAV $\beta(t)$:
\begin{equation}\label{vs:continuous-beta}
\beta_t=Z\left[-\frac{\mathrm d}{\mathrm dt}E(t)+\beta\bigl((\sigma\bm{u})_t,\sigma\bm{u}\bigr)\right],
  \qquad \beta(0)=\beta_0,
\end{equation}
where $Z>0$ is a user-defined stabilization parameter. A detailed calculation shows that \eqref{vs:continuous-beta} has a unique solution $\beta(t)\equiv 1$ for $\beta_0=1$. Unless otherwise specified, we set $Z=10^{-3}$ in the numerical simulations. The present formulation 
offers a more transparent energy interpretation, as will be discussed later. The original momentum equation can now be recast as the following doubly augmented continuous system:
\begin{equation}\label{vs:continuous-beta-system}
   \beta\sigma(\sigma\bm{u})_t+N(\rho,\bm{u},\bm{u})-\alpha L(\nu,\bm{u},\bm{u})+\nabla p=\bm{f},\quad
    \nabla\cdot\bm{u}=0,
\end{equation}
where $\alpha(t)$ and $\beta(t)$ are defined by \eqref{sav:alpha-linear} and \eqref{vs:continuous-beta}, respectively, with initial values $\alpha_0=\beta_0=1$.
In what follows, starting from the reformulated system \eqref{vs:continuous-beta-system} together with \eqref{sav:eq:rho}, we shall propose a variable-step prediction--filter--projection scheme.

\subsection{Second-order numerical scheme}\label{subsec:variable-scheme}
Let  $\widehat{\bm w}^{n}:=\sigma^n\widehat{\bm u}^{n}$.
The variable-step prediction--filter--projection scheme is proposed as follows:

\medskip
\noindent\textbf{Step 0 (Initialization).}
At the first time level $t=t_1$, compute \{$\rho^1$, $\bm{u}^1$, $p^1$, $\alpha^1$\} from \eqref{startup:rho}--\eqref{startup:alpha} with the time step $\tau_1$. Set $\beta^0=1$, and compute $\beta^1$ from
\begin{equation}\label{vs:beta-startup}
\frac{\beta^1-\beta^0}{Z\tau_1}
=-\frac12\partial_\tau\|\bm w^1\|^2+\beta^0\bigl(\partial_\tau\bm w^1,\bm w^1\bigr)
=\frac{1}{2\tau_1}\|\bm w^1-\bm w^0\|^2.
\end{equation}
 With this initialization, the discrete energy balance law also holds at $n=1$.

\medskip
For $n\geq1$, the variable-step scheme advances the solution from $t_n$ to $t_{n+1}$ as follows.

\medskip
\noindent\textbf{Step 1 (Density prediction).}
Compute the predicted density $\overline{\rho}^{n+1}$ by solving
\begin{equation}\label{vs:density-predict}
\frac{\overline{\rho}^{n+1}-\rho^n}{\tau_{n+1}}+\nabla\cdot\left(\overline{\rho}^{n+1}\bm{u}^{*,n+1}\right)=0.
\end{equation}

\noindent\textbf{Step 2 (Density time filtering).}
Apply the variable-step time filter \eqref{vs:timefilter} to obtain
\begin{equation}\label{vs:density-filter}
\widehat{\rho}^{\,n+1}=\overline{\rho}^{n+1}-\frac{r_n}{1+2r_n}\left[\overline{\rho}^{n+1}-(1+r_n)\rho^n+r_n\rho^{n-1}\right].
\end{equation}

\noindent\textbf{Step 3 (Bound-preserving and mass-conserving projection).}
Apply the same projection \eqref{density:projection} to the variable-step time filtered density $\widehat{\rho}^{\,n+1}$:
\begin{equation}\label{vs:density:projection}
   \rho^{n+1}  = \mathbb{P}_{\maC}\,\widehat{\rho}^{n+1}, 
\end{equation}
where the projected density $\rho^{n+1}$ and the scalar Lagrange multiplier $\xi^{n+1}$ are determined by \eqref{proj-clamp} and \eqref{eq:F-xi-cont}, respectively.

\medskip
\noindent\textbf{Step 4 (Velocity--pressure prediction).} Set $\sigma^{n+1}:=\sqrt{\rho^{n+1}}$. Compute the predicted velocity $\widehat{\bm{u}}^{n+1}$ and pressure ${p}^{n+1}$ by solving 
\begin{align}
 & \beta^{*,n+1}\sigma^{n+1}
\frac{\sigma^{n+1}\widehat{\bm{u}}^{n+1}-\sigma^n\bm{u}^n}{\tau_{n+1}}
+N\bigl(\rho^{n+1},\bm{u}^{*,n+1},\widehat{\bm{u}}^{n+1}\bigr) 
 \nonumber\\
   & \qquad \qquad ~
   -\alpha^{*,n+1}L\bigl(\nu^{n+1},\widehat{\bm{u}}^{n+1},\bm{u}^{*,n+1}\bigr) +\nabla p^{n+1}  = \bm{f}^{n+1},   \label{vs:momentum-predict} \\
  & \nabla\cdot\widehat{\bm{u}}^{n+1}=0. \label{vs:momentum-div}
\end{align}

\medskip
\noindent\textbf{Step 5 (Weighted-velocity time filtering).}
Apply the variable-step time filter to the weighted velocity by
\begin{equation}\label{vs:weighted-velocity-filter}
\sigma^{n+1}\bm{u}^{n+1}=\sigma^{n+1}\widehat{\bm{u}}^{\,n+1}-\frac{r_n}{1+2r_n}\left[\sigma^{n+1}\widehat{\bm{u}}^{\,n+1}-(1+r_n)\sigma^n\bm{u}^n+r_n\sigma^{n-1}\bm{u}^{n-1}\right].
\end{equation}
 
\noindent\textbf{Step 6 (D-CSAV prediction).} Compute the predicted auxiliary variables $\widehat{\alpha}^{\,n+1}$ and $\widehat{\beta}^{\,n+1}$ from
\begin{align}
  & \frac{\widehat{\alpha}^{\,n+1}-\alpha^n}{\tau_{n+1}}
   =Q\bigl[-2\bigl(\nu^{n+1}\nabla^{\mathrm s}\bm{u}^{n+1},\nabla^{\mathrm s}\bm{u}^{n+1}\bigr) -\alpha^{*,n+1}\maL\bigl(\nu^{n+1},\widehat{\bm{u}}^{\,n+1},\bm{u}^{*,n+1};\widehat{\bm{u}}^{\,n+1}\bigr)\bigr], \label{vs:alpha-predict} \\
 & \frac{\widehat{\beta}^{n+1}-\beta^n}{Z\tau_{n+1}}
  =-\frac12D_2\|\bm w^{n+1}\|^2
  +\beta^{*,n+1}\bigl(D_2\bm w^{n+1},\widehat{\bm w}^{n+1}\bigr). \label{vs:beta-predict}
\end{align}

\noindent \textbf{Step 7 (D-CSAV time filtering).}
Apply the variable-step time filter to obtain
\begin{align}
& \alpha^{n+1}=\widehat{\alpha}^{\,n+1}-\frac{r_n}{1+2r_n}\left[\widehat{\alpha}^{\,n+1}-(1+r_n)\alpha^n+r_n\alpha^{n-1}\right], \label{vs:alpha-filter}\\ 
& \beta^{n+1}
  =\widehat{\beta}^{n+1}
  -\frac{r_n}{1+2r_n}
  \left[
    \widehat{\beta}^{n+1}-(1+r_n)\beta^n+r_n\beta^{n-1}
  \right]. \label{vs:beta-filter}
\end{align}
 
\subsection{Physics-preserving analysis and adaptive time-stepping strategy}\label{subsec:variable-analysis}
We first establish the physics-preserving properties of the variable-step
scheme and then present an efficient adaptive time-stepping strategy.

Set $D_2v^1:=\partial_\tau v^1$. For $n\geq2$, \eqref{vs:timefilter} is equivalently written as
\begin{equation}\label{vs:filter-D2}
  \widehat v^n-v^{n-1}=\tau_n D_2v^n,
\end{equation}
and the variable-step BDF2 operator admits the discrete convolution form
\cite{liaoZhang2021,liao2023}
\begin{equation*}
   D_2 v^n=\sum_{j=1}^n d_{n-j}^{(n)}\partial_\tau v^j,
\end{equation*}
where the discrete convolution kernels $\{d_{j}^{(n)}\}$ satisfy $d_0^{(1)}=1$, and for $n \ge 2$,
\begin{equation*}\label{vs:BDF-kernels}
  d_0^{(n)}=\frac{1+2r_{n-1}}{1+r_{n-1}},\quad d_1^{(n)}=-\frac{r_{n-1}}{1+r_{n-1}},\quad d_j^{(n)}=0 \quad \text{for}~ 2\le j \le n-1.
\end{equation*}

It is worth noting that the convolution coefficients \(d_{j}^{(n)}\) depend on the local time-step ratios and vary with the time level. As a result, the resulting variable-step \(D_2\) operator cannot be handled by the same direct telescoping argument as in the constant-step case. To recover the original first-order difference structure suitable for energy analysis, we recursively define the DOC kernels $\{\theta_{n-j}^{(n)}\}$ by 
\begin{equation*}\label{vs:DOC-recursion}
\theta_0^{(n)}:=\frac1{d_0^{(n)}},\quad
\theta_{n-j}^{(n)}:=-\frac1{d_0^{(j)}}\sum_{i=j+1}^n\theta_{n-i}^{(n)}d_{i-j}^{(i)}, \quad  1\leq j\leq n-1.
\end{equation*}
As proved in \cite{liaoZhang2021,liao2023}, the DOC kernels are positive and satisfy the following convolution identity and discrete orthogonality relation:
\begin{equation}\label{vs:DOC-properties}
  \theta_{n-j}^{(n)}>0,~ 1\leq j\leq n;\,
  \sum_{j=1}^n\theta_{n-j}^{(n)}D_2v^j=\partial_\tau v^n;\,
\sum_{j=k}^n \theta_{n-j}^{(n)}d_{j-k}^{(j)}
   =\delta_{n,k}, ~ 1\leq k\leq n,
\end{equation}
for $1\le n \le N$, where $\delta_{n,k}$ is the Kronecker delta. 

The following theorem summarizes the main physics-preserving properties of the variable-step scheme.
\begin{theorem}\label{thm:variable-energy}
Let $\{\rho^n,\bm{u}^n,p^n\}_{n\geq1}$ be solutions of the variable-step prediction--filter--projection scheme \eqref{vs:beta-startup}--\eqref{vs:beta-filter}. Then the following properties hold.
\begin{itemize}
    \item[(i)] Bound preservation: If $\rho_{\min} \leq \rho^0(\bm{x}) \leq \rho_{\max}$, then
    \begin{equation*} 
	  \rho_{\min} \leq \rho^{n} \leq \rho_{\max}, \quad 1 \le n \le N.
     \end{equation*}

   \item[(ii)] Mass conservation: For any $1 \le n \le N$,  
      \begin{equation*} 
		 \int_{\Omega}  \rho^{n} \, d \bm{x} = \int_{\Omega}  \rho^{0} \, d \bm{x}.
	\end{equation*}	

     \item[(iii)] Unconditional modified-energy dissipation: If $\bm{f}=\bm{0}$, then
  \begin{equation*}
    \frac{\overline{E}_{\mathrm m}^n-\overline{E}_{\mathrm m}^{n-1}}{\tau_n}
    = -2\sum_{j=1}^n \theta_{n-j}^{(n)} \big\|\sqrt{\nu^j}\nabla^{\mrs}\bm{u}^j\big\|^2  \leq 0,
    \quad 1\leq n\leq N,
  \end{equation*}
  where $\overline{E}_{\mathrm m}^n  =\overline{E}_{\mathrm m,1}^n
    +\overline{E}_{\mathrm m,2}^n  +\overline{E}_{\mathrm m,3}^n$ with
  \begin{equation*}
    \overline{E}_{\mathrm m,1}^n   :=\frac12\|\bm w^n\|^2,
    \quad
    \overline{E}_{\mathrm m,2}^n  :=\frac{\alpha^n-1}{Q},
    \quad
    \overline{E}_{\mathrm m,3}^n    :=\frac{\beta^n-1}{Z}.
  \end{equation*}
\end{itemize}
\end{theorem}
\begin{proof}  The first two assertions follow immediately from the definition of $\maC$ (see \eqref{adm:set}) in the $L^2$-projection step \eqref{vs:density:projection}. Next, we show the last assertion.

For $2\leq n\leq N$, taking the inner product of \eqref{vs:momentum-predict} with $\widehat{\bm{u}}^n$, and proceeding as in the proof of Theorem~\ref{thm:pme-conservation}, with the aid of the filter identity \eqref{vs:filter-D2}, we obtain 
\begin{equation}\label{vs:momentum-energy}
  \beta^{*,n} \bigl(D_2\bm w^n,\widehat{\bm w}^n \bigr)
   -\alpha^{*,n} \maL \bigl(\nu^n,\widehat{\bm{u}}^n,\bm{u}^{*,n}; \widehat{\bm{u}}^n \bigr) 
  =\bigl(\bm{f}^n,\widehat{\bm{u}}^n\bigr).
\end{equation}
Furthermore, applying \eqref{vs:filter-D2} to the D-CSAV prediction equations \eqref{vs:alpha-predict} and \eqref{vs:beta-predict} yields
\begin{equation}\label{vs:energy}
  \begin{aligned}
  & -\alpha^{*,n} \maL\bigl(\nu^n,\widehat{\bm u}^n,\bm{u}^{*,n}; 
     \widehat{\bm u}^n\bigr)
   = \frac{1}{Q}D_2\alpha^n 
     + 2\, \bigl\|\sqrt{\nu^n}\nabla^{\mathrm{s}}\bm{u}^n\bigr\|^2.\\
  & \beta^{*,n}\bigl(D_2\bm w^n,\widehat{\bm w}^n\bigr)
   = \frac{1}{2}D_2\|\bm w^n\|^2 + \frac{1}{Z}D_2\beta^n.
   \end{aligned}
\end{equation}
Now, substituting \eqref{vs:energy} into the left-hand side of \eqref{vs:momentum-energy}, and applying the definition of $\overline E_{\mathrm m}^n$, yields
\begin{equation}\label{vs:energy-balance}
  D_2\overline E_{\mathrm m}^n + 2\, \bigl\|\sqrt{\nu^n}\nabla^{\mrs}\bm{u}^n\bigr\|^2
  =\bigl(\bm{f}^n,\widehat{\bm u}^n\bigr),\quad 2\leq n\leq N.
\end{equation}
In fact, \eqref{vs:energy-balance} also holds at the initial time level, i.e., for $n=1$. Setting $\widehat{\bm u}^1:=\bm{u}^1$ and recalling that $D_2v^1=\partial_\tau v^1$ and $\alpha^0=\beta^0=1$, the update formulas \eqref{startup:alpha} and \eqref{vs:beta-startup} imply
\begin{equation*}
  -\maL\bigl(\nu^1,\bm{u}^1,\bm{u}^0;\bm{u}^1\bigr) + \bigl(\partial_\tau\bm w^1,\bm w^1\bigr)
  =\frac{1}{Q}\partial_\tau\alpha^1+2\bigl\|\sqrt{\nu^1}\nabla^{\mrs}\bm{u}^1\bigr\|^2 + \frac12\partial_\tau\|\bm w^1\|^2+\frac{1}{Z}\partial_\tau\beta^1.
\end{equation*}
On the other hand, taking the inner product of \eqref{startup:up} with $\bm{u}^1$, and eliminating the convective and pressure terms as before, gives
\begin{equation*}
   \bigl(\partial_\tau\bm w^1,\bm w^1\bigr)
     - \maL\bigl(\nu^1,\bm{u}^1,\bm{u}^0;\bm{u}^1\bigr)
   =\bigl(\bm{f}^1,\bm{u}^1\bigr).
\end{equation*}
Therefore, combining the two equations above yields 
\eqref{vs:energy-balance} for all $1\leq n\leq N$.

Now assume that $\bm{f}=\bm{0}$. Multiplying the
level-$j$ form of \eqref{vs:energy-balance} by $\theta_{n-j}^{(n)}$, summing over $j=1,\ldots,n$, and then applying the convolution identity \eqref{vs:DOC-properties}, we obtain
\[
   \partial_\tau\overline E_{\mathrm m}^n
     +  2\sum_{j=1}^n \theta_{n-j}^{(n)}
     \bigl\|\sqrt{\nu^j}\nabla^{\mathrm{s}}\bm{u}^j\bigr\|^2
   = 0.
\]
Since $\theta_{n-j}^{(n)}>0$, the modified-energy dissipation law immediately follows. \qed
\end{proof}

\begin{remark} \label{rem:vart:uni}
Unlike the filter-dependent modified kinetic energy established in Theorem~\ref{thm:pme-conservation}, the variable-step modified energy contains the physical kinetic energy $\frac12\|\bm w^n\|^2$ exactly.
This is enabled by the auxiliary variable $\beta$ together with the $D_2$-consistent discretization in \eqref{vs:beta-predict}, where the energy derivative in \eqref{vs:continuous-beta} is approximated by
$D_2E^{n+1}$ rather than $\partial_\tau E^{n+1}$. Consequently, the physical kinetic energy and the auxiliary-variable corrections are governed by the same $D_2$ operator. If $\partial_\tau E^{n+1}$ were used instead, rewriting the resulting energy relation in a unified $D_2$ form would introduce the residual
\[
  D_2E^{n+1}-\partial_\tau E^{n+1}
  =
  \frac{r_n}{1+r_n}
  \left(
    \partial_\tau E^{n+1}
    -\partial_\tau E^n
  \right),
\]
which has no definite sign and prevents the DOC identity from being applied directly to the modified energy as a whole. Thus, the $D_2$-consistent discretization in \eqref{vs:beta-predict} is essential to the unconditional modified-energy dissipation law, which is closer to the ``original” one than that of the constant-step scheme. Finally, since $\alpha=\beta\equiv1$ at the continuous level, the terms
$(\alpha^n-1)/Q$ and $(\beta^n-1)/Z$ measure the deviations of the discrete auxiliary variables from their continuous values and may be further reduced by relaxation techniques \cite{jiang2022improving,liu2024novel,zhang2024novel}.
\end{remark}

To further improve computational efficiency, particularly for long-term flow simulations, an efficient and reliable adaptive time-stepping strategy should be developed for the proposed variable-step scheme. Such a strategy would take smaller steps during rapid transients and larger steps when the designed the energy-based indicator varies slowly, thereby reducing computational cost while maintaining temporal resolution. 

Given the underlying energy structure of the variable-density, variable-viscosity incompressible flow, we propose an energy-based time-step selection strategy. Specifically, for prescribed minimum and maximum time steps $\tau_{\min}$ and $\tau_{\max}$, we define the energy-based candidate time step as
\begin{equation}\label{vs:adaptive-indicator}
\tau_{\mathrm{E}}^{n+1} 
 :=\max\left\{
    \tau_{\min}, \frac{\tau_{\max}}{\sqrt{1+\kappa (C_{\mathrm{E}}^n)^2}}
 \right\}, 
\end{equation}
where $\kappa\geq 0$ controls the sensitivity of the time-step adjustment to the energy variation, and the associated energy-variation indicator $C_{\mathrm{E}}^n$ is defined by   
\begin{equation*}
    C_{\mathrm{E}}^n:=\Big|
    \frac{\overline{E}_{\mathrm m}^n-\overline{E}_{\mathrm m}^{n-1}}{\tau_n}
    \Big|.
\end{equation*}
Thus, larger energy variations shall yield smaller candidate time steps, while
$\tau_{\mathrm E}^{n+1}$ approaches $\tau_{\max}$ as $C_{\mathrm E}^n$ decreases.

To ensure a relatively smooth change in the time step, we restrict the next time step $\tau_{n+1}$ by imposing the step-ratio
condition $\frac12 \leq r_n:=\tau_{n+1}/\tau_n\leq2$ \cite{liao2023} to \eqref{vs:adaptive-indicator}:
\begin{equation}\label{vs:adaptive-update}
 \tau_{n+1}
  :=\min\left\{2\tau_n,\,
    \max\left\{\frac12\tau_n,\,\tau_{\mathrm{E}}^{n+1}
    \right\}  \right\}.
\end{equation}
It is clear that the update \eqref{vs:adaptive-update} guarantees $\tau_{\min}\leq \tau_{n+1} \leq\tau_{\max}$ for $1\leq n\leq N$,
provided that $\tau_{\min}\leq\tau_1\leq\tau_{\max}$.

The resulting energy-based adaptive time-stepping strategy is summarized in Algorithm~\ref{alg:adaptive_step}.
\begin{algorithm}[H]\caption{Adaptive time-stepping strategy}
\label{alg:adaptive_step}
\begin{algorithmic}[1]
\Require Numerical solutions up to $t_n$; parameters $\tau_{\min}$, $\tau_{\max}$, and $\kappa$.
\Ensure The next time step $\tau_{n+1}$ and corresponding numerical solutions.

\Statex \textbf{Step 1: Indicator evaluation and candidate time step selection.}
\State Compute $C_{\mathrm{E}}^{n}$ and $\tau_{\mathrm{E}}^{n+1}$ from \eqref{vs:adaptive-indicator}.

\Statex \textbf{Step 2:  Step restriction.}
\State Compute the next time step $\tau_{n+1}$ from \eqref{vs:adaptive-update}.

\Statex \textbf{Step 3: Solution update.}
\State Advance the variable-step scheme from $t_n$ to $t_{n+1}$ using \eqref{vs:density-predict}--\eqref{vs:beta-filter} with the determined time step $\tau_{n+1}$ and the adjacent-step ratio $r_n={\tau_{n+1}}/{\tau_{n}}$.
\end{algorithmic}
\end{algorithm} 
 
\section{Numerical experiments} \label{sec:test}
In this section, we present several numerical experiments to examine the accuracy and physics-preserving properties of the proposed schemes. The density is approximated by the continuous $P_3$ elements, while the velocity and pressure are discretized using the stable Taylor--Hood $P_3$--$P_2$ pair. All computations are performed using FreeFEM++ \cite{hecht2012}. 

For spatial discretization, unlike the standard approach that imposes bounds at the Lagrange interpolation nodes, the bound-preserving $L^2$-projection is performed in the Bernstein–B\'{e}zier coefficient space \cite{kuzmin2020}. This guarantees pointwise density bounds of the finite element solution, owing to the non-negativity and partition-of-unity properties of the Bernstein–B\'{e}zier basis. The numerical errors for $\zeta=\{\bm u,p,\rho\}$, defined as $e_\zeta:=\zeta^N-\zeta(T)$, and the divergence error $e_{\mathrm{div}}  := \nabla\cdot\bm{u}^N$, are computed at the final time $T$. Specifically, the velocity error is measured in the $L^2$- and $H^1$-norms, the pressure and density errors in the $L^2$-norm, and the divergence error in the $L^2$- and $H^{-1}$-norms. The CSAV errors are defined as  $|e_\alpha| := \max_{1\leq n\leq N}|\alpha^n-1|$ and $|e_\beta| := \max_{1\leq n\leq N}|\beta^n-1|$.

\begin{example}\label{num:e1}
In this example, we evaluate both the temporal convergence and the preservation of physical properties of the proposed prediction--filter--projection framework, using both constant and variable time steps. The manufactured solutions, adapted from \cite{liJCP2024}, are defined on the unit disk $\Omega=\{(x,y):x^2+y^2\leq1\}$ and over the time interval $[0, 0.4]$ as follows:
\begin{equation*}
  \begin{aligned}
    & \rho(x,y,t)  = 2 + 0.2\left[  
       x\cos\!\left((1-x^2-y^2)\sin t\right)  +y\sin\!\left((1-x^2-y^2)\sin t\right) \right], \\
    & \bm{u}(x,y,t)  = \left(   -y(1-x^2-y^2)\cos t,\,    x(1-x^2-y^2)\cos t  \right)^{\top}, \\
    & p(x,y,t)   =  \sin x\sin y\sin t.
  \end{aligned}
\end{equation*}
The variable viscosity is prescribed as $\nu(x,y,t)=0.1+(x^2+y^2)(1+\sin t)$.
\end{example}

We first test the temporal accuracy of the proposed constant-step scheme. To isolate temporal discretization errors, we use a fixed spatial grid with 497,011 degrees of freedom (DoFs) for the $P_3$ element and 221,161 DoFs for the $P_2$ element, so that the spatial errors are negligible compared with the temporal errors. As shown in Table~\ref{tab:rate_variable}, the velocity, pressure, density, and divergence errors all exhibit second-order temporal convergence in their respective norms. Moreover, the CSAV error $|e_\alpha|$ remains close to zero with high precision and also achieves second-order accuracy under temporal refinement. 
\begin{table}
[pos=!h]\footnotesize
  \centering
  \caption{Temporal errors and convergence rates under constant time steps.} \label{tab:rate_variable}
  \begin{tabular*}{0.9\columnwidth}{@{\extracolsep\fill}ccccc@{\extracolsep\fill}}
    \toprule
   $N$   & $10$ & $20$ & $40$ & $80$ \\
    \midrule
      $\|e_{u}\|$
        & 1.5790e-4 & 3.8633e-5 & 9.6147e-6 & 2.4007e-6 \\ 
      Rate
        & ---        & 2.03       & 2.01       & 2.00    \\[2pt]
      $\|e_{u}\|_{1}$
        & 1.1211e-3 & 2.7074e-4 & 6.6861e-5 & 1.6628e-5 \\ 
      Rate
        & ---        & 2.05       & 2.02       & 2.01     \\[2pt]
      $\|e_{p}\|$
        & 1.2128e-4 & 3.1667e-5 & 8.0714e-6 & 2.0365e-6 \\ 
      Rate
        & ---        & 1.94       & 1.97       & 1.99     \\[2pt]
      $\|e_{\rho}\|$
        & 5.3849e-5 & 1.3535e-5 & 3.3903e-6 & 8.4828e-7 \\
      Rate
        & ---        & 1.99       & 2.00       & 2.00       \\[2pt]
      $\|e_{\mathrm{div}}\|$
        & 1.0339e-5 & 2.6916e-6 & 6.8591e-7 & 1.7310e-7 \\
      Rate
        & ---        & 1.94       & 1.97       & 1.99       \\[2pt]
      $\|e_{\mathrm{div}}\|_{-1}$
        & 2.6559e-6 & 6.9211e-7 & 1.7645e-7 & 4.4539e-8 \\
      Rate
        & ---        & 1.94       & 1.97       & 1.99       \\[2pt]
      $|e_\alpha|$
        & 1.5823e-6 & 4.4264e-7 & 1.1654e-7 & 2.9876e-8 \\
      Rate
        & ---        & 1.84       & 1.93       & 1.96 \\       
    \bottomrule
  \end{tabular*}
\end{table}

\begin{figure} [pos=!h]
    \centering
        \includegraphics[width=0.32\textwidth,keepaspectratio]{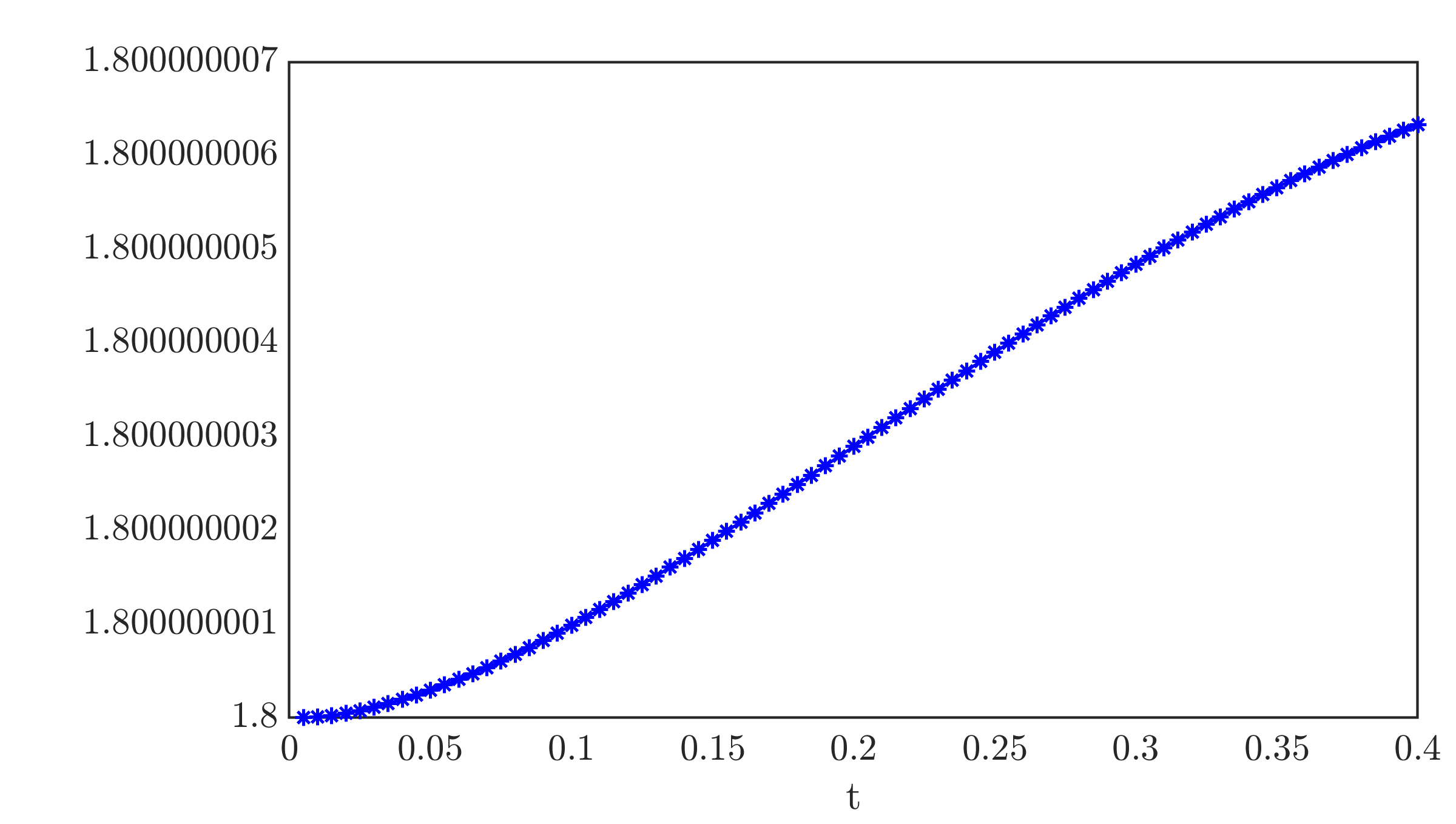}
    \hspace{-0.01\textwidth}%
        \includegraphics[width=0.32\textwidth,keepaspectratio]{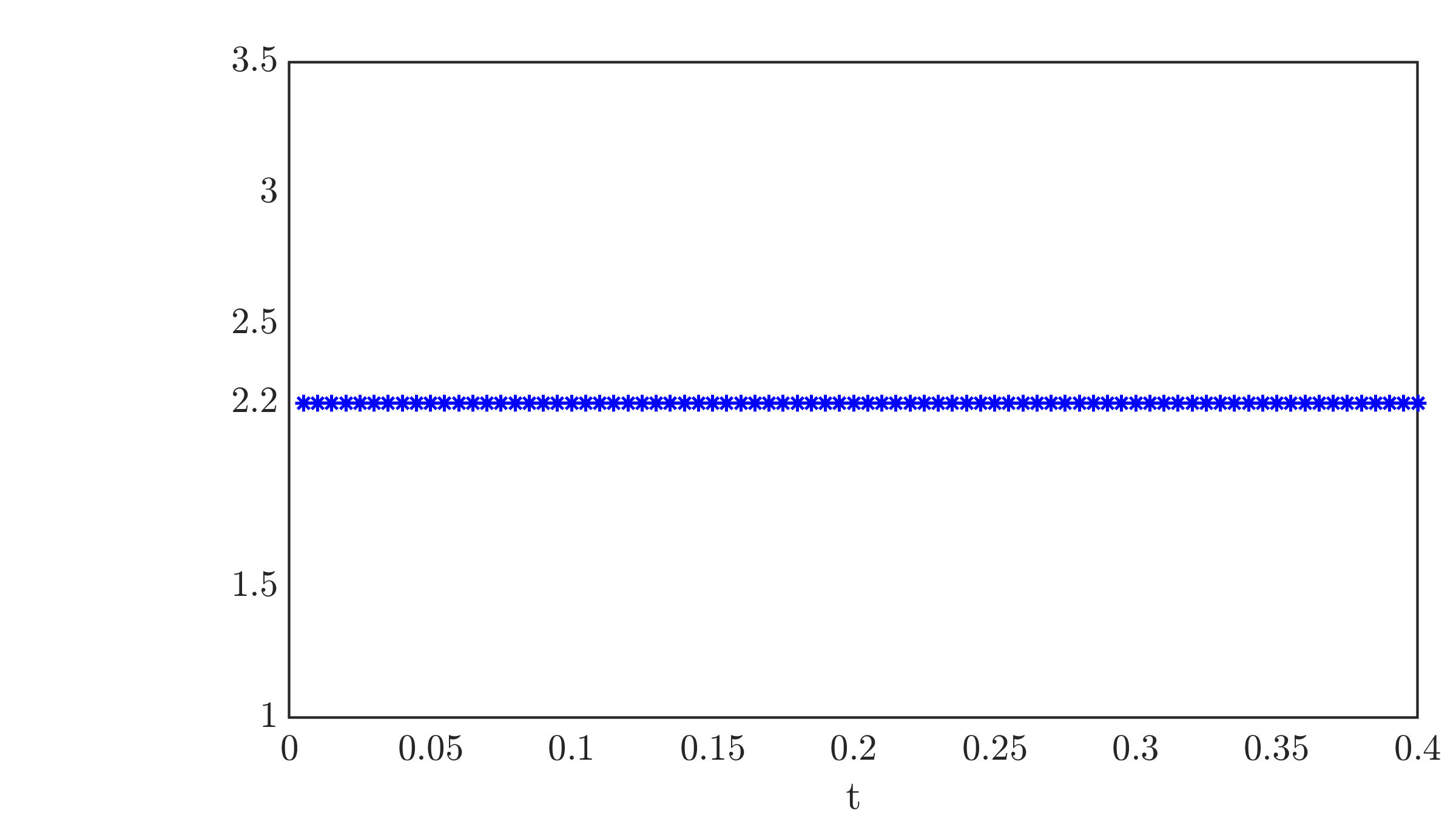}
    \hspace{-0.01\textwidth}%
        \includegraphics[width=0.32\textwidth,keepaspectratio]{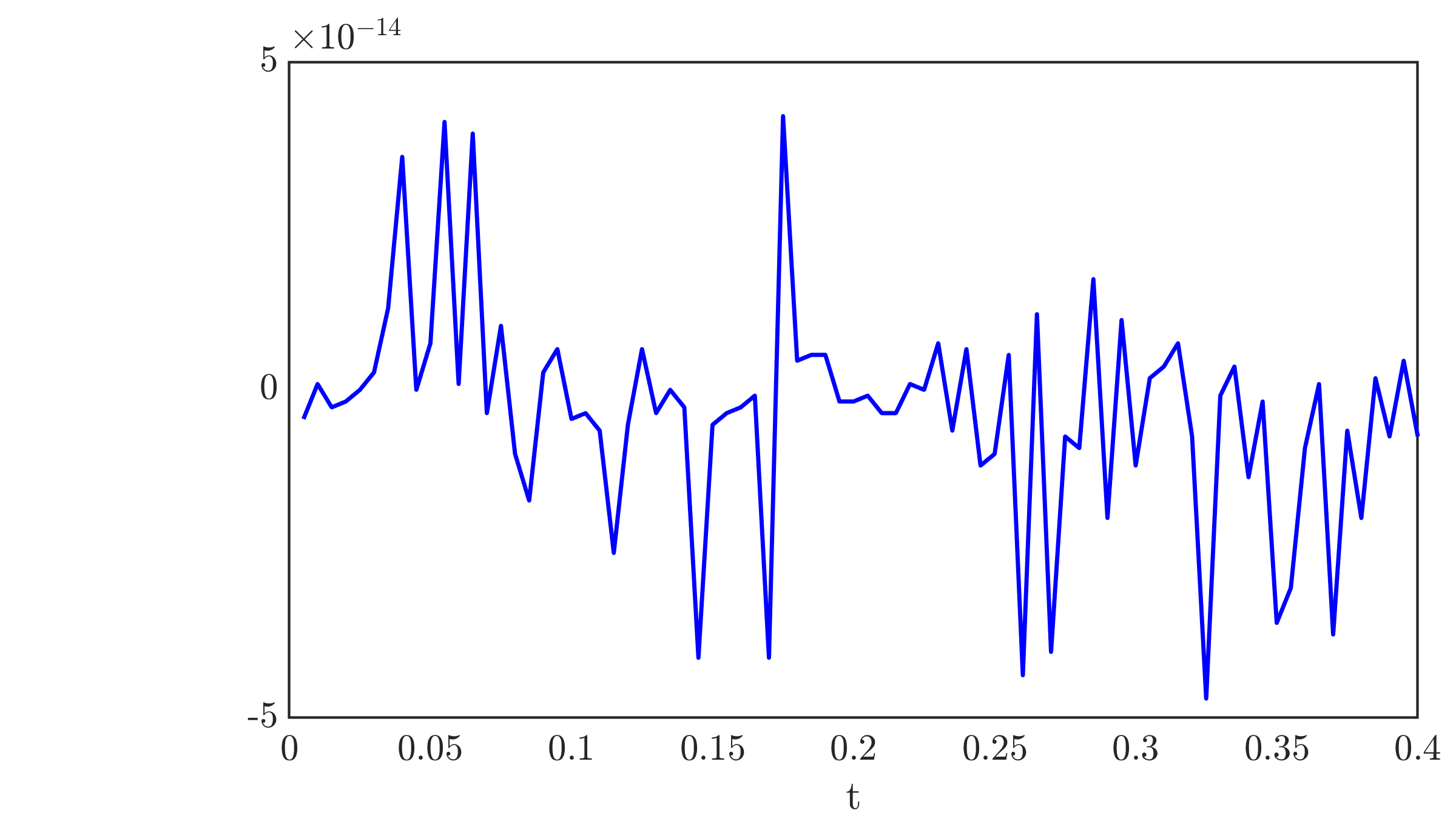}
    \setlength{\abovecaptionskip}{0.0cm} 
	  \setlength{\belowcaptionskip}{0.0cm}    
    \caption{Evolution of density bounds (min: left; max: middle) and mass error (right) under constant time steps.}      \label{fig:physics:var}
\end{figure}

\begin{figure}[pos=!h] 
    \centering
        \includegraphics[width=0.32\textwidth,keepaspectratio]{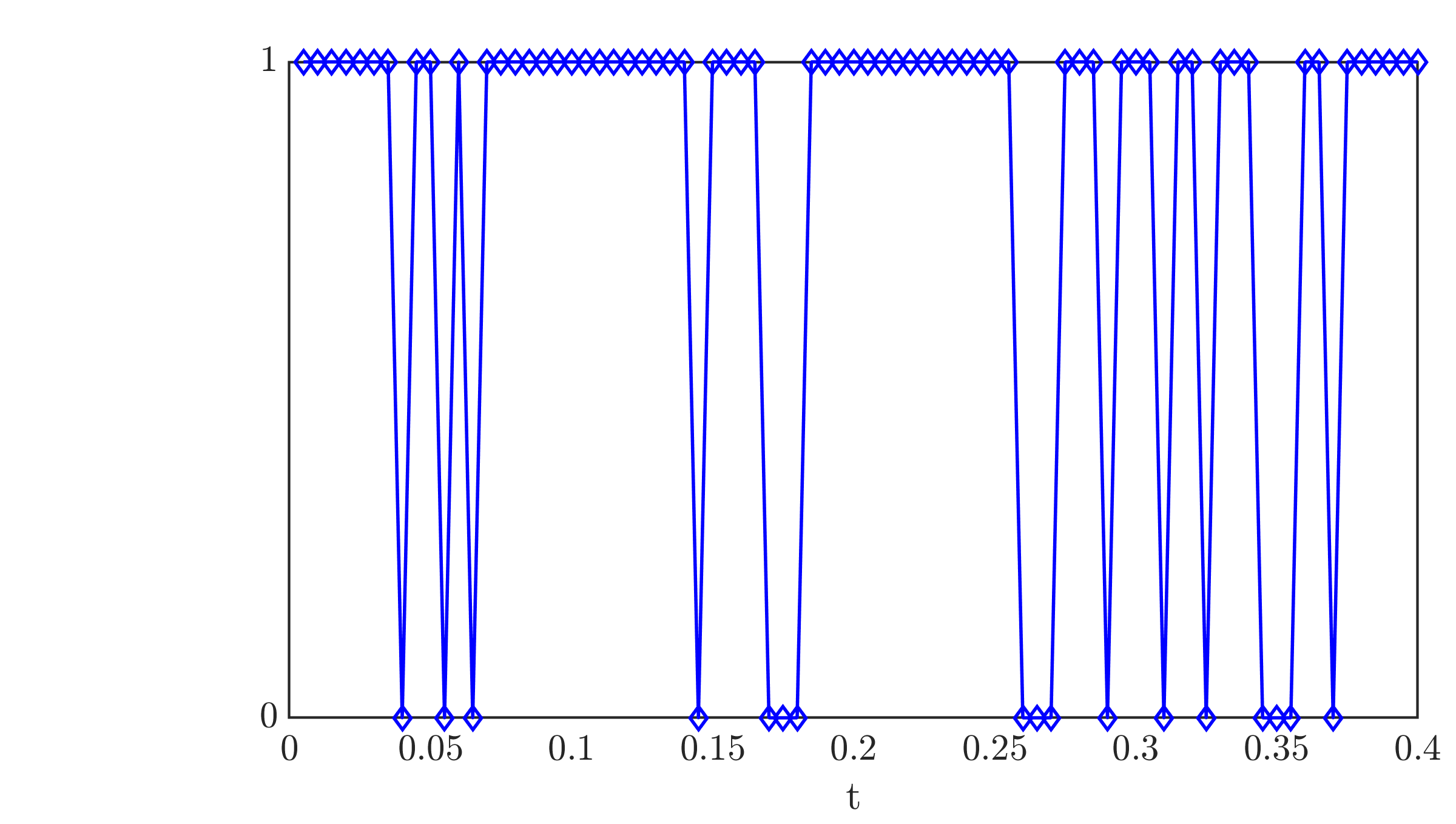}
    \hspace{-0.01\textwidth}%
        \includegraphics[width=0.32\textwidth,keepaspectratio]{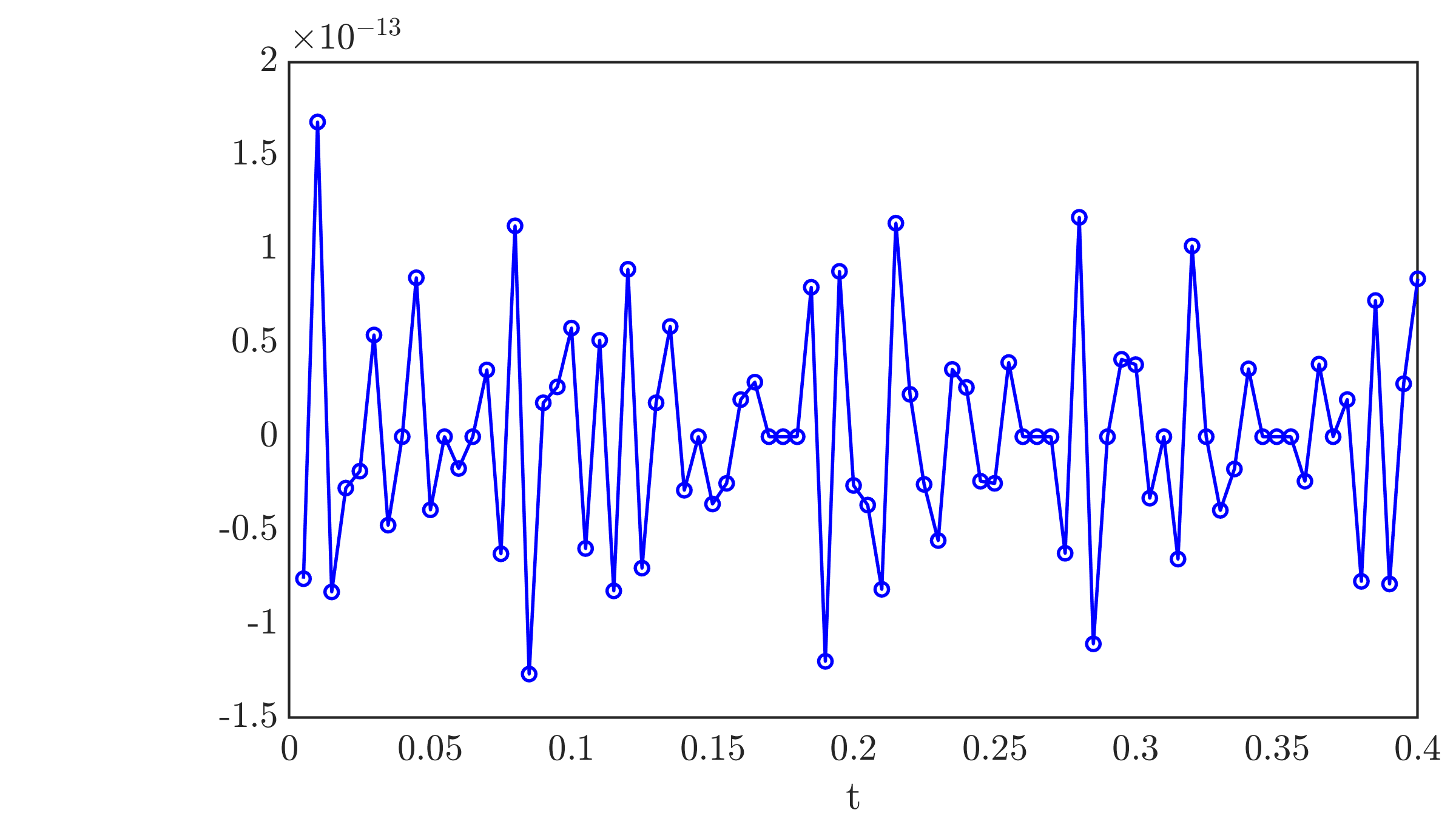}
    \hspace{-0.01\textwidth}%
        \includegraphics[width=0.32\textwidth,keepaspectratio]{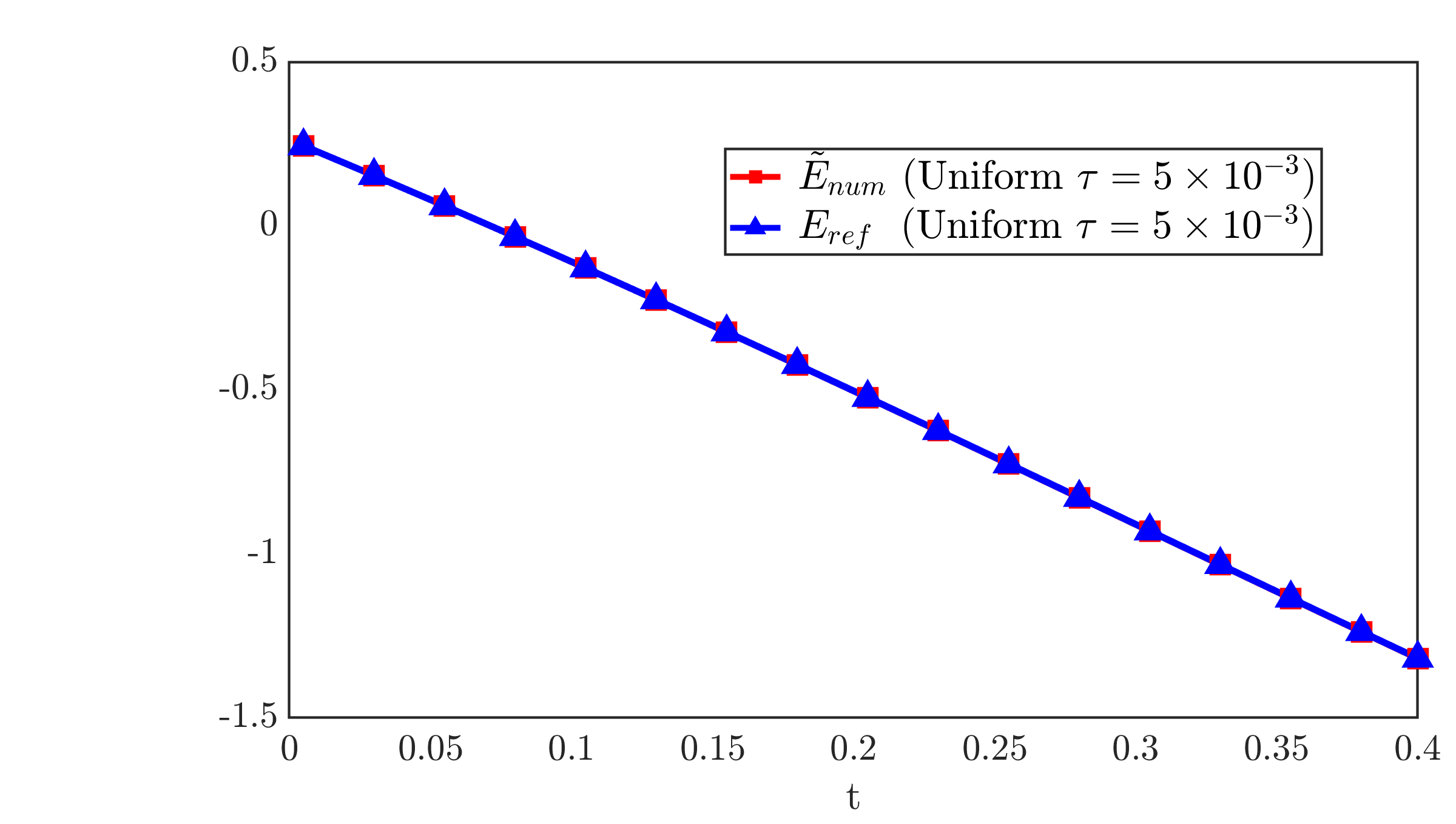}
    \setlength{\abovecaptionskip}{0.0cm} 
	  \setlength{\belowcaptionskip}{0.0cm}    
    \caption{Evolution of the nonlinear iteration count (left), the Lagrange multiplier \(\xi\) (middle), and energy (right) under constant time steps.} \label{fig:Ite_xi:var}
\end{figure}

Next, we further examine the physics-preserving properties using $N=80$ as a representative case.  Fig.~\ref{fig:physics:var} shows that the minimum and maximum density values remain within the prescribed initial bounds $[\rho_{\min}, \rho_{\max}]=[1.8, 2.2]$, and the mass error, defined by $|\int_{\Omega}  \rho^{n} \, d \bm{x}- \int_{\Omega}  \rho^{0} \, d \bm{x}|$, stays on the order of $10^{-14}$, which are in agreement with Theorem~\ref{thm:pme-conservation}.  Since $\bm f\neq\bm 0$ in this test, there is a cumulative energy induced by the external force defined by $W^{n+1}:=\tau\sum_{j=1}^{n+1}(\bm f^j,\widehat{\bm u}^j)$, where $\widehat{\bm u}^{n+1}$ is given by \eqref{2th:up}--\eqref{2th:div}. In this sense, the corresponding discrete energy established in Theorem~\ref{thm:pme-conservation} is modified to
\[
   \widetilde E_{\mathrm{num}}^{n+1}
    := \widehat E_{\mathrm{m},1}^{n+1}
       + \frac{3\alpha^{n+1}-\alpha^n}{2Q} -W^{n+1}
       -\frac{1}{Q},
\]
and it also satisfies the corresponding modified-energy dissipation law
\begin{equation} \label{Enum:W}
   \frac{\widetilde E_{\mathrm{num}}^{n+1} 
      -\widetilde E_{\mathrm{num}}^{n}}{\tau}
   \leq -2 \bigl \|\sqrt{\nu^{n+1}}\nabla^s\bm{u}^{n+1}\bigr\|^2  \leq 0.
\end{equation}
To assess the discrete energy behavior, we compare this modified energy $\widetilde E_{\mathrm{num}}^{n+1}$ with the corresponding reference energy $E_{\mathrm{ref}}^{n+1}:=\frac12\|\bm w^{n+1}\|^2-W^{n+1}$. As shown in Fig.~\ref{fig:Ite_xi:var}, the two energy curves  decrease monotonically and nearly overlap, which further supports that the modified energy is very close to the original one. Moreover, it shows that at each time level, the nonlinear iteration count for solving the nonlinear equation \eqref{eq:F-xi-cont}  does not exceed one, and the Lagrange multiplier $\xi$, used to enforce the mass constraint \eqref{2th:mass}, remains of order $10^{-13}$. This implies that $\rho^{n+1}$ obtained from \eqref{proj-clamp} can be viewed as a certain cut-off of the time-filtered solution $\widehat{\rho}^{n+1}$, while mass conservation is still preserved. 

\begin{figure}
    \centering
     \hspace{-0.01\textwidth}%
      \includegraphics[width=0.32\textwidth,keepaspectratio]{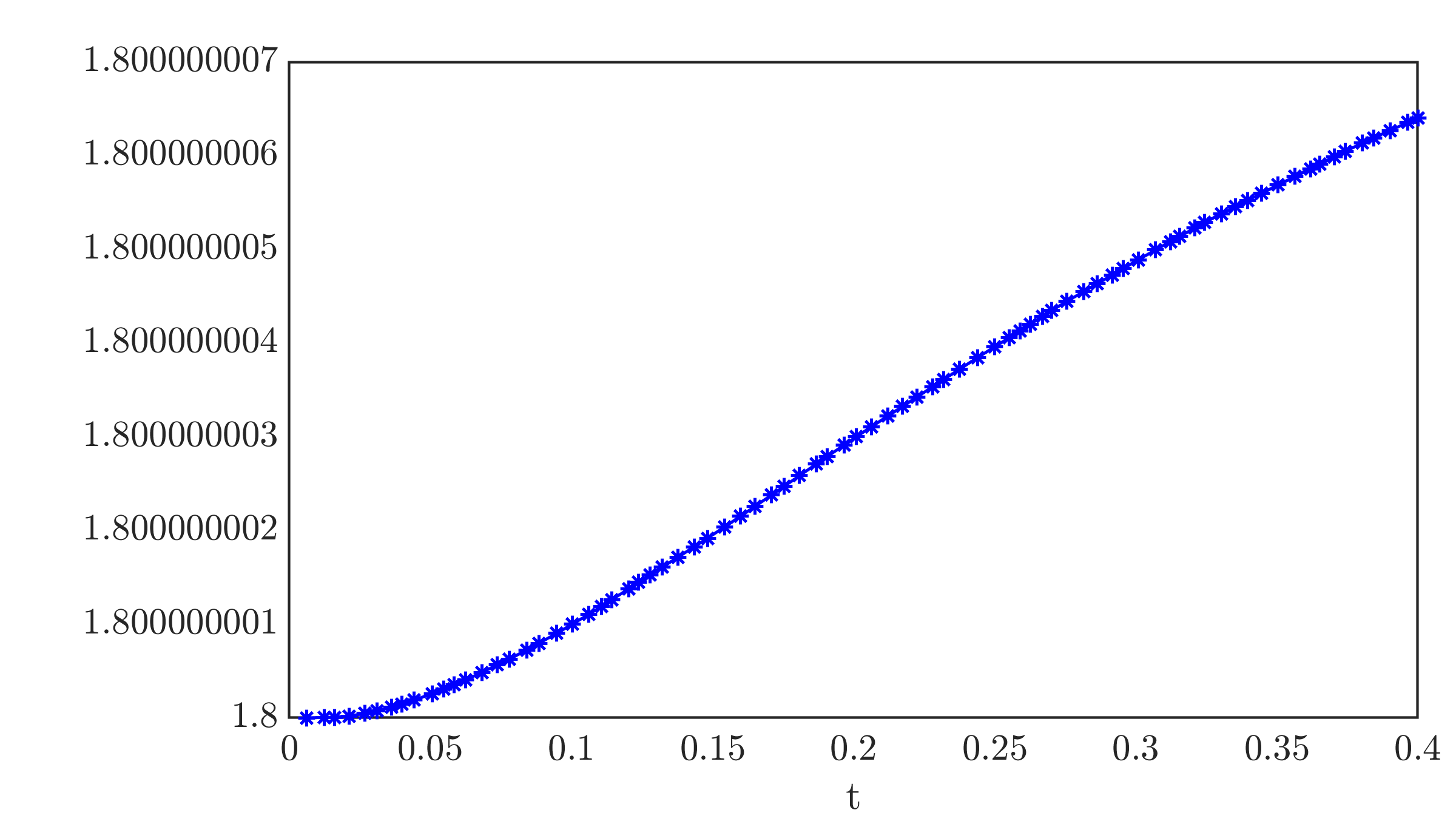}
    \hspace{0.03\textwidth}%
        \includegraphics[width=0.32\textwidth,keepaspectratio]{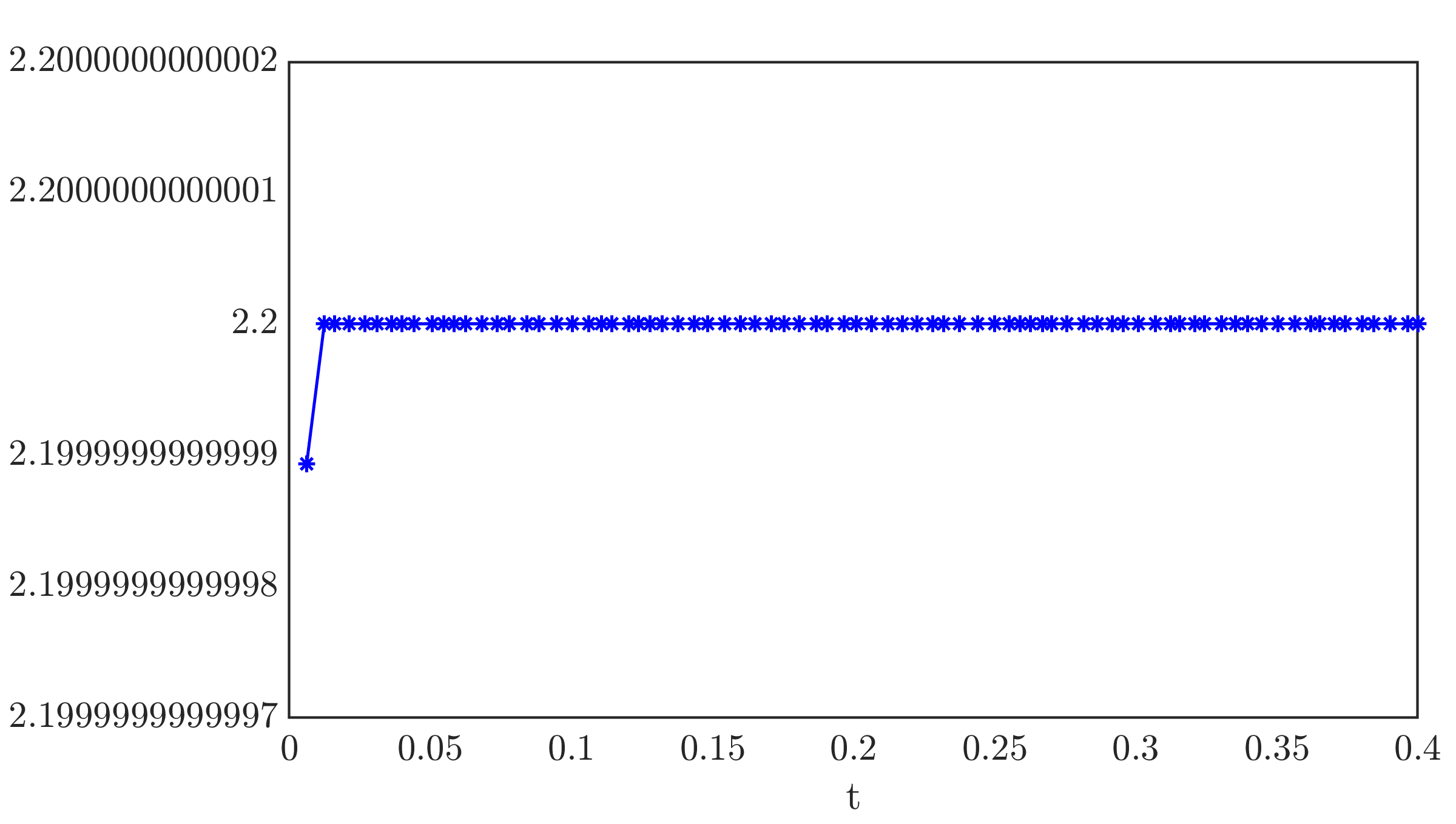}
    \hspace{-0.01\textwidth}%
        \includegraphics[width=0.32\textwidth,keepaspectratio]{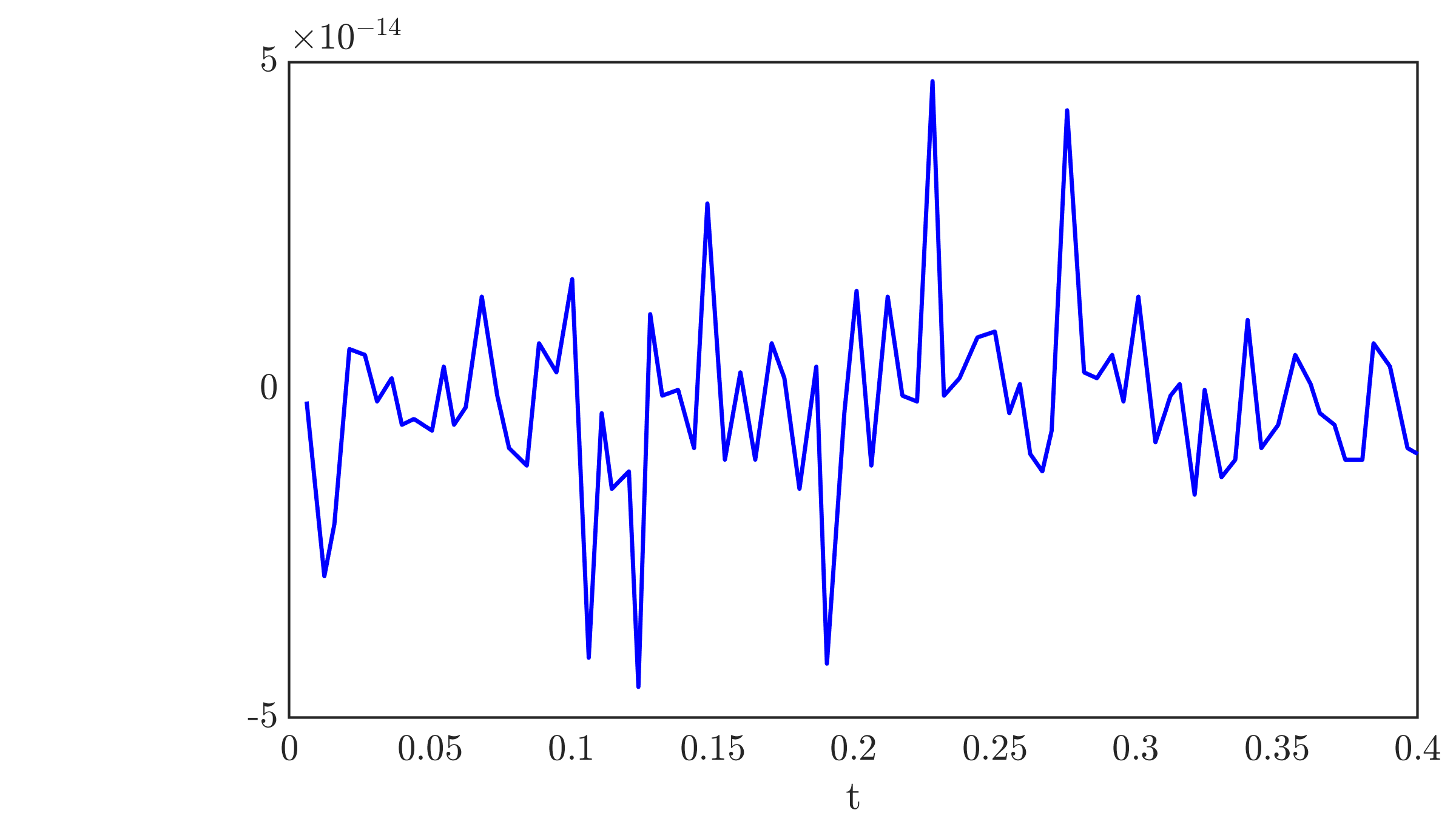}
    \setlength{\abovecaptionskip}{0.0cm} 
	  \setlength{\belowcaptionskip}{0.0cm}    
    \caption{Evolution of density bounds (min: left; max: middle) and mass error (right) under variable time steps.}    
    \label{fig:physics:varynugradedstep}
\end{figure}

\begin{figure}
    \centering

    \hspace{-0.01\textwidth}%
        \includegraphics[width=0.32\textwidth,keepaspectratio]{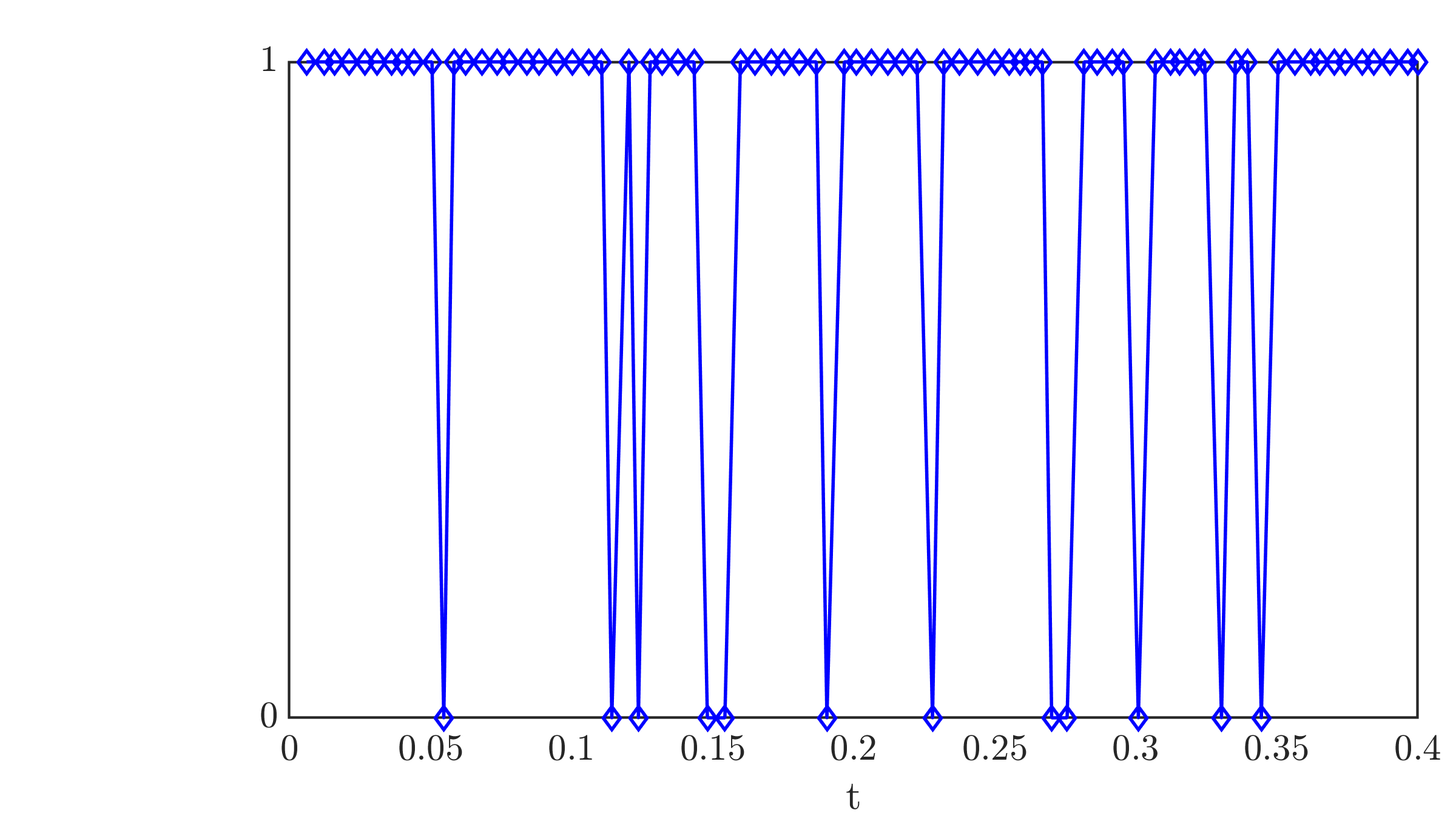}
    \hspace{0.03\textwidth}%
        \includegraphics[width=0.32\textwidth,keepaspectratio]{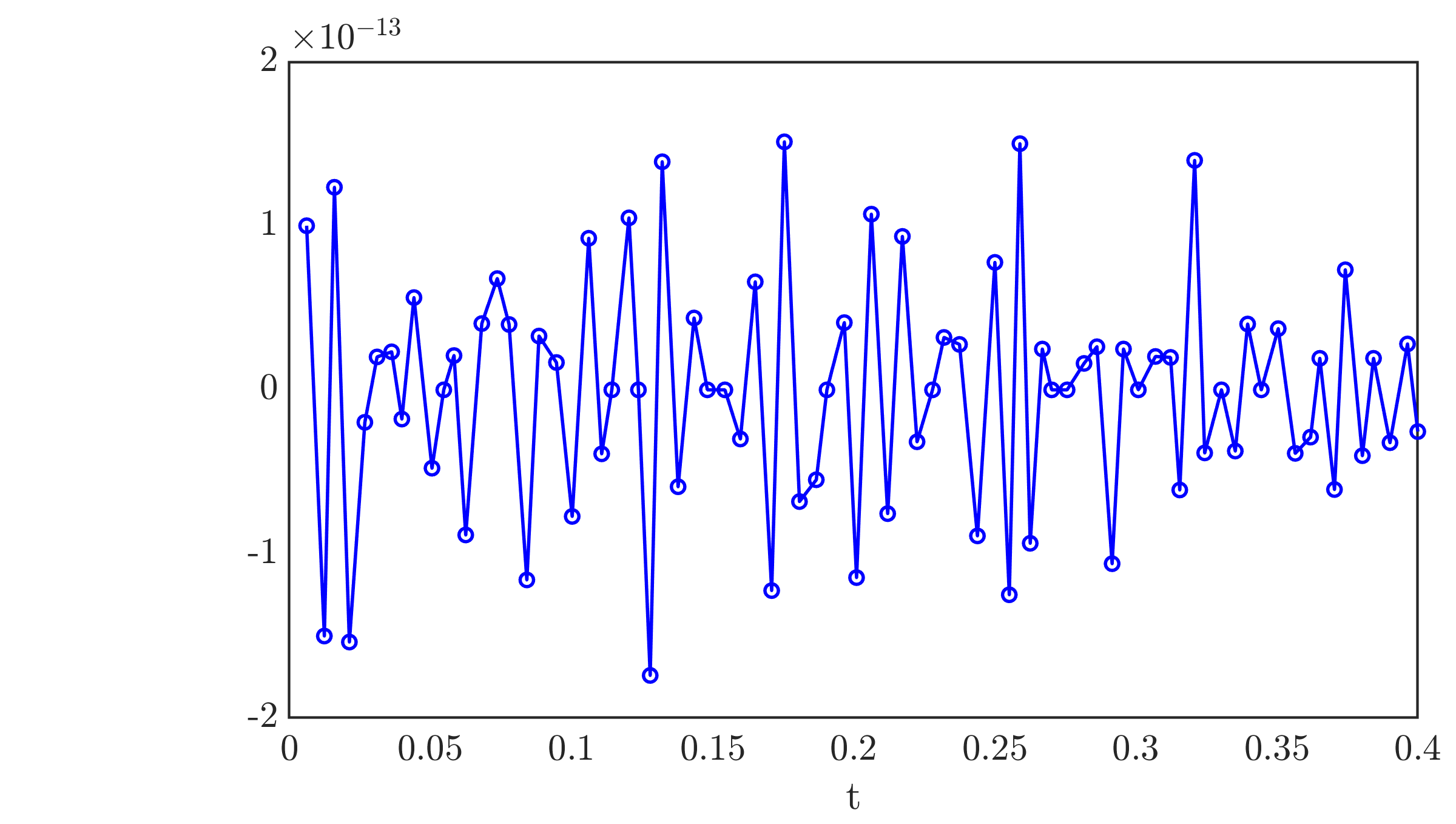}
    \hspace{-0.01\textwidth}%
        \includegraphics[width=0.32\textwidth,keepaspectratio]{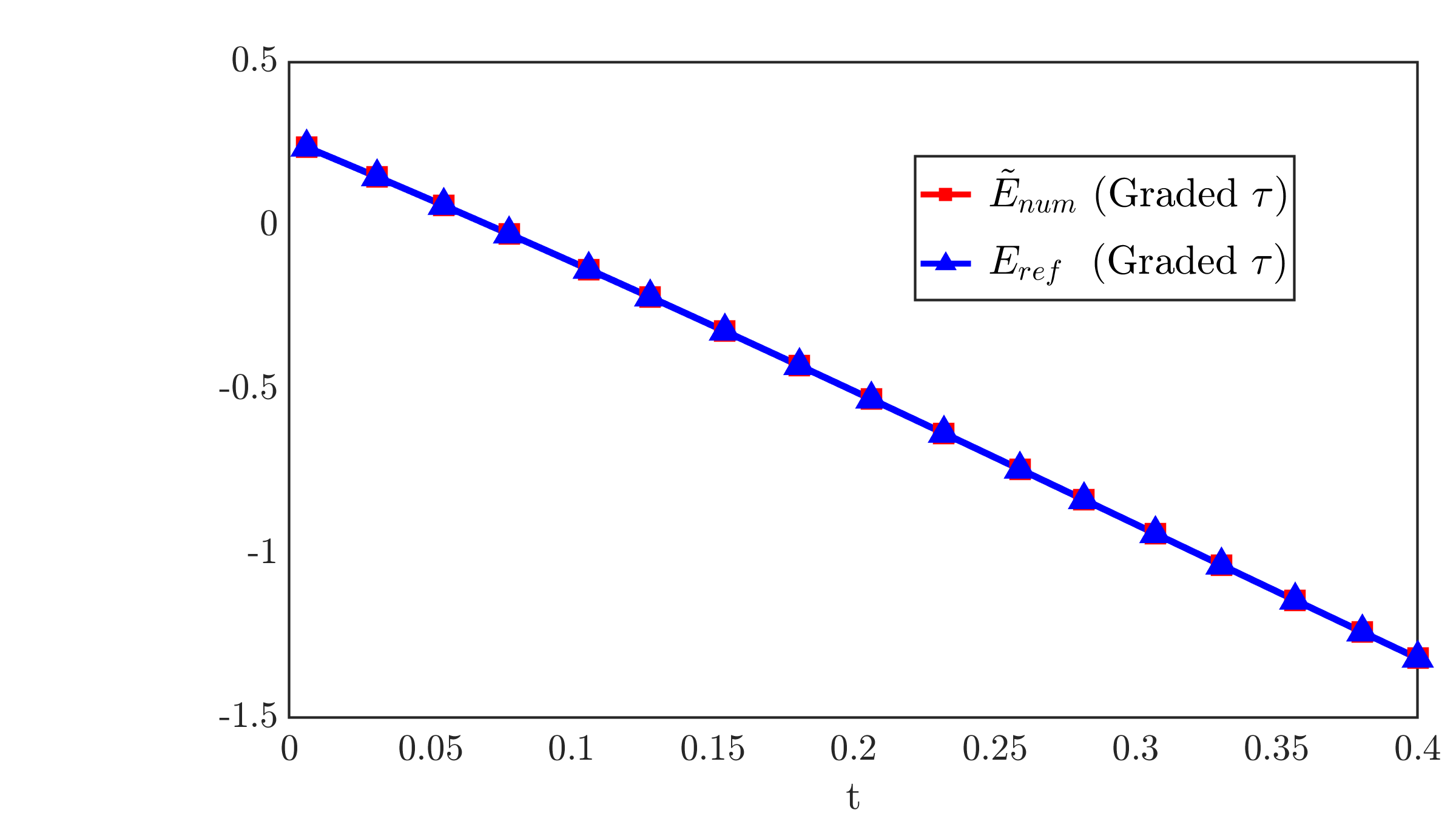}
    \setlength{\abovecaptionskip}{0.0cm} 
	  \setlength{\belowcaptionskip}{0.0cm}    
    \caption{Evolution of the nonlinear iteration count (left), the Lagrange multiplier \(\xi\) (middle), and energy (right) under variable time steps.}    
    \label{fig:Ite_xi:varnugradedstep}
\end{figure}

Finally, we evaluate the performance of the proposed variable-step scheme under the same setup as the constant-step case, using nonuniform graded temporal grids $ t_n=T \bigl(\frac{n}{N}\bigr)^{1.5}$. The largest adjacent-step ratio is $2^{1.5}-1\approx 1.828$, which satisfies the requirement $1/2\le r_n\le 2$  \cite{liao2023}. As shown in Table~\ref{tab:rate_variable_graded}, all solution variables achieve second-order temporal accuracy in their respective norms on the graded grids. In particular, both auxiliary-variable errors $|e_\alpha|$ and $|e_\beta|$ exhibit second-order convergence. Moreover, Fig.~\ref{fig:physics:varynugradedstep} demonstrates that at each time level, the projected density solutions are still in the initial bounds $[\rho_{\min}, \rho_{\max}]=[1.8, 2.2]$, with mass errors on the order of $10^{-14}$, all of which is consistent with Theorem~\ref{thm:variable-energy} for the variable-step scheme. Furthermore, similar observations regarding the nonlinear iteration count and the Lagrange multiplier are evident in Fig.~\ref{fig:Ite_xi:varnugradedstep}, where at most one nonlinear iteration is required per time level and  \(\xi\) remains on the order of \(10^{-13}\). Meanwhile, a similar energy behavior for the variable-step scheme is also observed.
\begin{table}[htbp] 
  \footnotesize
  \centering
  \caption{Temporal errors and convergence rates under variable time steps.}
  \label{tab:rate_variable_graded}
  \begin{tabular*}{0.9\columnwidth}{@{\extracolsep\fill}ccccc@{\extracolsep\fill}}
    \toprule
         $N$   & $10$ & $20$ & $40$ & $80$ \\
    \midrule
      $\|e_{u}\|$
        & 3.1613e-4 & 8.0972e-5 & 2.0477e-5 & 5.1483e-6 \\
      Rate
        & ---        & 2.00       & 2.00       & 2.00       \\[2pt]
      $\|e_{u}\|_{1}$
        & 1.9122e-3 & 4.8443e-4 & 1.2199e-4 & 3.0611e-5 \\
      Rate
        & ---        & 2.02       & 2.01       & 2.00       \\[2pt]
      $\|e_{p}\|$
        & 3.4078e-4 & 8.9534e-5 & 2.2875e-5 & 5.7778e-6 \\
      Rate
        & ---        & 1.97       & 1.99       & 1.99       \\[2pt]
      $\|e_{\rho}\|$
        & 3.1160e-5 & 8.5444e-6 & 2.2408e-6 & 5.7406e-7 \\
      Rate
        & ---        & 1.90       & 1.95       & 1.97       \\[2pt]
      $\|e_{\mathrm{div}}\|$
        & 2.1189e-5 & 5.7924e-6 & 1.5098e-6 & 3.8524e-7 \\
      Rate
        & ---        & 1.91       & 1.96       & 1.98       \\[2pt]
      $\|e_{\mathrm{div}}\|_{-1}$
        & 5.4326e-6 & 1.4884e-6 & 3.8830e-7 & 9.9114e-8 \\
      Rate
        & ---        & 1.90       & 1.96       & 1.98       \\[2pt]
      $|e_\alpha|$
        & 2.4817e-6 & 6.4926e-7 & 1.6531e-7 & 4.1669e-8 \\
      Rate
        & ---        & 1.97       & 1.99       & 2.00       \\[2pt]
      $|e_\beta|$
        & 1.3490e-7 & 3.6639e-8 & 9.5528e-9 & 2.4394e-9 \\
      Rate
        & ---        & 1.92       & 1.96       & 1.98       \\
    \bottomrule
  \end{tabular*}
\end{table}

\begin{example}[{Rayleigh--Taylor instability}]\label{subsec:e2} In this example, we test the Rayleigh--Taylor (RT) instability, a standard benchmark for variable-density incompressible flows \cite{guermond2000,espinozaM2AN2025}. Two initially quiescent fluids occupy $\Omega=(0,a/2)\times(-2a,2a)$, with the heavier fluid of density $\rho_2=3\rho_1$ initially above the lighter fluid of density $\rho_1$. The dynamic viscosity is given by $\nu=\nu_0\rho$, and the external force is represented by the gravity $\bm f=(0,-\rho g)^{\mathrm T}$. The initial density is given by
\begin{equation*}
   \rho^0(x,y)=\frac{\rho_2+\rho_1}{2} +\frac{\rho_2-\rho_1}{2}
     \tanh\left(\frac{y-\eta(x)}{0.01a}\right),
\end{equation*}
where 
$
\eta(x) = -0.1a\cos\left(\frac{2\pi x}{a}\right)
$
denotes the initial position of the perturbed interface.
Using $a$, $\rho_1$, $\sqrt{ag}$, and $\sqrt{a/g}$ as the reference length, density, velocity, and time scales, respectively, the Reynolds and Atwood numbers are defined as
\[
  \mathrm{Re}:=\frac{\rho_1a\sqrt{ag}}{\nu(\rho_1)} =\frac{a\sqrt{ag}}{\nu_0}=1000,
   \qquad
  \mathrm{At}:=\frac{\rho_2-\rho_1}{\rho_2+\rho_1}=\frac12.
\]
In the dimensionless simulations, we set $a=g=\rho_1=1$, which yields $\nu_0=10^{-3}$. 
The numerical viscosity at time level $t^{n+1}$ is evaluated as $\nu^{n+1}=\nu_0\rho^{n+1}$ (cf. Remark~\ref{rem:mu:1}).

As for the boundary conditions, no-slip conditions are imposed on the horizontal walls, while symmetry (free-slip) conditions are enforced on the vertical walls. To improve accuracy and efficiency, we adapt the spatial mesh using FreeFEM++'s \texttt{adaptmesh} routine, with the metric based on the computed fields  $\rho_h$, $u_{1h}$, and $u_{2h}$ and the local element size restricted to $0.01\leq h\leq0.05$. 
\end{example}

\begin{figure}
    \centering
    \includegraphics[width=0.7\textwidth,height=0.45\textheight,keepaspectratio]{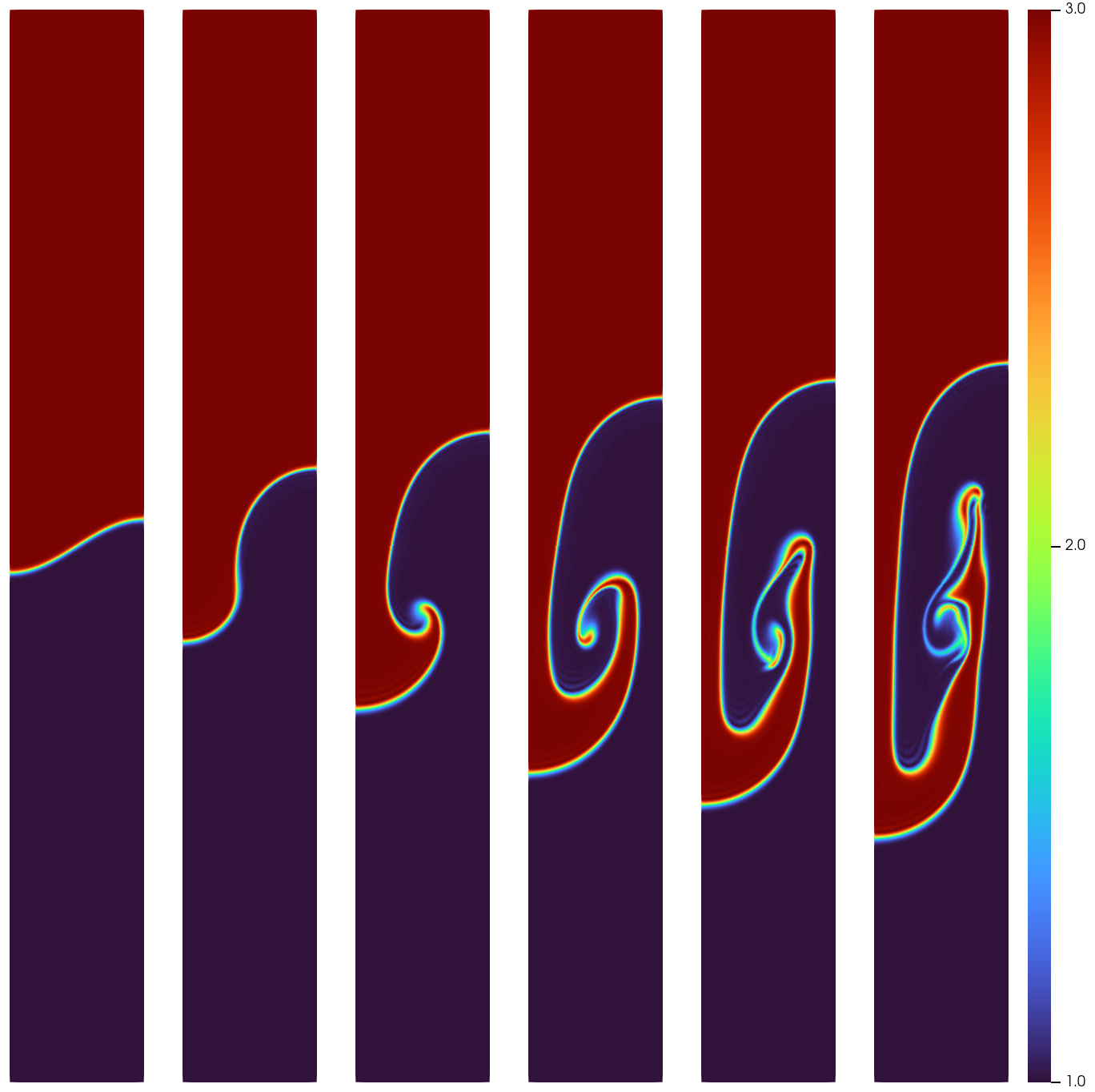}
    \setlength{\abovecaptionskip}{0.0cm} 
	  \setlength{\belowcaptionskip}{0.0cm}    
    \caption{Density evolution for the RT instability with $\tau=10^{-2}$ at $t=0,\,1.4,\,2.1,\,2.8,\,3.15,$ and $3.5$.}\label{fig:RT_var_dt1e2}
\end{figure}

\begin{figure}
    \centering
    \includegraphics[width=0.7\textwidth,height=0.45\textheight,keepaspectratio]{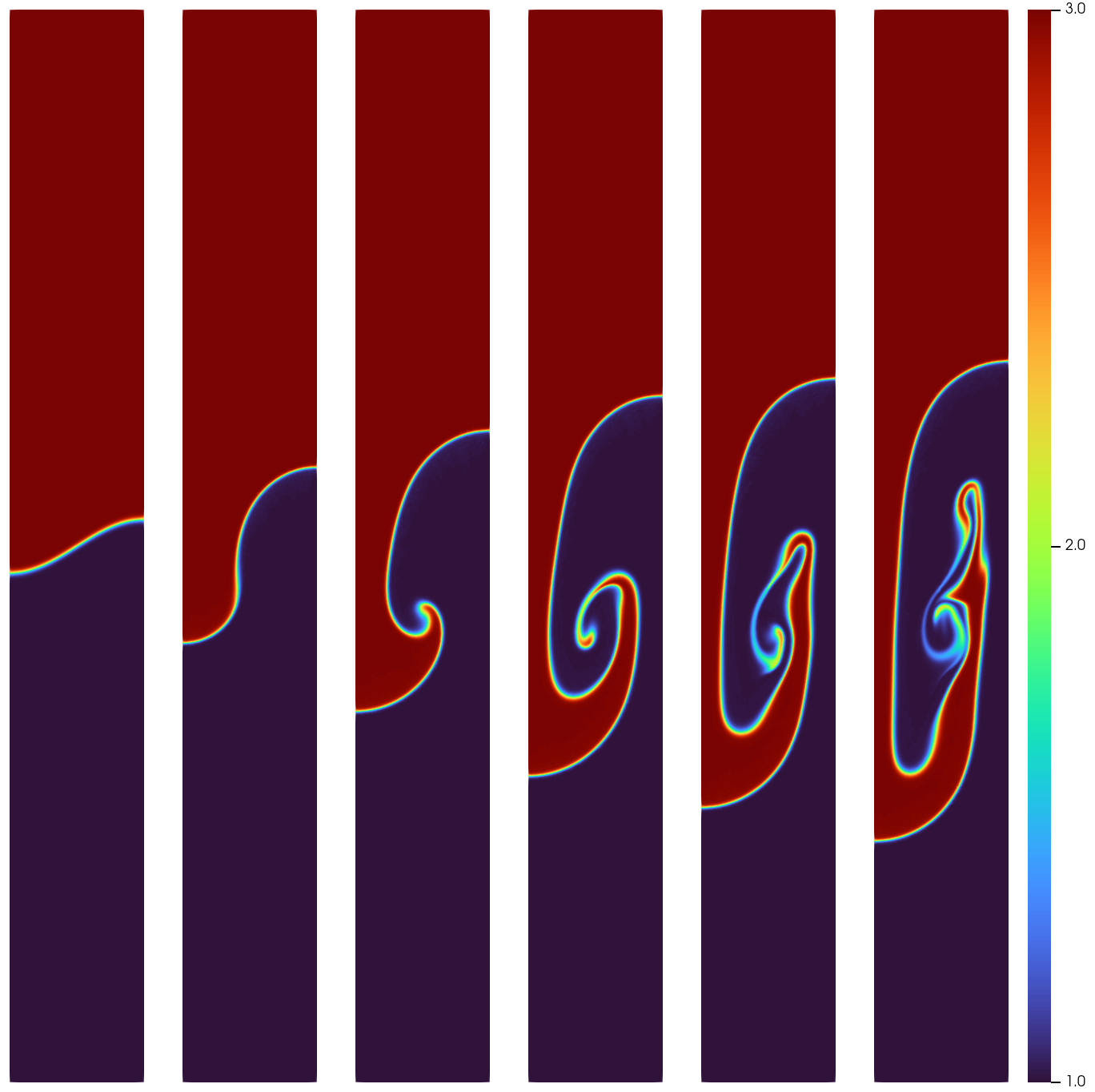}
    \setlength{\abovecaptionskip}{0.0cm} 
	  \setlength{\belowcaptionskip}{0.0cm}
    \caption{Density evolution for the RT instability with $\tau=5 \times 10^{-4}$ at $t=0,\,1.4,\,2.1,\,2.8,\,3.15,$ and $3.5$.}
    \label{fig:RT_var_dt5e4}
\end{figure}

\begin{figure}
    \centering
    \includegraphics[width=0.7\textwidth,keepaspectratio]{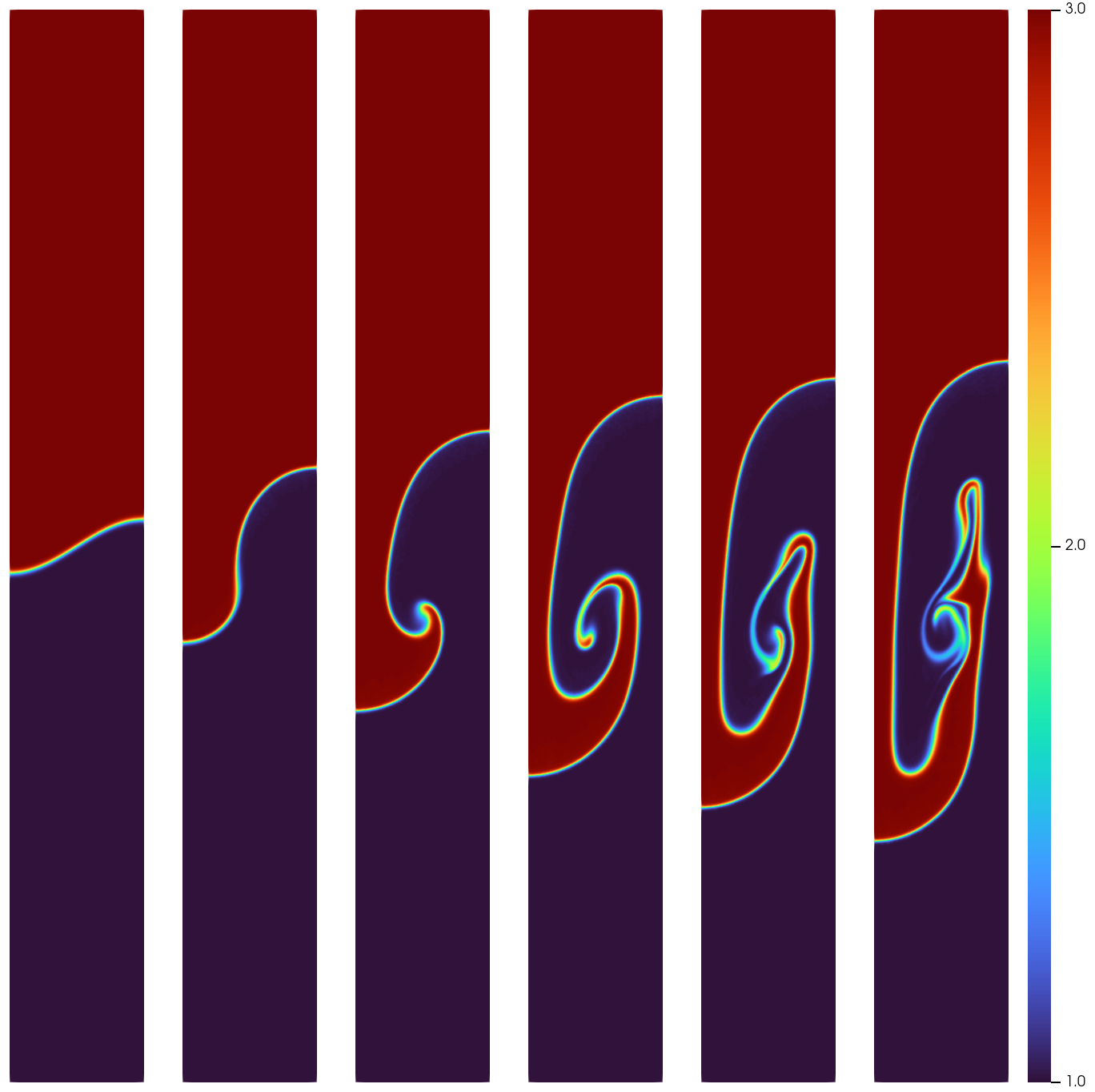}
    \setlength{\abovecaptionskip}{0.0cm} 
	  \setlength{\belowcaptionskip}{0.0cm}    
    \caption{Density evolution for the RT instability with adaptive time-stepping at $t=0,\,1.4,\,2.1,\,2.8,\,3.15,$ and $3.5$.}
    \label{fig:RT_var_ada}
\end{figure}

\begin{figure}[pos=!h]
    \centering
        \includegraphics[width=0.32\textwidth,keepaspectratio]{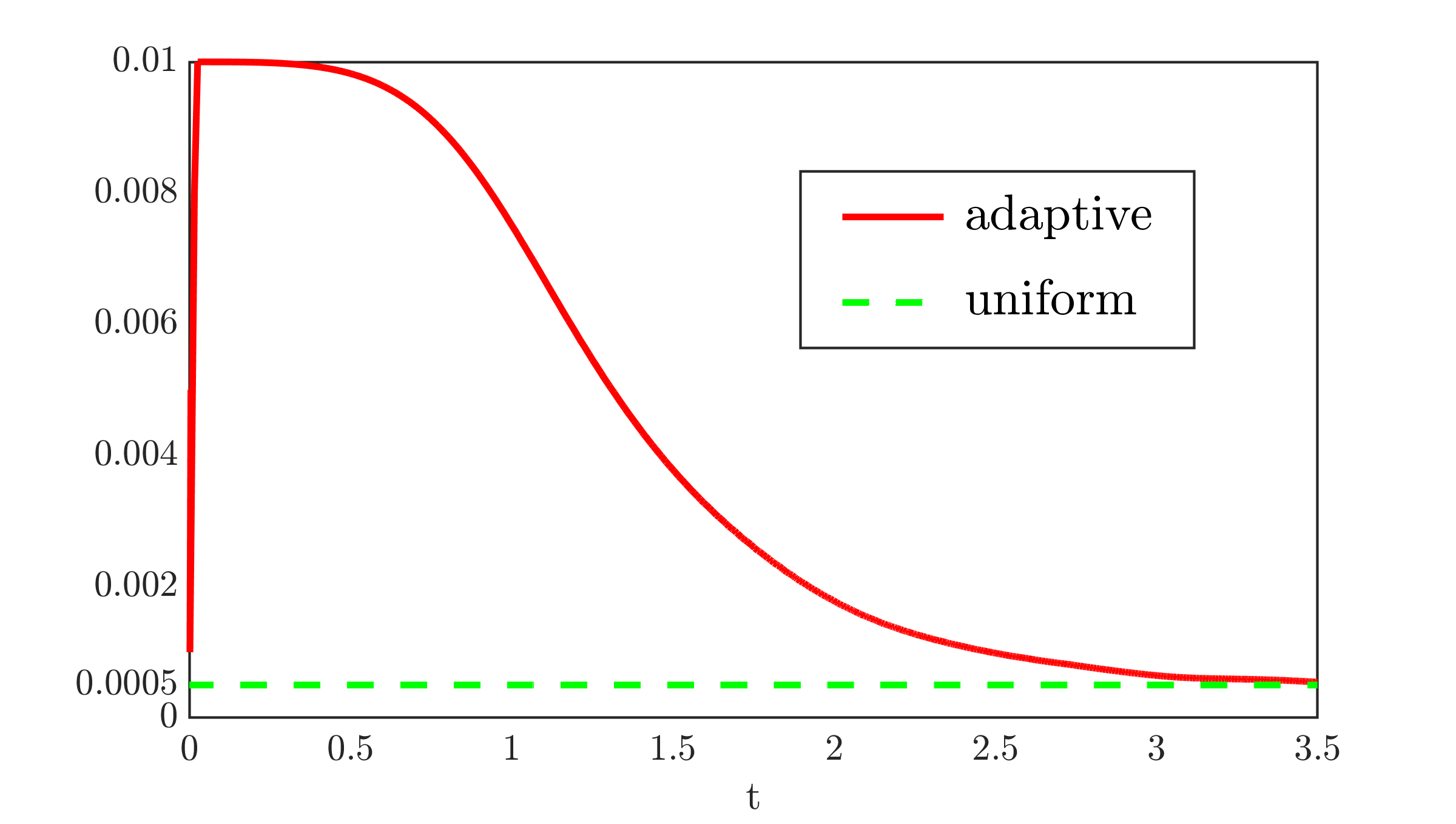}
    \hspace{-0.015\textwidth}
        \includegraphics[width=0.32\textwidth,keepaspectratio]{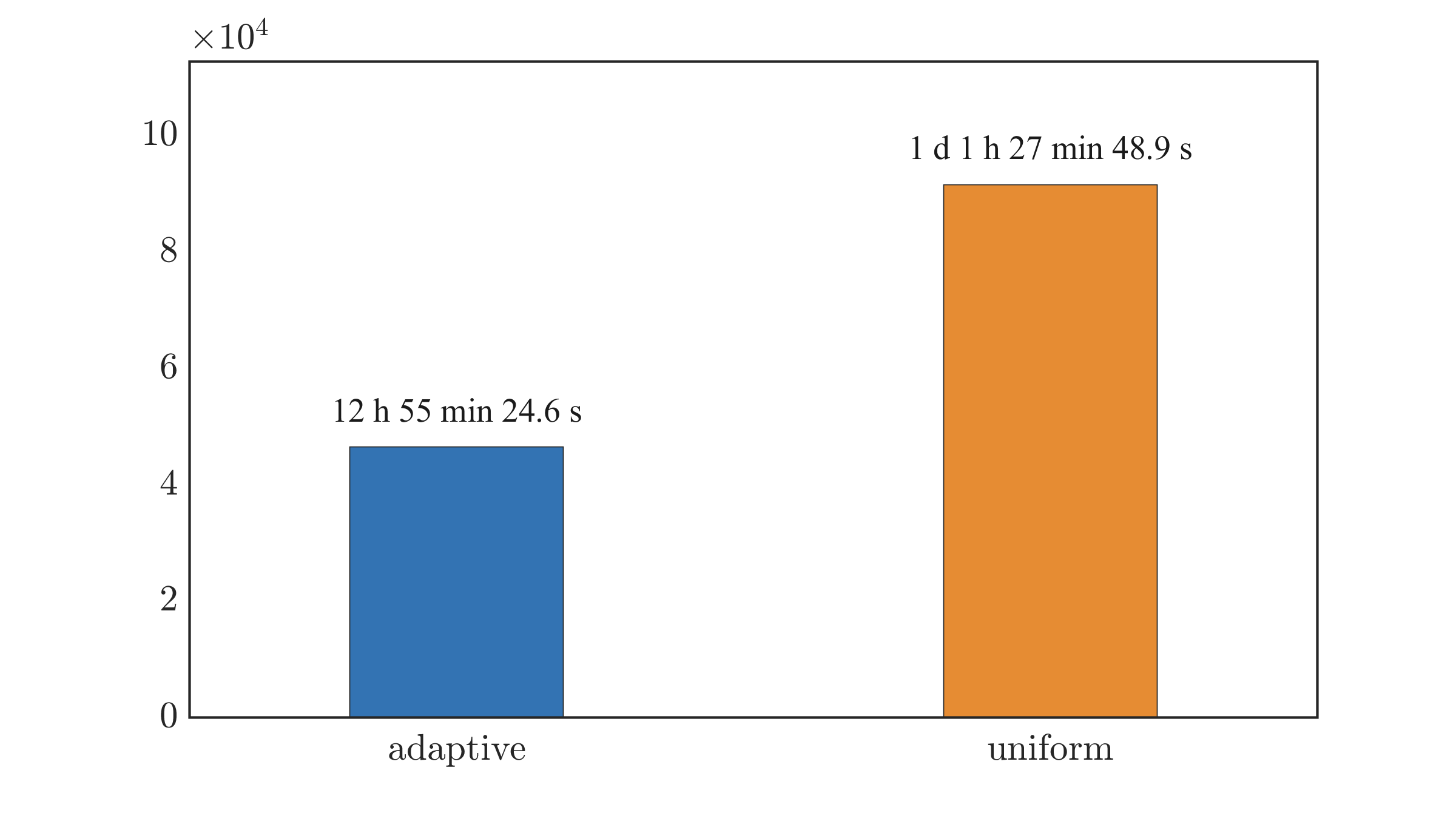}
    \hspace{-0.015\textwidth}
        \includegraphics[width=0.32\textwidth,keepaspectratio]{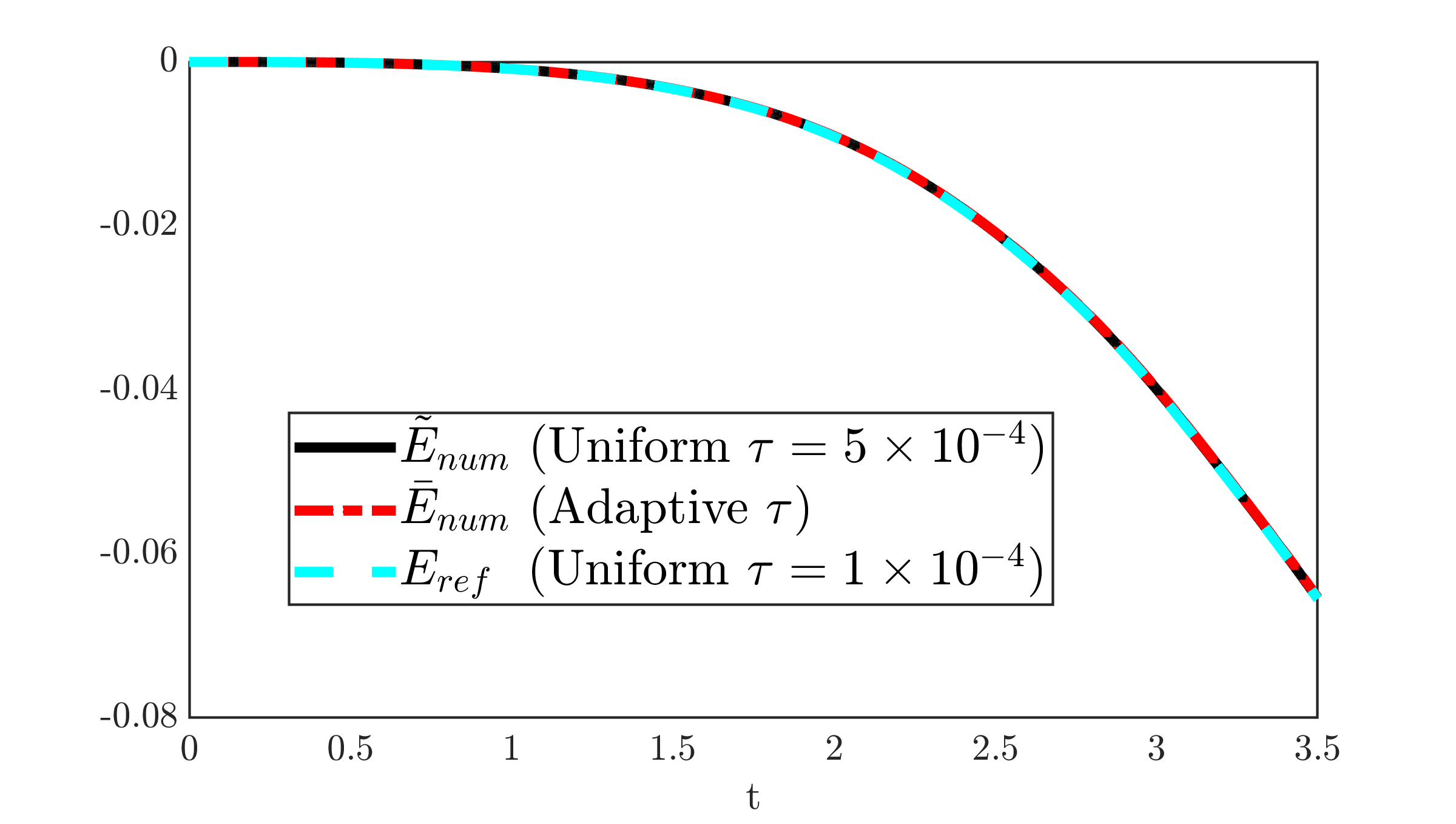}
    \setlength{\abovecaptionskip}{0.0cm} 
	  \setlength{\belowcaptionskip}{0.0cm}    
    \caption{Evolution of the adaptive time steps (left), CPU times (middle), and energies (right) for the RT instability.}  \label{fig:dtEnergy:RT_var}
\end{figure}

We first assess the proposed adaptive time-stepping strategy by comparing the density evolutions obtained with constant time steps $\tau=10^{-2}$ and $\tau= 5 \times 10^{-4}$. For the time-adaptivity computation, we set $\tau_{\min}=5 \times 10^{-4}$, $\tau_{\max}=10^{-2}$ and  $\kappa=1.1 \times 10^5$ in \eqref{vs:adaptive-indicator}. As shown in Figs.~\ref{fig:RT_var_dt1e2}--\ref{fig:RT_var_ada}, the three solutions remain close at early times. At $t=3.5$, however, the solution with $\tau=10^{-2}$ exhibits pronounced local thinning of the interface, whereas the adaptive solution remains close to that obtained with the finer step $\tau=5 \times 10^{-4}$. This indicates that the adaptive strategy reproduces the interfacial evolution of the finer constant step throughout the computation. Fig.~\ref{fig:dtEnergy:RT_var} reports the evolution of the adaptive time steps, CPU times, and energies. Three observations follow: (i) the adaptive time step stays near $\tau_{\max}=10^{-2}$ initially, then decreases progressively toward $\tau_{\min}=5\times10^{-4}$ till the end of the simulation; (ii) the adaptive strategy reduces the CPU time by about 49\% relative to the finer step $\tau=5\times10^{-4}$; and (iii) the modified discrete energies obtained with $\tau=5 \times 10^{-4}$ and with adaptive strategy both closely track the reference energy $E_{\mathrm{ref}}$ computed with $\tau=10^{-4}$, with all three curves decreasing monotonically.

\begin{figure}[pos=!h]
    \centering
        \includegraphics[width=0.32\textwidth,keepaspectratio]{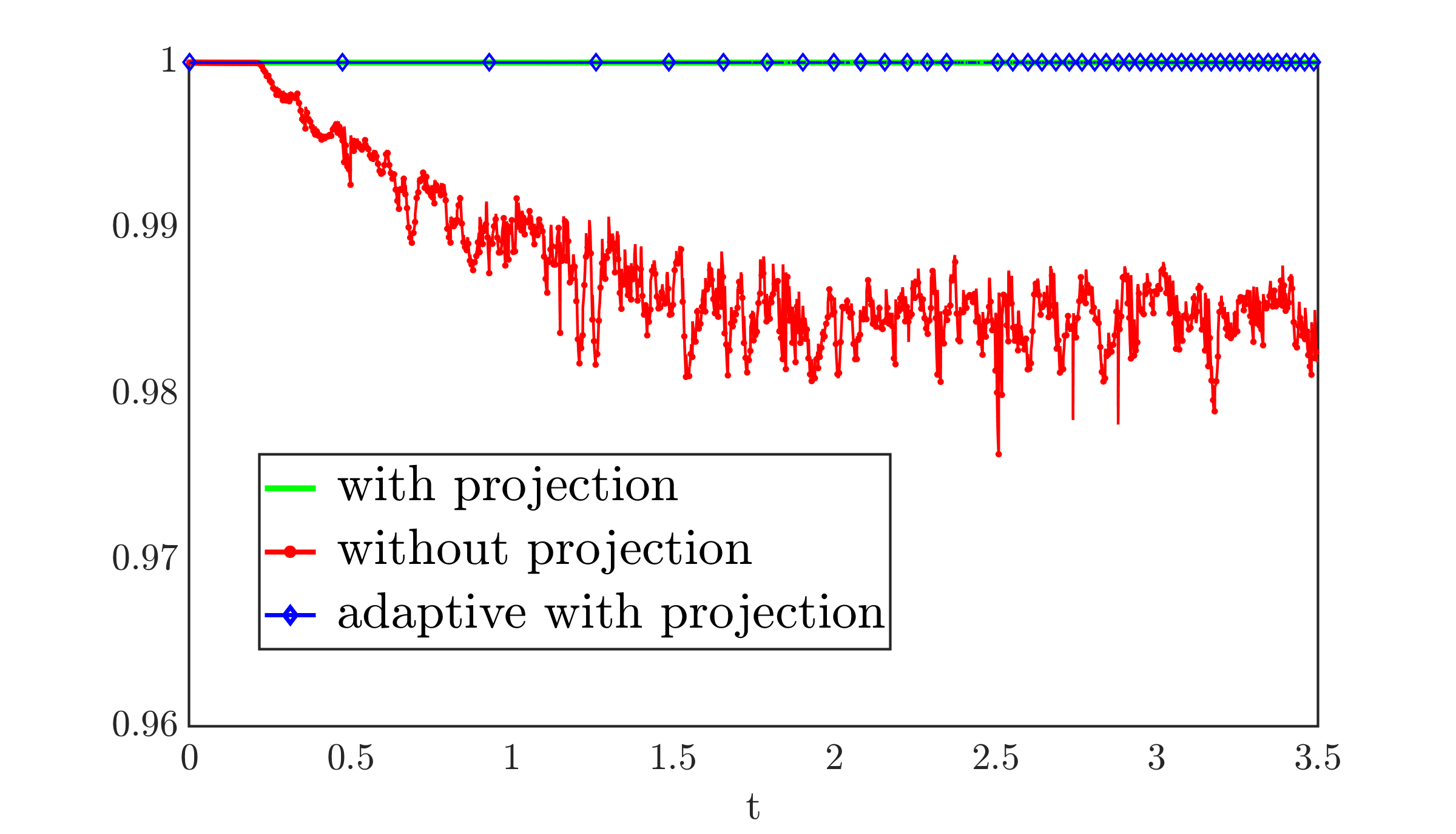}
    \hspace{-0.01\textwidth}%
        \includegraphics[width=0.32\textwidth,keepaspectratio]{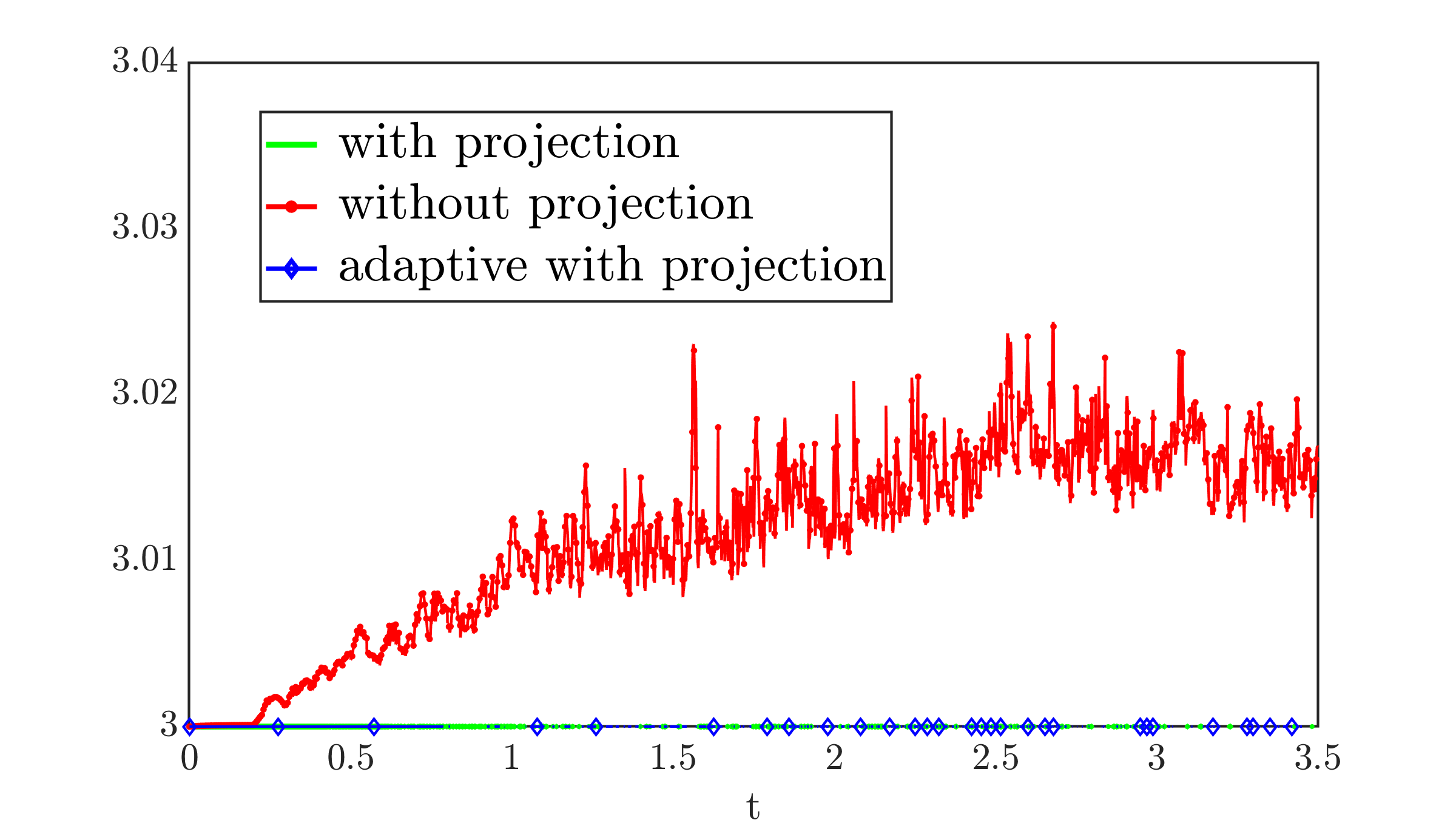}
    \hspace{-0.01\textwidth}%
        \includegraphics[width=0.32\textwidth,keepaspectratio]{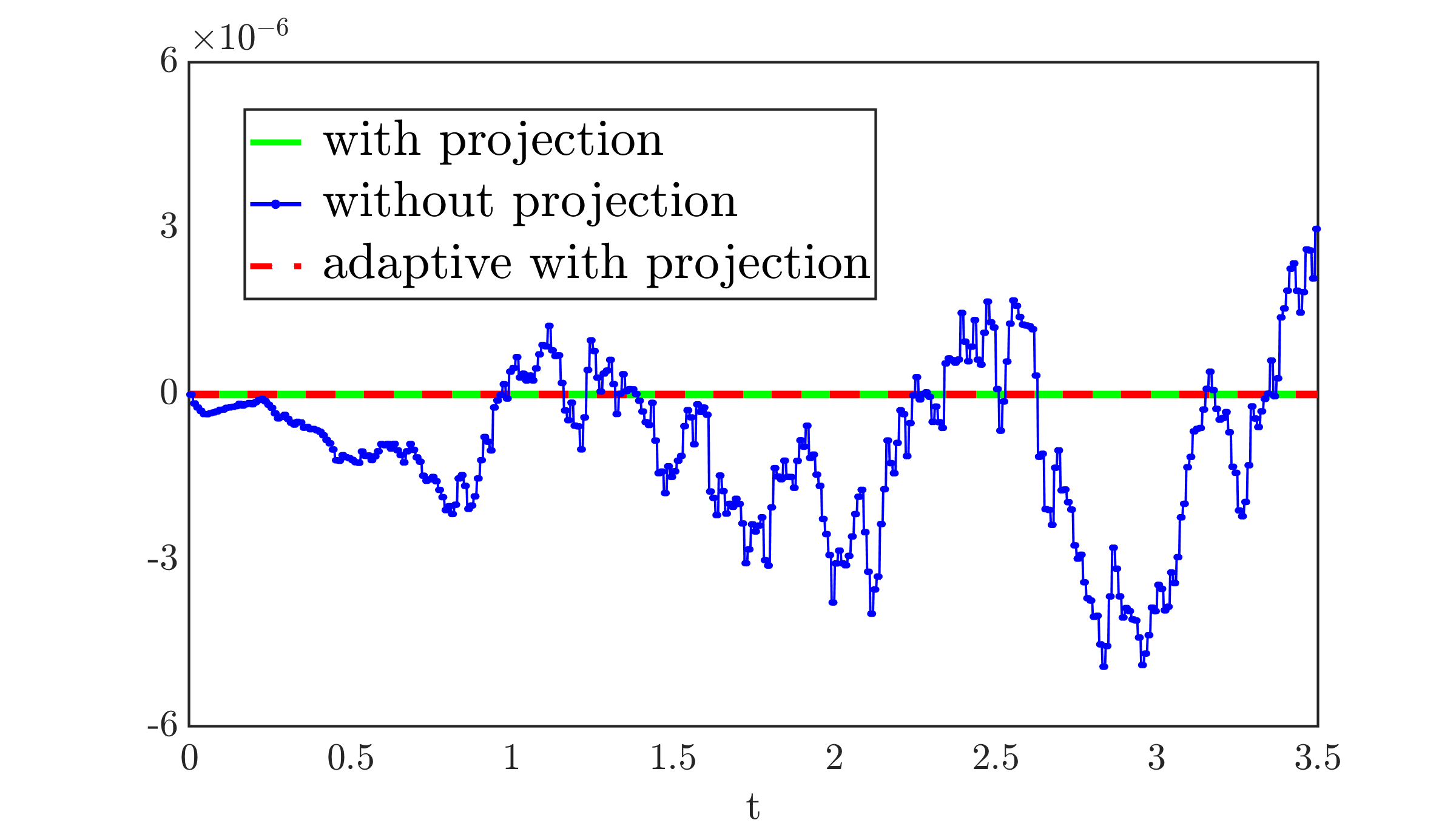}
    \setlength{\abovecaptionskip}{0.0cm} 
	  \setlength{\belowcaptionskip}{0.0cm}    
    \caption{Effect of the projection on density bounds  (min: left; max: middle) and mass error (right).}
    \label{fig:BP:RT_var}
\end{figure}

Next, we examine the effect of the projection step on density-bound preservation and mass conservation. In this simulation, the prescribed initial density bounds are $[\rho_{\min},\rho_{\max}]=[1,3]$. As shown in Fig.~\ref{fig:BP:RT_var}, without the projection \eqref{density:projection} or \eqref{vs:density:projection}, the minimum drops to about $0.98$ and the maximum density reaches about $3.02$, indicating clear overshoots and undershoots. Moreover, the mass error only reaches the level of $10^{-6}$ as time marching. With the projection applied, however, both schemes yield density solutions strictly within $[1,3]$, and the mass is conserved exactly, with the error near machine precision. These results directly confirm the bound-preserving and mass-conserving properties established in Theorems~\ref{thm:pme-conservation} and \ref{thm:variable-energy} for the proposed prediction--filter--projection schemes.

\begin{example}[Rising bubble] In this example, we consider the rising-bubble benchmark with a large density ratio and pronounced interfacial deformation \cite{liJCP2021,liJCP2024}. The physical parameters are listed in Table~\ref{tab:parameter}. The computational domain is $\Omega=[-0.01m,0.01m]\times[0.003m,0.02m]$, with homogeneous no-slip boundary conditions on the horizontal walls and symmetry (free-slip) conditions on the vertical walls. The initial air-water interface is regularized by
\begin{equation*}
  \rho^0(x,y)=\rho_{\mathrm{air}}+\frac{\rho_{\mathrm{water}}-\rho_{\mathrm{air}}}{2}\left[1+\tanh\left(\frac{d(x,y)-r}{\epsilon}\right)\right],
\end{equation*}
where $r=2.5 \times 10^{-3}$, $\epsilon=2.5 \times 10^{-4}$, and $d(x,y)=\sqrt{x^2+(y-0.0075)^2}$. Following \cite{liJCP2021,liJCP2024}, the dynamic viscosity varies linearly with the density:
\begin{equation*}
\nu(\rho):=\nu_{\mathrm{air}}+\frac{\rho-\rho_{\mathrm{air}}}{\rho_{\mathrm{water}}-\rho_{\mathrm{air}}}\bigl(\nu_{\mathrm{water}}-\nu_{\mathrm{air}}\bigr),
\end{equation*}
and is evaluated directly from the updated density, as discussed in Remark~\ref{rem:mu:1}.
The fluid is initially at rest, i.e., $\bm{u}^0(x,y)=\bm{0}$, and the buoyancy force is
$
\bm{f}=\left(0,\,(\rho-\rho_{\mathrm{water}})g\right)^\top
$ with $g=-9.80665m/s^2$. The local spatial mesh size is restricted to $0.01/200\leq h\leq0.01/64$, and the constant time step is $\tau=10^{-5}$.
\begin{table}[!t]
  \small
  \centering
  \caption{Physical parameters of the rising-bubble test.} \label{tab:parameter}
  \begin{tabular}{lllc}
    \toprule
    Parameter & Air & Water & Unit \\
    \midrule
    Density ($\rho$)     & 1.161       & 995.65      & $\mathrm{kg/m^3}$ \\
    Viscosity ($\nu$)    & $1.86 \times 10^{-5}$ & $7.977 \times 10^{-4}$ & $\mathrm{kg/(m \cdot s)}$ \\
    \bottomrule
  \end{tabular}
\end{table}
\end{example}

\begin{figure}[pos=!h]
    \centering
    \includegraphics[width=0.92\textwidth,keepaspectratio]{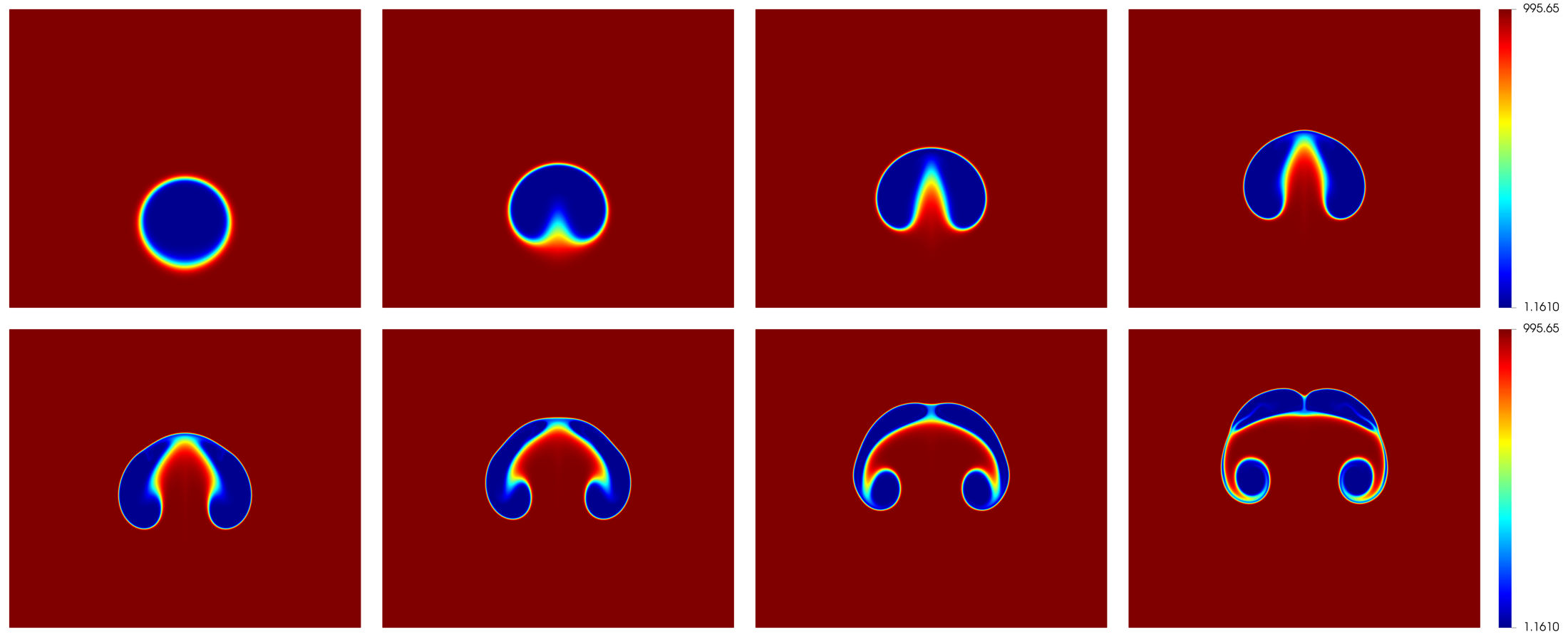}
    \setlength{\abovecaptionskip}{0.0cm} 
	  \setlength{\belowcaptionskip}{0.0cm}    
   \caption{Density evolution for the rising bubble at $t=0.01,\,0.02,\,0.03,\,0.04$ (top) and $0.05,\,0.06,\,0.07,\,0.08$ (bottom).}
    \label{fig:rising:rho}
\end{figure}

\begin{figure}[pos=!h]
    \centering
        \includegraphics[width=0.48\textwidth,keepaspectratio]{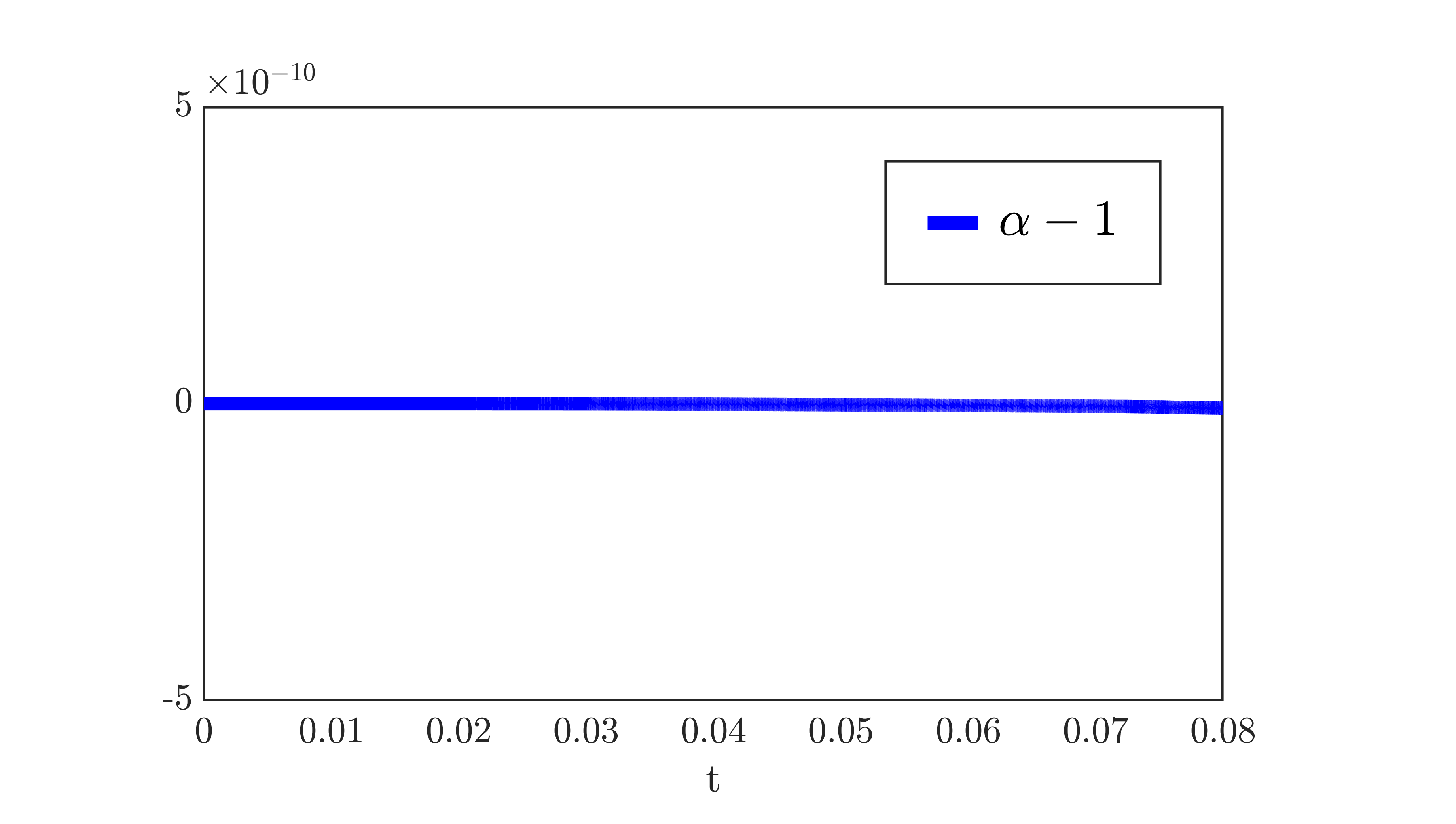}
    \hspace{-0.08\textwidth}
        \includegraphics[width=0.48\textwidth,keepaspectratio]{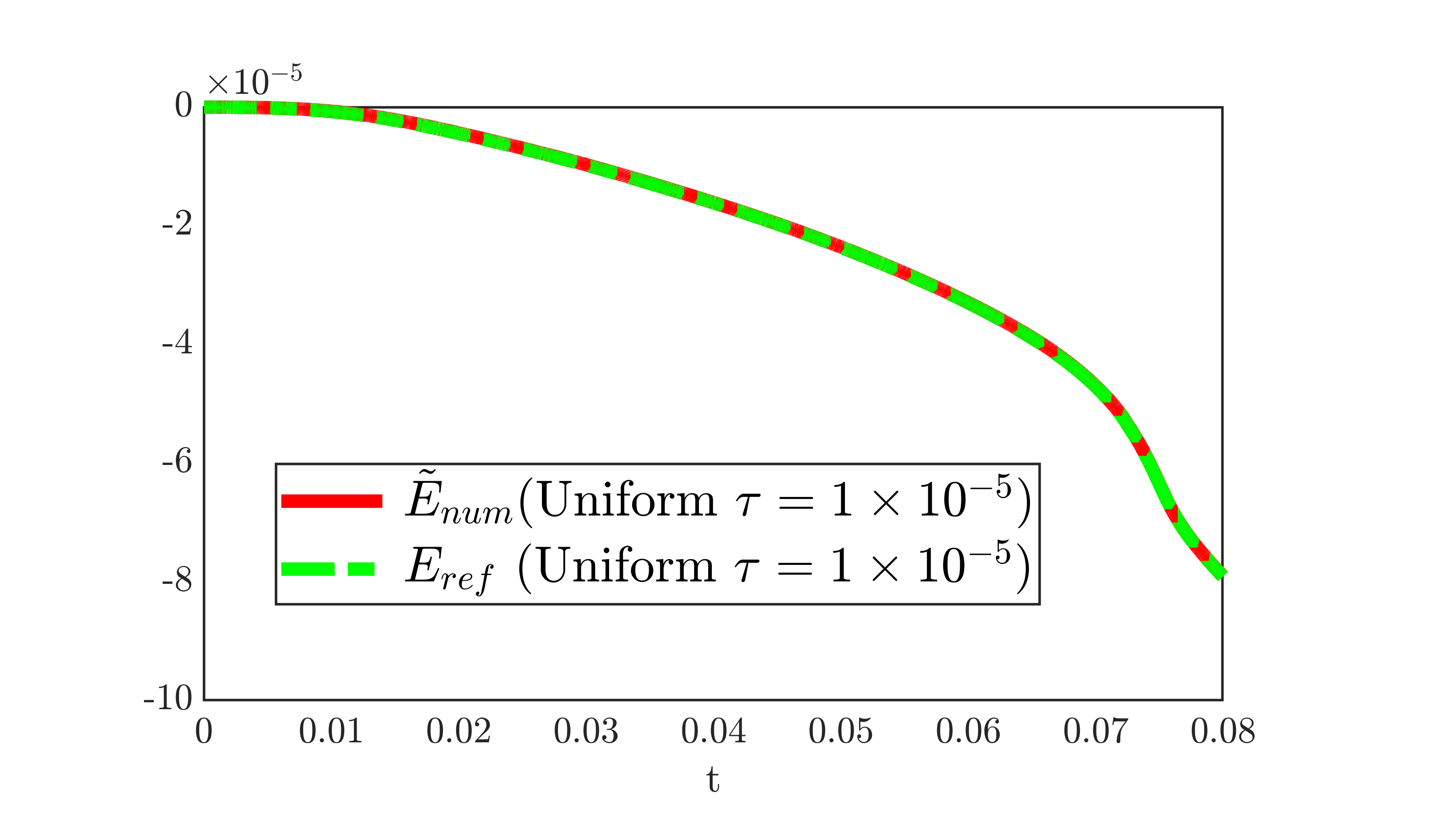}
    \setlength{\abovecaptionskip}{0.0cm} 
	  \setlength{\belowcaptionskip}{0.0cm}    
    \caption{Evolution of auxiliary-variable error (left) and energy (right) for the rising bubble.}
    \label{fig:rising:BP}
\end{figure}

This test features a large density ratio $\rho_{\mathrm{water}}/\rho_{\mathrm{air}}\approx858$ and a viscosity ratio $\nu_{\mathrm{water}}/\nu_{\mathrm{air}}\approx43$.  Fig.~\ref{fig:rising:rho} presents the density evolution from $t=0.01s$ to $0.08s$  for the prediction--filter--projection scheme \eqref{density:predict}--\eqref{2th:TF:M}. As can be observed, the initially circular bubble rises and undergoes substantial deformation while retaining the overall symmetry, and the computed interface shapes agree well qualitatively with those reported in \cite{liJCP2021,liJCP2024}. Moreover, Fig.~\ref{fig:rising:BP} shows the auxiliary-variable error and energy evolution. We can see that the error $\alpha-1$ remains on the order of $10^{-10}$ throughout the simulation, confirming that the discrete auxiliary variable stays very close to its continuous value $\alpha(t) \equiv 1$. Moreover, the modified discrete energy $\widetilde E_{\mathrm{num}}$ nearly coincides with the reference energy $E_{\mathrm{ref}}$ throughout the evolution, and both curves decrease monotonically, consistent with the modified-energy dissipation law established in \eqref{Enum:W}.
 
\begin{example}[Falling droplet] In this example, we consider a heavy droplet falling through a lighter surrounding fluid and subsequently merging with the heavy-fluid layer below. The computational domain is $\Omega=(0,L)\times(0.7L,2L)$ with $L=1$.  Homogeneous free-slip boundary conditions are imposed on all boundaries.  The fluid is initially at rest, and gravity acts in the negative $y$-direction through $\bm{f}=(0,-\rho g)^{\mathrm T}$ with $g=1$. The constant dynamic viscosity $\nu$ is chosen such that 
\begin{equation*}
  Re:=\frac{\rho_{\min}L^{3/2}g^{1/2}}{\nu}=3132.
\end{equation*}
Initially, a circular heavy droplet of radius $0.2$ is centered at $(0.5,1.75)$, and the initial density is
\begin{equation*}
 \rho^{0}(x,y)=
    \begin{cases}
      100, & 0.7\leq y\leq1\ \text{or}\ r\leq 0.2,\\
      1, & 1<y\leq2\ \text{and}\ r>0.2,
    \end{cases}
  \qquad r:=\sqrt{(x-0.5)^2+(y-1.75)^2}.
\end{equation*}  
\end{example}

Since this initial density is discontinuous, its direct interpolation into the continuous $P_3$ finite element space may generate spurious oscillations and nonphysical undershoots. To avoid this, we apply a smoothing step to the initial condition by solving the following  regularized $L^2$-projection problem 
\begin{equation*}
    \bigl(\widetilde{\rho}^{0},w\bigr)
       +\eta_{\mathrm s}\bigl(h^{2}\nabla\widetilde{\rho}^{0},\nabla w\bigr)
       =\bigl(\rho^{0},w\bigr),\qquad \forall\,w\in M_h,
\end{equation*}
as considered in \cite{lundgrenNazarov2024}.
Here, $M_h$ is the continuous $P_3$ finite element space, 
and the smoothing parameter is set to $\eta_{\mathrm s}=7$. In this example, the spatial discretization employs 117,691 and 52,461 DoFs for the $P_3$ and $P_2$ finite element spaces, respectively. Moreover, we set $Q=10^{-5}$ and $Z=10^{-5}$.

\begin{figure}
    \centering
        \includegraphics[width=0.9\textwidth,keepaspectratio]{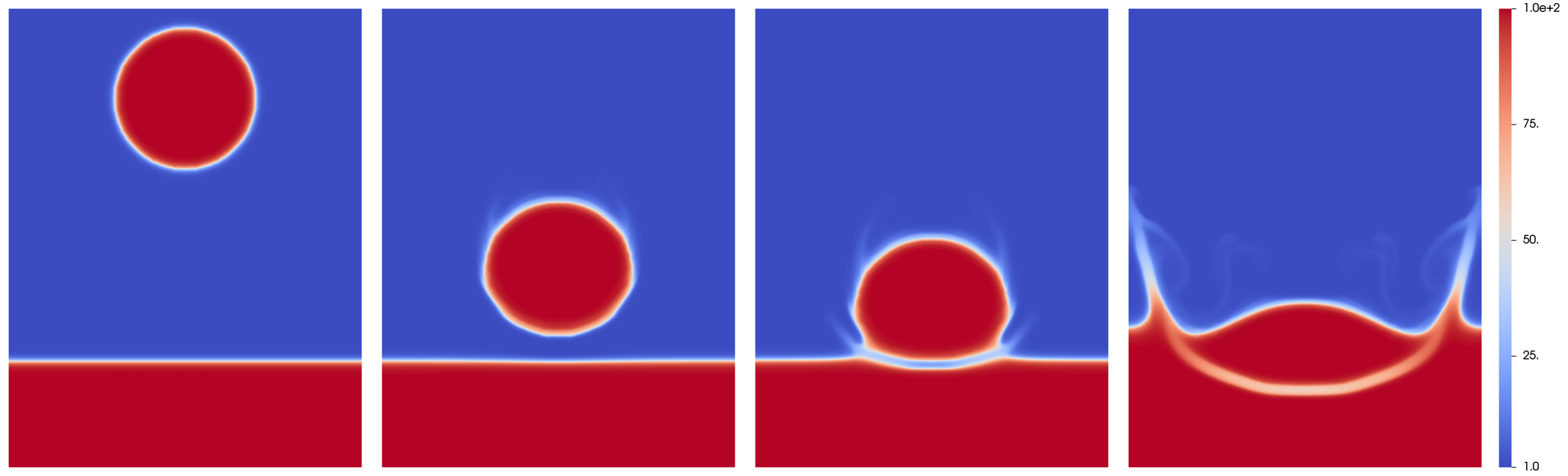}\\
        \includegraphics[width=0.9\textwidth,keepaspectratio]{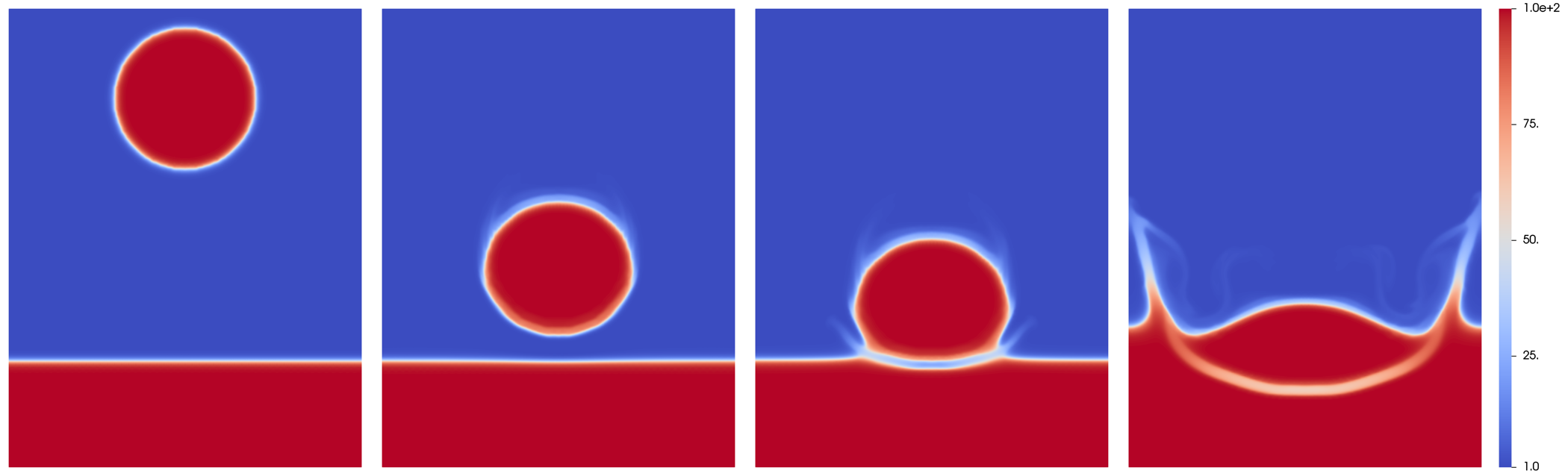}
    \hspace{-0.02\textwidth}%
    \setlength{\abovecaptionskip}{0.0cm} 
	  \setlength{\belowcaptionskip}{0.0cm}    
    \caption{Density evolution for the falling droplet at \(t=0.1,\,1.0,\,1.1,\) and \(1.3\) for both uniform (top) and adaptive (bottom) time-stepping.}
    \label{fig:falling:rho}
\end{figure}
\begin{figure}
    \centering
        \includegraphics[width=0.32\textwidth,keepaspectratio]{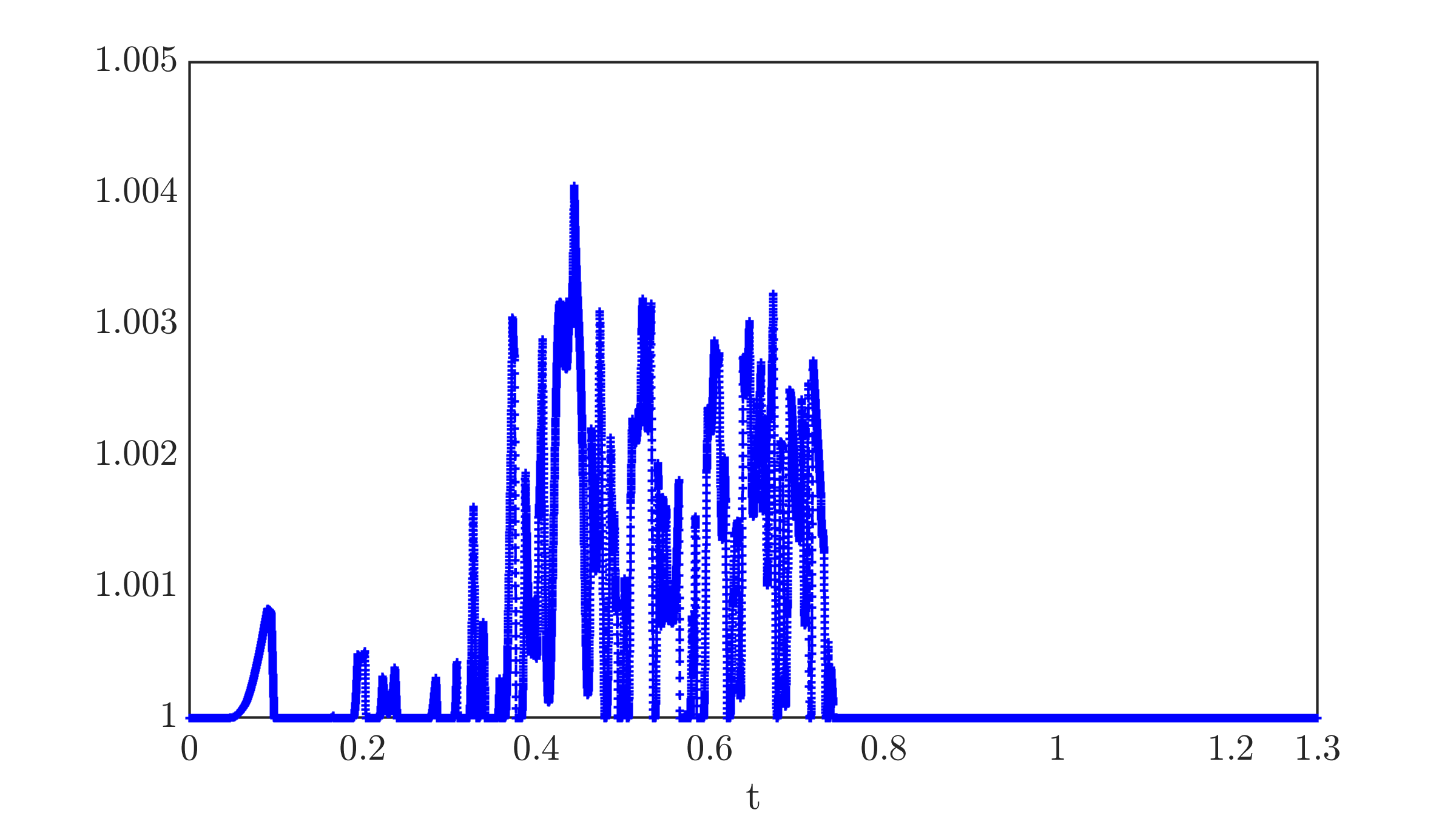}
    \hspace{-0.01\textwidth}
        \includegraphics[width=0.32\textwidth,keepaspectratio]{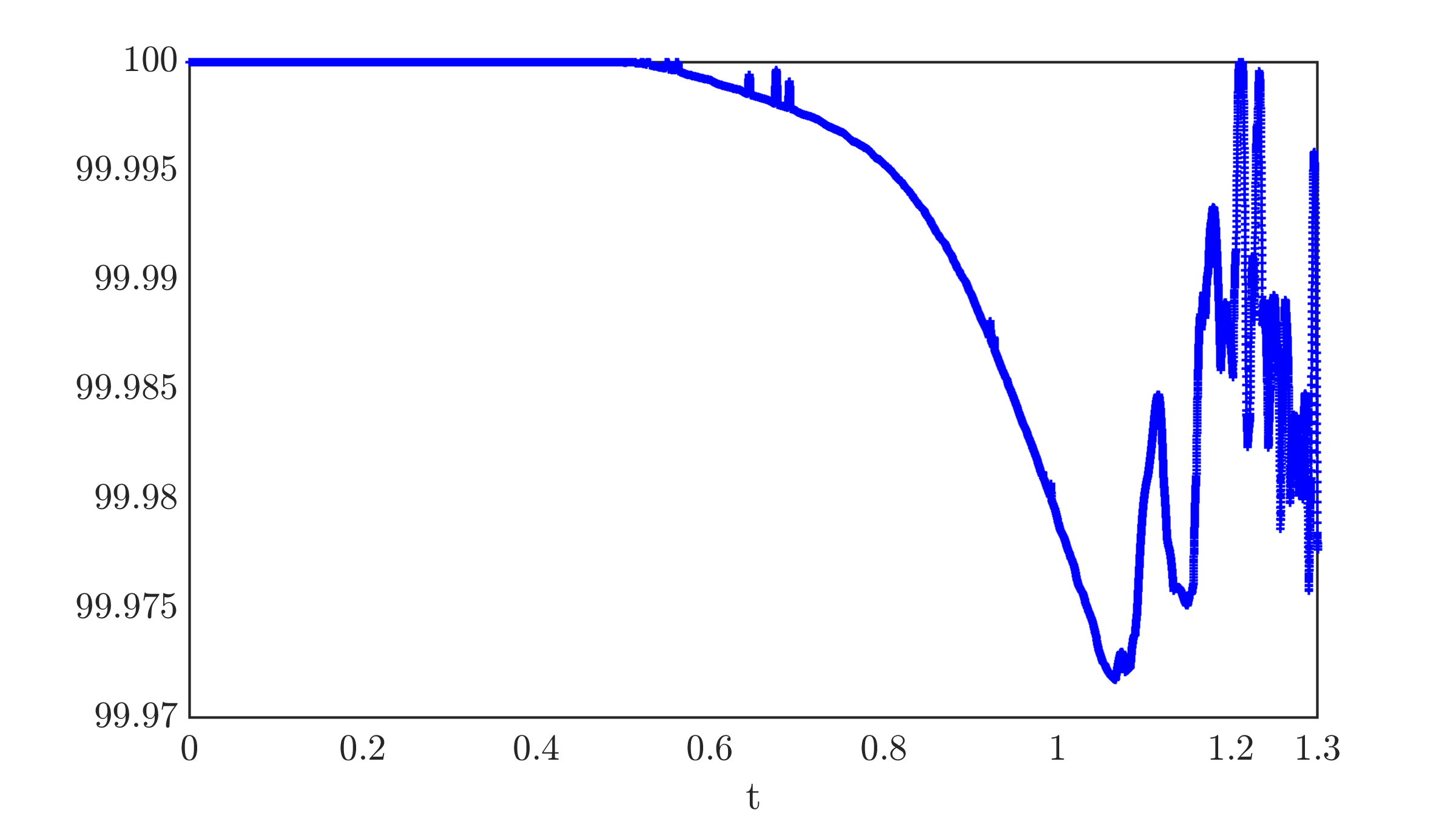}
    \hspace{-0.01\textwidth}%
        \includegraphics[width=0.32\textwidth,keepaspectratio]{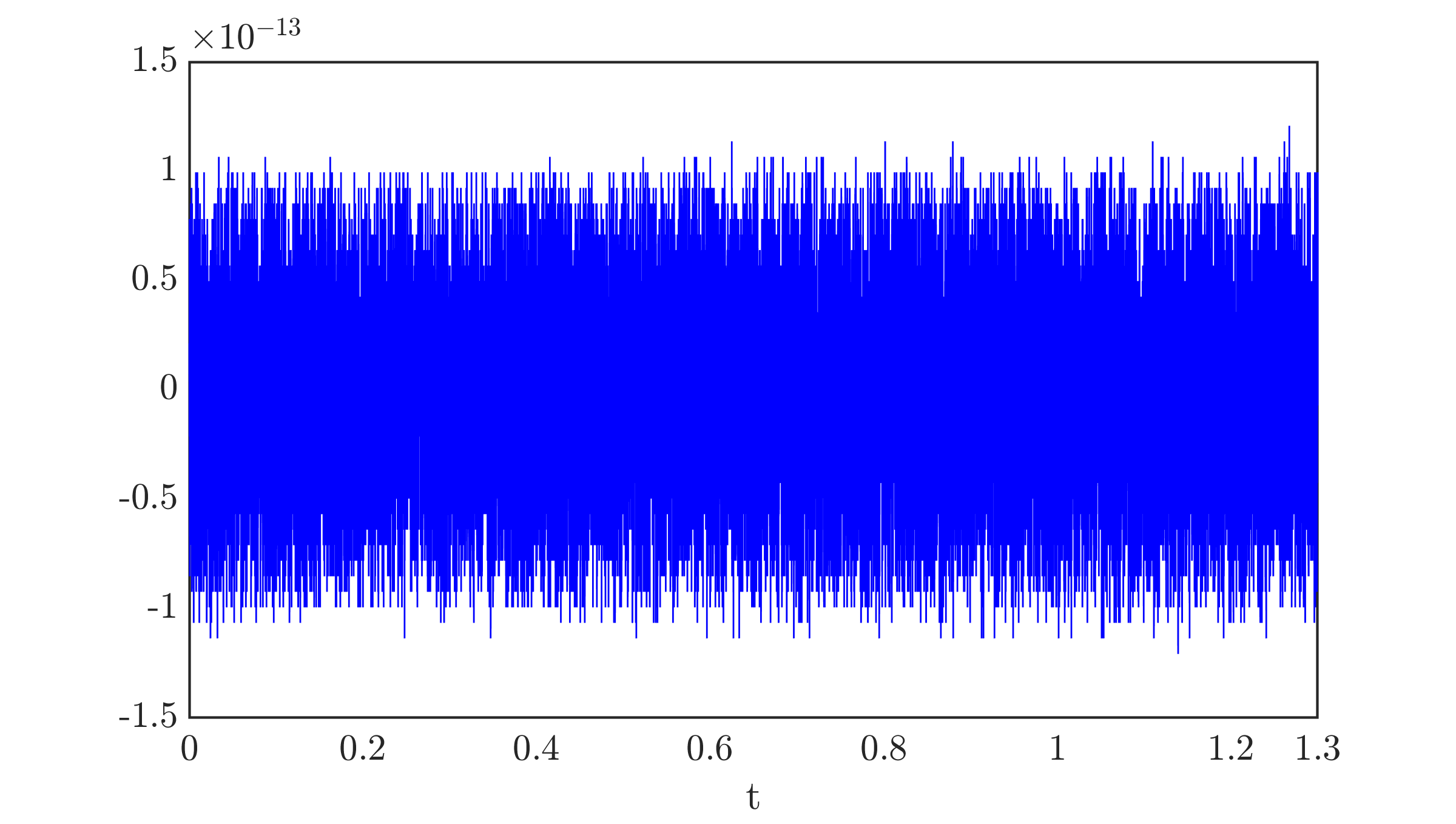}
    \setlength{\abovecaptionskip}{0.0cm} 
	  \setlength{\belowcaptionskip}{0.0cm}    
    \caption{Evolution of density bounds (min: left; max: middle) and mass error (right) for the falling droplet.}    
    \label{fig:physics:falling bubbles}
\end{figure}
\begin{figure}[pos=!h]
    \centering
        \includegraphics[width=0.48\textwidth,keepaspectratio]{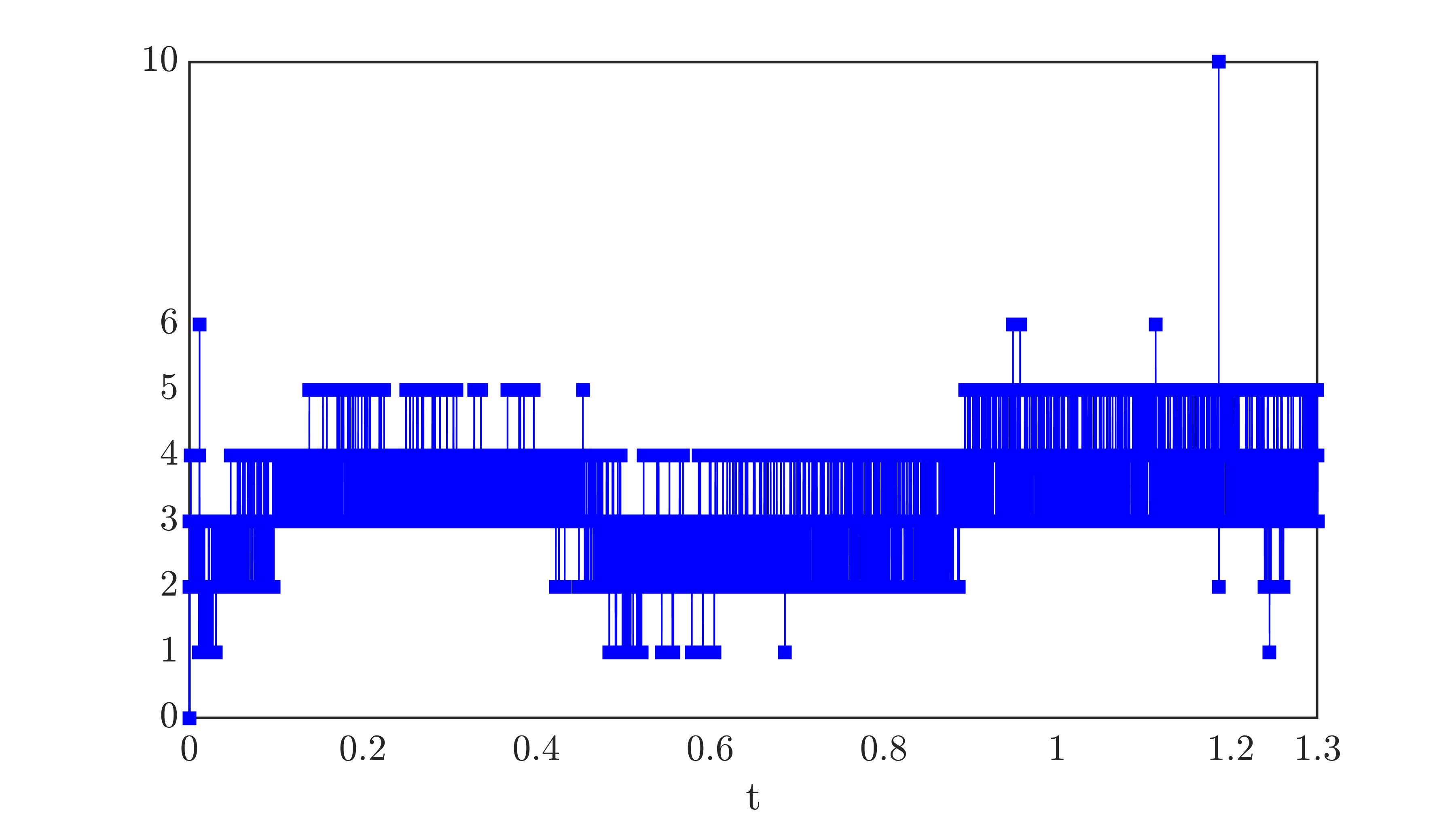}\hspace{-0.03\textwidth}%
       \includegraphics[width=0.48\textwidth,keepaspectratio]{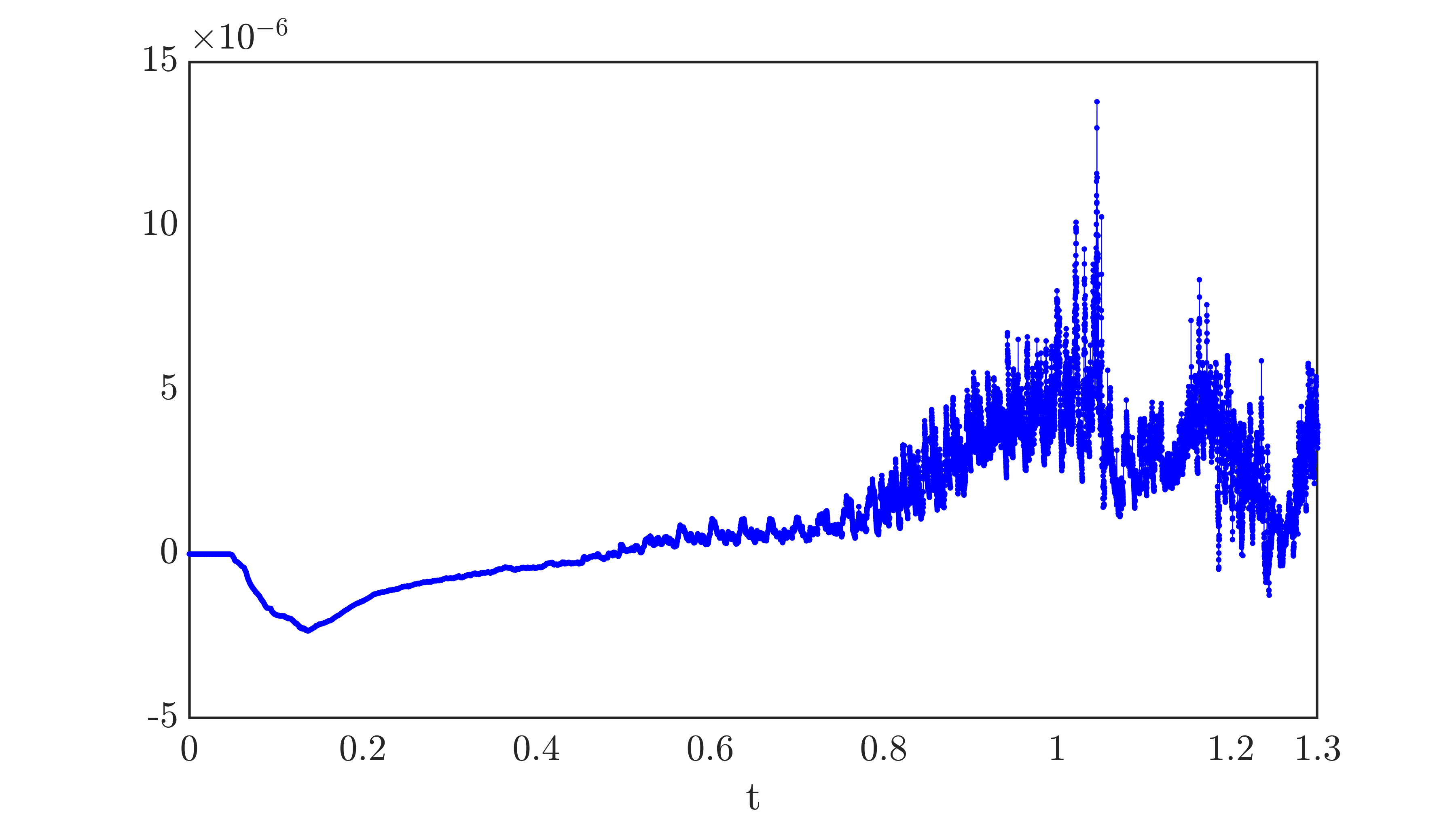}
    \setlength{\abovecaptionskip}{0.0cm} 
	  \setlength{\belowcaptionskip}{0.0cm}       
   \caption{Evolution of the nonlinear iteration count (left) and the Lagrange multiplier $\xi$ (right) for the falling droplet.}
\label{fig:Ite_xi:falling_bubbles}
    \label{fig:Ite_xi:falling bubbles}
\end{figure}

Fig.~\ref{fig:falling:rho} presents the density evolutions obtained with the uniform and adaptive time-stepping strategies, where the constant time step is $6\times10^{-5}$ and, for the adaptive computation, $\tau_{\min}=1\times10^{-5}$, $\tau_{\max}=1\times10^{-3}$, and $\kappa=1.3$. It can be seen that the adaptive solution closely agrees with the constant-step solution at the reported times, demonstrating that the time-adaptivity strategy reproduces the constant-step interfacial evolution while permitting larger time steps during the early stage. Consistent with \cite{lundgrenNazarov2024}, the heavy droplet falls through the lighter fluid and eventually merges with the underlying heavy-fluid layer; in the absence of surface tension, it undergoes noticeable deformation during its descent. Next, we further examine the physics-preserving properties. As shown in Fig.~\ref{fig:physics:falling bubbles}, despite the discontinuous initial density and strong interfacial deformation, the numerical density remains within the initial bounds $[1,100]$, while the mass error stays at the $10^{-13}$ level. The slight inward shifts of the minimum and maximum density values may stem from the regularized initial projection and residual-based stabilization. Moreover, as shown in Fig.~\ref{fig:Ite_xi:falling_bubbles}, the Lagrange multiplier $\xi$ in the projection step remains at the $10^{-5}$ level, while the nonlinear iteration count is typically between two and five, with only a few time levels outside this range and a maximum of ten. 

\begin{figure}[pos=!h]
    \centering
    \includegraphics[width=0.45\textwidth,keepaspectratio]{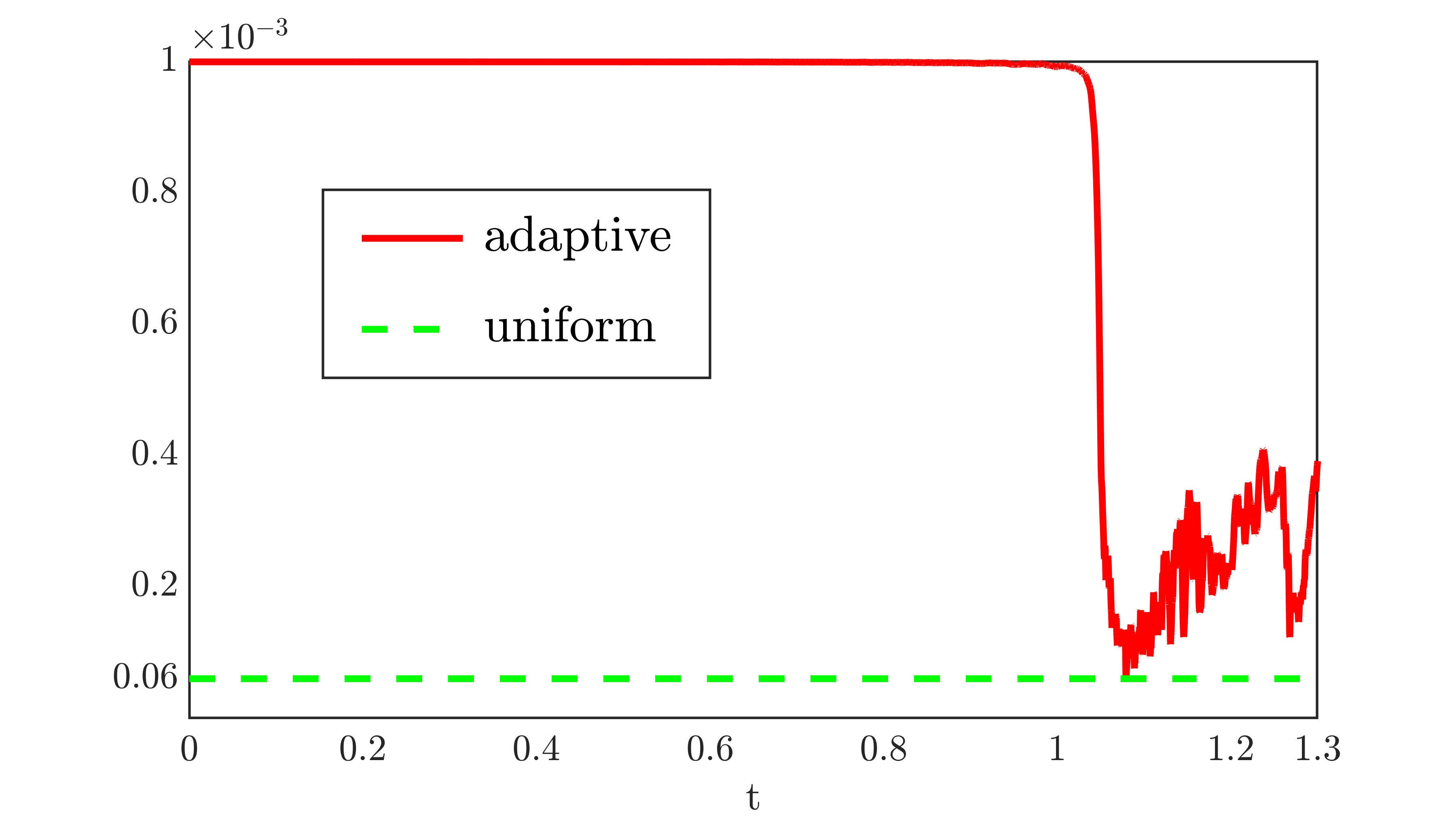}
     \hspace{-0.01\textwidth}%
  \includegraphics[width=0.45\textwidth,keepaspectratio]{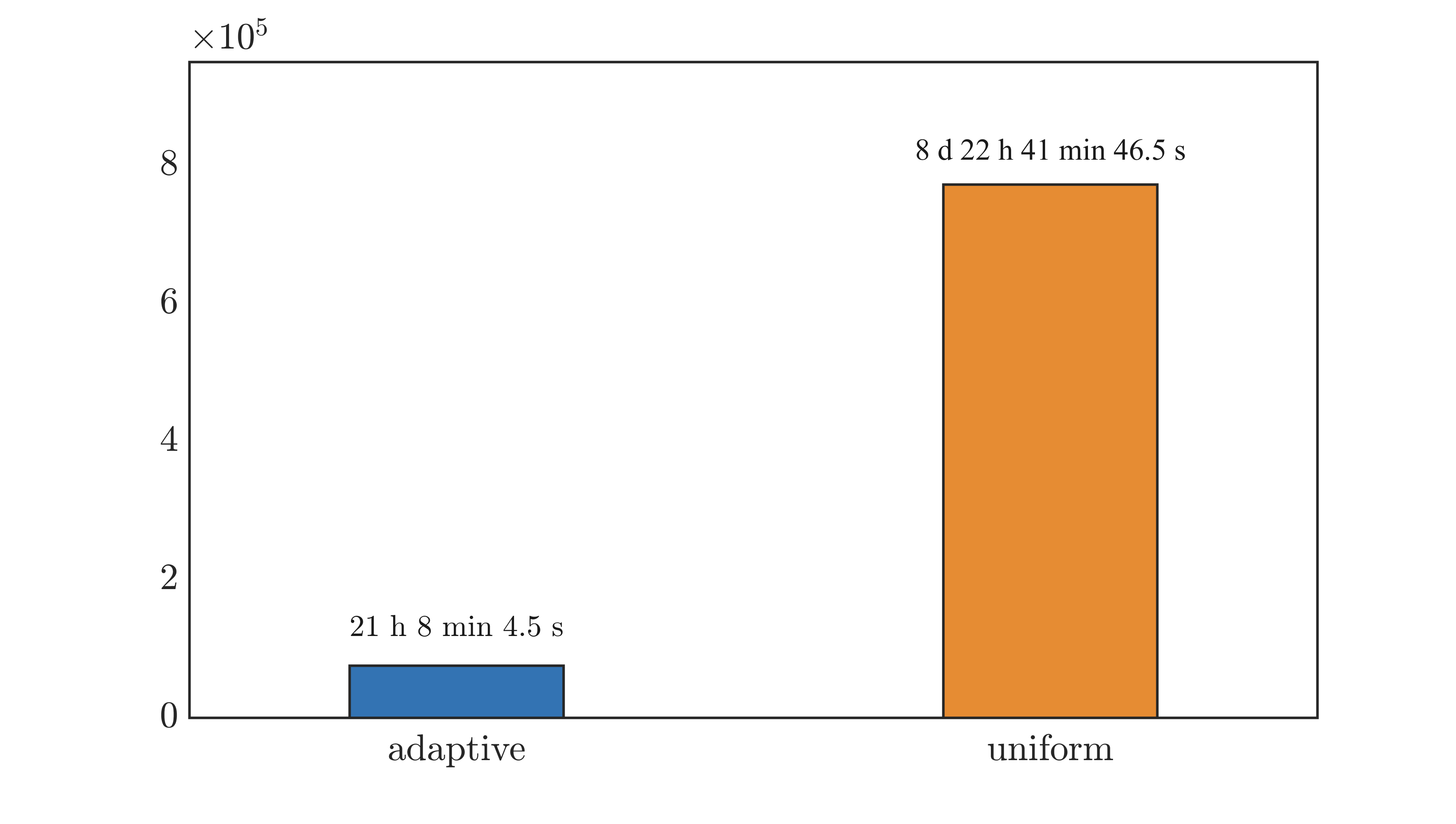}
   \hspace{-0.01\textwidth}%
    \includegraphics[width=0.45\textwidth,keepaspectratio]{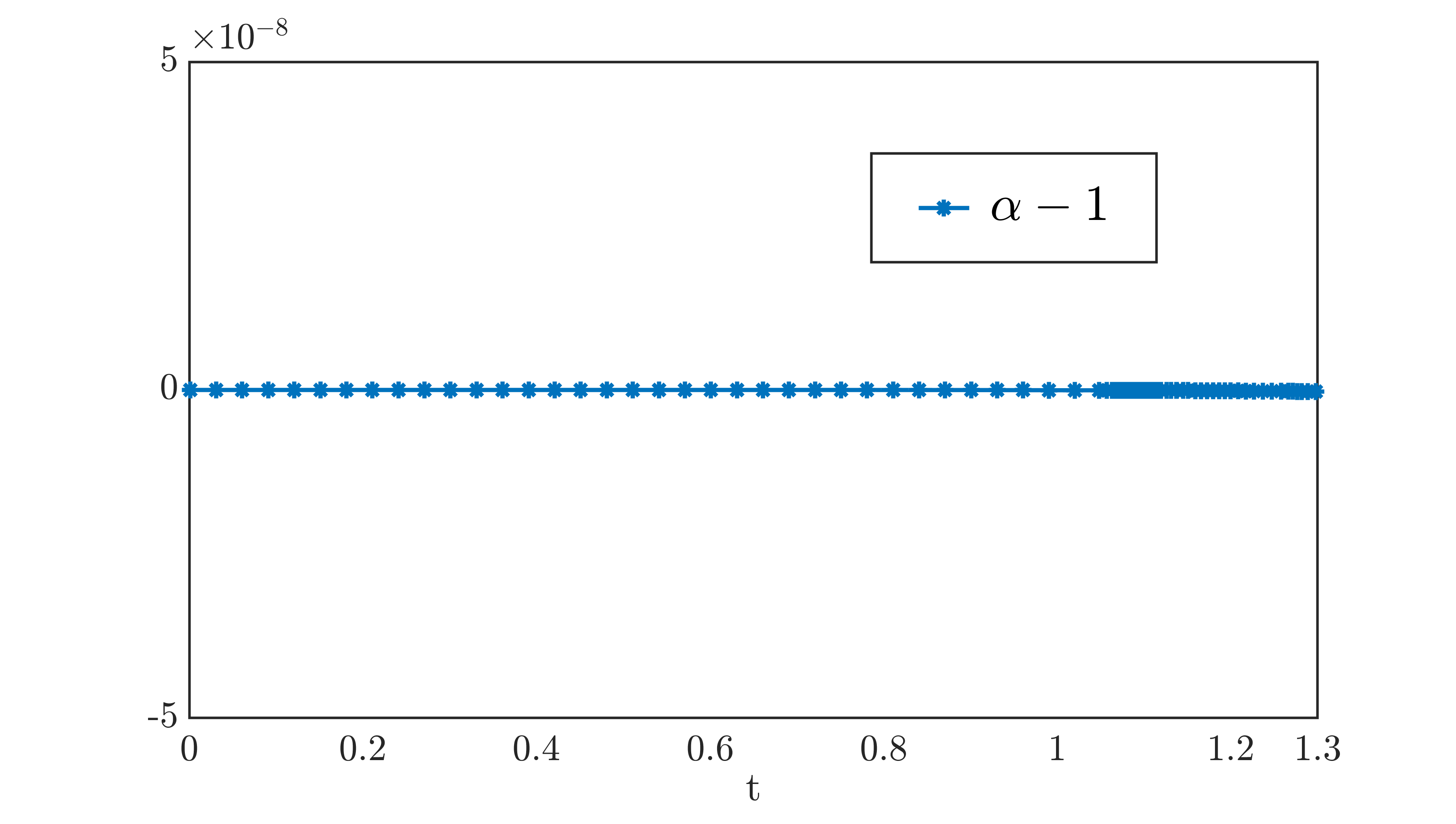}
    \hspace{-0.01\textwidth}%
    \includegraphics[width=0.45\textwidth,keepaspectratio]{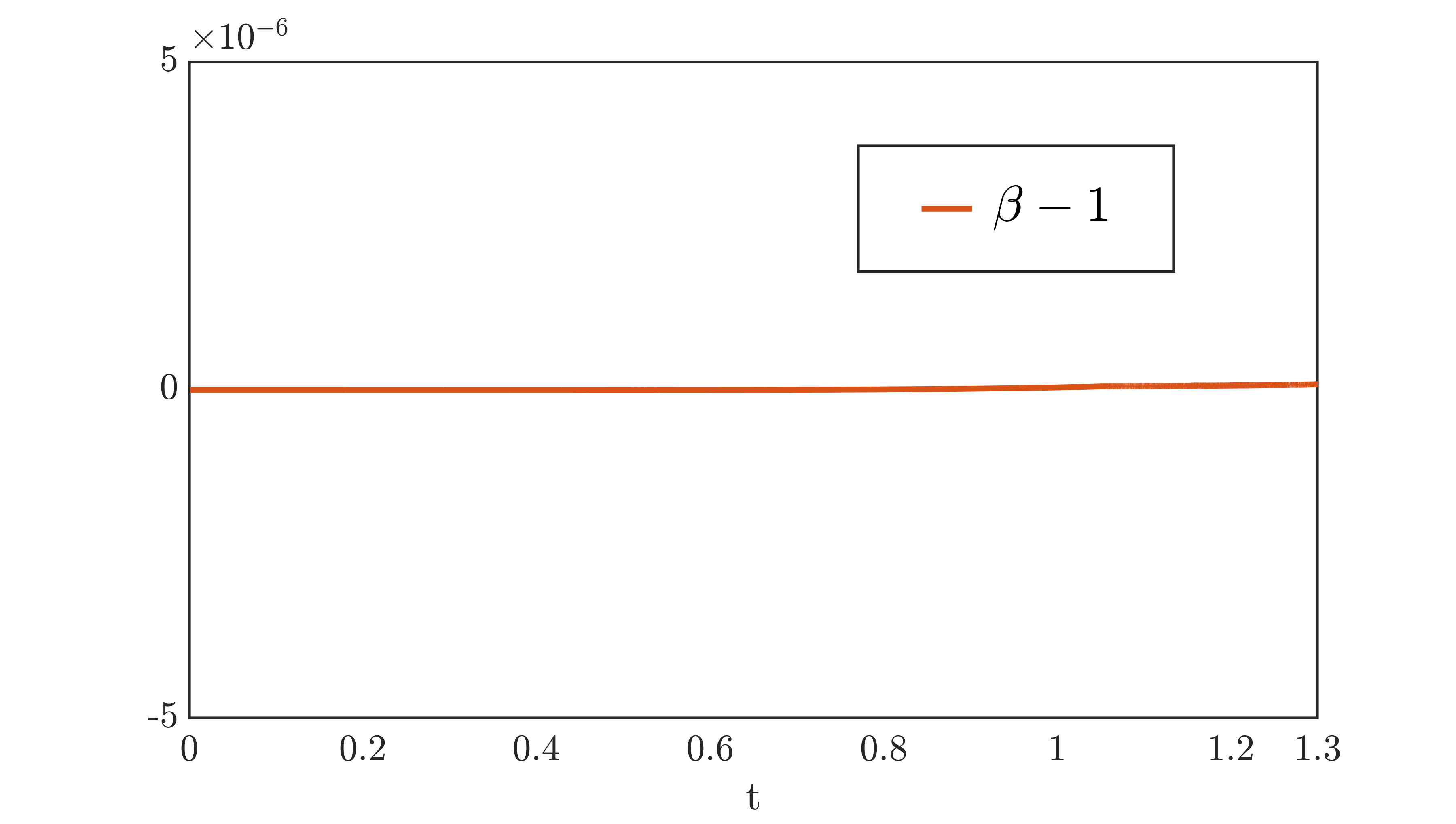}
    \hspace{-0.01\textwidth}%
    \includegraphics[width=0.45\textwidth,keepaspectratio]{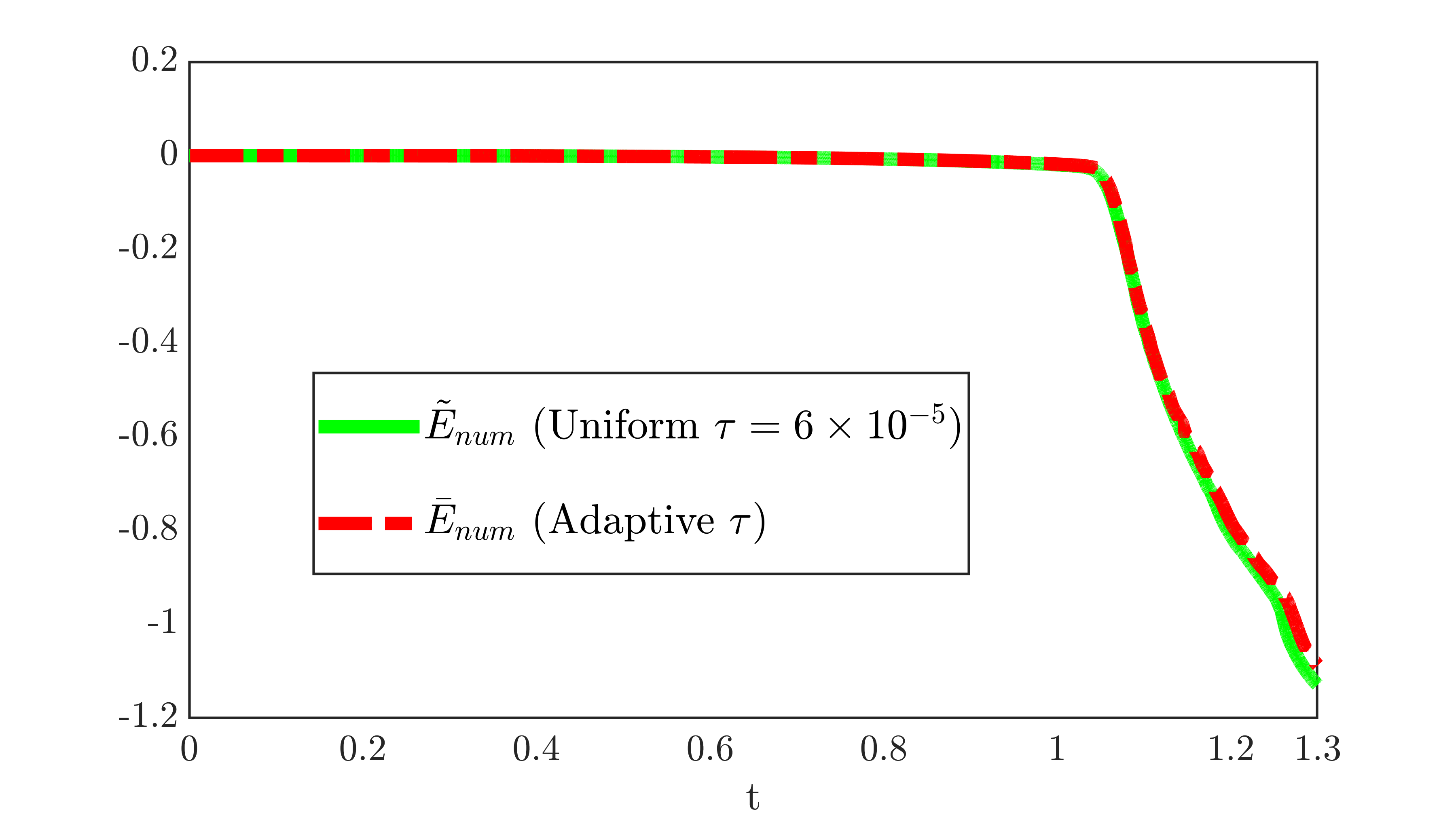}
    \hspace{-0.01\textwidth}%
    \includegraphics[width=0.45\textwidth,keepaspectratio]{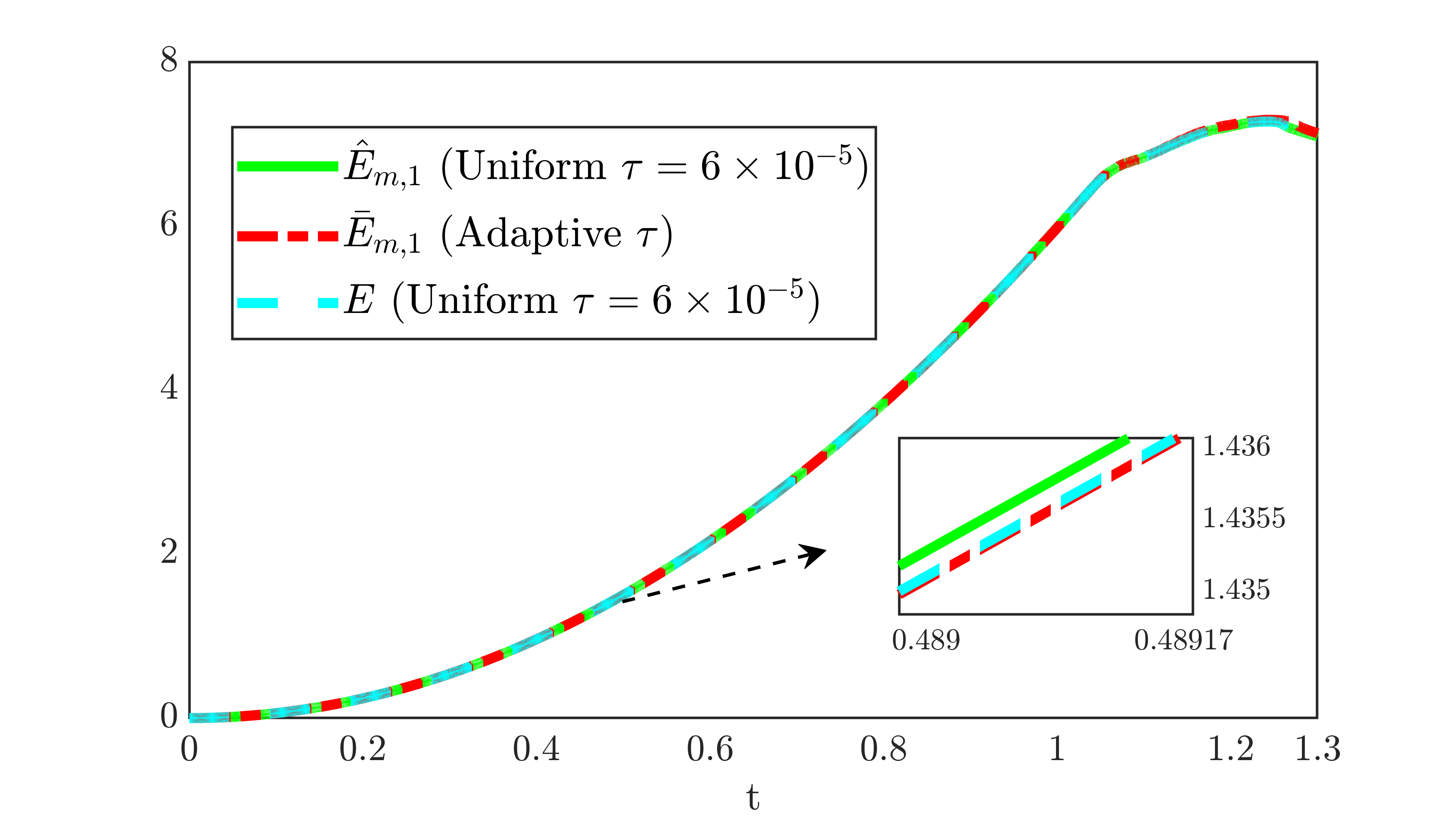}
    \hspace{-0.02\textwidth}%
    \setlength{\abovecaptionskip}{0.0cm} 
	  \setlength{\belowcaptionskip}{0.0cm}    
    \caption{
    Evolution of the adaptive time steps and CPU times, errors of $\alpha-1$ and $\beta-1$, and total and kinetic energies for the falling droplet  (from left to right, from top to bottom).}
    \label{fig:physics:falling bubbles ada}
\end{figure}
Finally, Fig.~\ref{fig:physics:falling bubbles ada} summarizes the evolution of the adaptive time steps,  CPU times,  auxiliary-variable errors $\alpha-1$ and $\beta-1$, and the discrete total and kinetic energies. Four observations follow: (i) the adaptive step remains close to $\tau_{\max}=1\times 10^{-3}$ during most of the early stage, then decreases rapidly after approximately $t=1.0$ as the droplet comes into contact with the underlying fluid interface, and subsequently varies at much smaller values; (ii) this adaptive strategy leads to a substantial reduction in computational cost, with the adaptive computation requiring about 21 hours versus nearly 9 days for the constant-step computation; (iii) throughout the simulation, the errors $\alpha-1$ and $\beta-1$ remain at the $10^{-8}$ and $10^{-6}$ levels, respectively; and (iv) the discrete total energies obtained with the uniform and adaptive time-stepping strategies nearly overlap and exhibit the same dissipative behavior. Moreover, the kinetic energies $\widehat E_{\mathrm m,1}$ computed with the constant step, $\overline E_{\mathrm m,1}$ with the adaptive step, and the original kinetic energy $E$ defined by \eqref{model:Energy} are nearly indistinguishable on the full scale; the enlarged inset shows that the adaptive $\overline E_{\mathrm m,1}$ remains closer to the original $E$, whereas the constant-step $\widehat E_{\mathrm m,1}$ exhibits a slightly larger deviation. This is consistent with Remark~\ref{rem:vart:uni}, since $\overline E_{\mathrm m,1}$ coincides exactly with the physical kinetic energy, whereas $\widehat E_{\mathrm m,1}$ is filter-dependent. 

\section{Conclusion} \label{sec:conclusion}

We proposed a second-order prediction--filter--projection framework for variable-density, variable-viscosity incompressible Navier--Stokes equations. A constrained $L^2$-projection enforces the density bounds and
mass conservation, while the GL reformulation, D-CSAV approach, and density-weighted time filter provide a second-order velocity update with decoupled velocity components. An additional scalar auxiliary variable extends the framework to variable time steps and leads to an efficient energy-based adaptive time-stepping strategy.

Both the constant- and variable-step schemes preserve the prescribed density bounds, mass conservation, and unconditional modified-energy dissipation laws. In addition, the energy analysis avoids the viscosity-dependent time-step restriction and the uniform bound on the discrete density variation required in previous studies, with stability constants independent of $\|\nabla\nu\|_{L^\infty(\Omega)}$. In particular, the introduction of an additional SAV and its compatible temporal discretization in the variable-step scheme allows the corresponding stability analysis to be extended to nonuniform grids and yields a modified energy closer to the ``original” one than that of the constant-step scheme. 

Numerical results confirm the second-order temporal accuracy and the physics-preserving properties, and the benchmark tests demonstrate the effectiveness and efficiency of the adaptive time-stepping strategy for complex interfacial flows with large density ratios. Future work will focus on rigorous error analysis and fully discrete physics-preserving schemes.

  \printcredits


\section*{Data availability}
The FreeFem++ code used in this study is available at: \url{https://github.com/HFu20/Navier-Stokes}.

\section*{Declaration of competing interest}
The authors declare that they have no known competing financial interests or personal relationships that could have appeared 
to influence the work reported in this paper.

\section*{Acknowledgements}
This work was partially supported by the National Natural Science Foundation of China (No.~12131014), the Shandong Provincial Natural Science Foundation (No.~ZR2024MA023), and the Shandong Provincial Key Laboratory of Stochastic System Control and Scientific Computing.

\FloatBarrier

\bibliographystyle{model1-num-names}

\bibliography{Ref_BP}

@article{liJCP2021,
  author  = {Li, Maojun and Cheng, Yongping and Shen, Jie and Zhang, Xiangxiong},
  title   = {A bound-preserving high order scheme for variable density incompressible {N}avier--{S}tokes equations},
  journal = {J.~Comput.~Phys.},
  volume  = {425},
  year    = {2021},
  pages   = {109906}
}

@article{liJCP2024,
  author  = {Li, Buyang and Shen, Jie and Yang, Zongze and Zhang, Yanrong},
  title   = {Efficient energy-stable schemes for incompressible flows with variable density},
  journal = {J.~Comput.~Phys.},
  volume  = {517},
  year    = {2024},
  pages   = {113365}
}

@article{zhangjcp2025,
  author  = {Zhang, Jinpeng and Luo, Li and Wang, Xiaoping},
  title   = {A fully decoupled second-order-in-time and unconditionally energy-stable scheme for the {Cahn}--{Hilliard}--{Navier}--{Stokes} equations with variable density},
  journal = {J.~Comput.~Phys.},
  volume  = {532},
  year    = {2025},
  pages   = {113943}
}

@article{guermond2000,
  author  = {Guermond, Jean-Luc and Quartapelle, Luigi},
  title   = {A projection {FEM} for variable density incompressible flows},
  journal = {J.~Comput.~Phys.},
  volume  = {165},
  year    = {2000},
  pages   = {167--188}
}

@article{tong2024pnp,
  author  = {Tong, Fenghua and Cai, Yongyong},
  title   = {Positivity preserving and mass conservative projection method for the {Poisson}--{Nernst}--{Planck} equation},
  journal = {SIAM~J.~Numer.~Anal.},
  volume  = {62},
  year    = {2024},
  pages   = {2004--2024}
}

@book{facchinei2003vol2,
  author    = {Facchinei, Francisco and Pang, Jong-Shi},
  title     = {Finite-Dimensional Variational Inequalities and Complementarity Problems, Volume~II},
  publisher = {Springer},
  address   = {New York},
  year      = {2003}
}

@article{espinozaM2AN2025,
  author  = {Espinoza-Contreras, N. and Barrenechea, G.~R. and Castillo, E. and Pacheco, D.},
  title   = {Unconditionally stable, linearised {IMEX} schemes for incompressible flows with variable density},
  journal = {ESAIM Math.~Model.~Numer.~Anal.},
  volume  = {59},  
  year    = {2025},
  pages   = {2739--2761}
}

@article{barrenecheaSISC2024,
  author  = {Barrenechea, G. and Castillo, E. and Pacheco, D.},
  title   = {Implicit-explicit schemes for incompressible flow problems with variable viscosity},
  journal = {SIAM J.~Sci.~Comput.},
  volume  = {46},
  year    = {2024},
  pages   = {A2660--A2682}
}

@article{cheng2022II,
  author  = {Cheng, Q. and Shen, J.},
  title   = {A new {Lagrange} multiplier approach for constructing structure preserving schemes, {II}. Bound preserving},
  journal = {SIAM J.~Numer.~Anal.},
  volume  = {60},
  number  = {3},
  year    = {2022},
  pages   = {970--998}
}

@article{cheng2022I,
  author  = {Cheng, Q. and Shen, J.},
  title   = {A new {Lagrange} multiplier approach for constructing structure preserving schemes, {I}. Positivity preserving},
  journal = {Comput.~Methods Appl.~Mech.~Eng.},
  volume  = {391},
  year    = {2022},
  pages   = {114585}
}

@article{yang2021IJNME,
  author  = {Yang, X.},
  title   = {A novel fully decoupled scheme with second-order time accuracy and unconditional energy stability for the {Navier}--{Stokes} equations coupled with mass-conserved {Allen}--{Cahn} phase-field model of two-phase incompressible flow},
  journal = {Int.~J.~Numer.~Methods Eng.},
  volume  = {122},
  number  = {5},
  year    = {2021},
  pages   = {1283--1306}
}

@article{yang2021,
  author  = {Yang, X.},
  title   = {A novel decoupled second-order time marching scheme for the two-phase incompressible {Navier}--{Stokes}/{Darcy} coupled nonlocal {Allen}--{Cahn} model},
  journal = {Comput.~Methods Appl.~Mech.~Eng.},
  volume  = {377},
  year    = {2021},
  pages   = {113597}
}

@article{KUANG2025,
title = {Strang splitting structure-preserving high-order compact difference schemes for nonlinear convection diffusion equations},
journal = {Commun. Nonlinear Sci. Numer. Simul.},
volume = {146},
pages = {108749},
year = {2025},
author = {Baolin Kuang and Shusen Xie and Hongfei Fu},
}

@book{temam2001,
  title={Navier-Stokes Equations: Theory and Numerical Analysis},
  author={Temam, Roger},
  year={2001},
  publisher={AMS Chelsea Publishing},
  address={Providence, RI}
}

@article{hecht2012,
  title={New development in {FreeFem++}},
  author={Hecht, Fr{\'e}d{\'e}ric},
  journal={J. Numer. Math.},
  volume={20},
  number={3--4},
  pages={251--266},
  year={2012},
  publisher={De Gruyter}
}

@article{dong2012,
  author  = {Dong, Suchuan and Shen, Jie},
  title   = {A time-stepping scheme involving constant coefficient matrices for phase-field simulations of two-phase incompressible flows with large density ratios},
  journal = {J.~Comput.~Phys.},
  volume  = {231},
  year    = {2012},
  pages   = {5788--5804}
}

@article{kuzmin2020,
  author  = {Kuzmin, Dmitri and Quezada de Luna, Manuel},
  title   = {Subcell flux limiting for high-order {B}ernstein finite
             element discretizations of scalar hyperbolic conservation laws},
  journal = {J.~Comput.~Phys.},
  volume  = {411},
  year    = {2020},
  pages   = {109411}
}

@article{lundgrenNazarov2024,
  author  = {Lundgren, Lukas and Nazarov, Murtazo},
  title   = {A high-order residual-based viscosity finite element method for incompressible variable density flow},
  journal = {J.~Comput.~Phys.},
  volume  = {497},
  pages   = {112608},
  year    = {2024}
}

@article{songCMAME2026,
  author  = {Song, Yanfei and Chen, Tanming and Zhang, Tong and Yue, Pengtao and Yang, Xiaofeng},
  title   = {Time-filtered second- and third-order accurate in time, fully decoupled and unconditionally energy-stable schemes for variable-density phase-field model of incompressible two-phase flows},
  journal = {Comput.~Methods Appl.~Mech.~Eng.},
  volume  = {457},
  pages   = {118980},
  year    = {2026}
}

@article{guzel2018time,
  author  = {G{\"u}zel, Ahmet and Layton, William},
  title   = {Time filters increase accuracy of the fully implicit method},
  journal = {BIT Numer. Math.},
  volume  = {58},
  number  = {2},
  pages   = {301--315},
  year    = {2018}
}

@article{decaria2020time,
  author  = {DeCaria, Victor and Layton, William and Zhao, Haiyun},
  title   = {A time-accurate, adaptive discretization for fluid flow problems},
  journal = {Int. J. Numer. Anal. Model.},
  volume  = {17},
  number  = {2},
  pages   = {254--280},
  year    = {2020}
}

@article{li2023filtered,
  author  = {Li, Ning and Wu, Jilian and Feng, Xinlong},
  title   = {Filtered time-stepping method for incompressible {Navier--Stokes} equations with variable density},
  journal = {J.~Comput.~Phys.},
  volume  = {473},
  pages   = {111764},
  year    = {2023}
}

@article{liao2023,
  author  = {Liao, Hong-lin and Tang, Tao and Zhou, Tao},
  title   = {Stability and convergence of the variable-step time filtered backward {E}uler scheme for parabolic equations},
  journal = {BIT Numer. Math.},
  volume  = {63},
  pages   = {39},
  year    = {2023}
}

@article{yangSISC2021,
  author  = {Yang, Xiaofeng},
  title   = {On a Novel Fully Decoupled, Second-Order Accurate Energy Stable Numerical Scheme for a Binary Fluid--Surfactant Phase-Field Model},
  journal = {SIAM J.~Sci.~Comput.},
  volume  = {43},
  number  = {2},
  pages   = {B479--B507},
  year    = {2021}
}

@article{liaoZhang2021,
  author  = {Liao, Hong-lin and Zhang, Zhimin},
  title   = {Analysis of adaptive {BDF2} scheme for diffusion equations},
  journal = {Math.~Comp.},
  volume  = {90},
  number  = {329},
  pages   = {1207--1226},
  year    = {2021}
}

@article{jiang2022improving,
  title={Improving the accuracy and consistency of the scalar auxiliary variable {(SAV)} method with relaxation},
  author={Jiang, Maosheng and Zhang, Zengyan and Zhao, Jia},
  journal={J. Comput. Phys.},
  volume={456},
  pages={110954},
  year={2022},
  publisher={Elsevier}
}

@article{liu2024novel,
  title={A Novel Energy-Optimized Technique of {SAV-based (EOP-SAV)} Approaches for Dissipative Systems},
  author={Liu, Zhengguang and Zhang, Yanrong and Li, Xiaoli},
  journal={J. Sci. Comput.},
  volume={101},
  pages={38},
  year={2024},
  publisher={Springer}
}

@article{zhang2024novel,
  title={A novel efficient generalized energy-optimized exponential {SAV scheme with variable-step BDFk} method for gradient flows},
  author={Zhang, Bingyin and Zhou, Chengxi and Fu, Hongfei},
  journal={Appl. Numer. Math.},
  volume = {210},
  pages = {39-63},
  year={2025},
  publisher={Elsevier}
}

@article{doanjcp2025,
  author  = {Doan, Cao-Kha and Hoang, Thi-Thao-Phuong and Ju, Lili and Lan, Rihui},
  title   = {Dynamically regularized Lagrange multiplier schemes with energy dissipation for the incompressible {Navier}--{Stokes} equations},
  journal = {J.~Comput.~Phys.},
  volume  = {521},
  year    = {2025},
  pages   = {113550}
}

@article{wangmathcomp2000,
  author  = {Wang, Cheng and Liu, Jian-Guo},
  title   = {Convergence of gauge method for incompressible flow},
  journal = {Math.~Comp.},
  volume  = {69},
  number  = {232},
  year    = {2000},
  pages   = {1385--1407}
}

@article{linjcp2019,
  author  = {Lin, Lianlei and Yang, Zhiguo and Dong, Suchuan},
  title   = {Numerical approximation of incompressible {Navier}--{Stokes} equations based on an auxiliary energy variable},
  journal = {J.~Comput.~Phys.},
  volume  = {388},
  year    = {2019},
  pages   = {1--22}
}

@article{kuang2025unconditional,
  author  = {Kuang, Baolin and Fu, Hongfei and Lan, Rihui and Ju, Lili},
  title   = {Unconditional energy stable and mass-conserving
             {Gauge--Uzawa} {MSAV} methods for the
             {Cahn--Hilliard--Navier--Stokes} model},
  journal = {J.~Comput.~Phys.},
  volume  = {549},
  year    = {2026},
  pages   = {114614}
}

@article{zhang2024unified,
  title={A unified framework of the {SAV-ZEC method for a mass-conserved Allen--Cahn} type two-phase ferrofluid flow model},
  author={Zhang, Guo-Dong and He, Xiaoming and Yang, Xiaofeng},
  journal={SIAM~J.~Sci.~Comput.},
  volume={46},
  number={2},
  pages={B77--B106},
  year={2024},
  publisher={SIAM}
}

@book{constantin1988navier,
  title     = {{Navier--Stokes Equations}},
  author    = {Constantin, Peter and Foias, Ciprian},
  publisher = {University of Chicago Press},
  address   = {Chicago},
  year      = {1988},
  series    = {Chicago Lectures in Mathematics},
  isbn      = {9780226115498}
}

@article{kwak2005computational,
  title   = {Computational challenges of viscous incompressible flows},
  author  = {Kwak, Dochan and Kiris, Cetin and Kim, Chang Sung},
  journal = {Comput.~Fluids},
  volume  = {34},
  number  = {3},
  pages   = {283--299},
  year    = {2005}
}

@article{shenXuYang2018,
  title   = {The Scalar Auxiliary Variable ({SAV}) Approach for Gradient Flows},
  author  = {Shen, Jie and Xu, Jie and Yang, Jiang},
  journal = {J.~Comput.~Phys.},
  volume  = {353},
  pages   = {407--416},
  year    = {2018}
}

@article{KuangFuLi2026,
  author  = {Kuang, Baolin and Fu, Hongfei and Li, Xiaoli},
  title   = {Error estimates of linear decoupled structure-preserving incremental viscosity splitting methods for the {Cahn--Hilliard--Navier--Stokes} system},
  journal = {ESAIM Math. Model. Numer. Anal.},
  volume    = {60},
  number  = {5},  
  pages     = {2449--2487},  
  year    = {2026}
}

@article{shapiroDrikakis2005,
  title   = {Artificial Compressibility, Characteristics-Based Schemes for Variable Density, Incompressible, Multi-Species Flows. {Part {I}: Derivation} of Different Formulations and Constant Density Limit},
  author  = {Shapiro, Evgeniy and Drikakis, Dimitris},
  journal = {J.~Comput.~Phys.},
  volume  = {210},
  number  = {2},
  pages   = {584--607},
  year    = {2005}
}

@article{liu2007convergence,
  title={Convergence of numerical approximations of the incompressible {Navier--Stokes} equations with variable density and viscosity},
  author={Liu, Chun and Walkington, Noel J},
  journal = {SIAM~J.~Numer.~Anal.},
  volume={45},
  number={3},
  pages={1287--1304},
  year={2007},
  publisher={SIAM}
}

@article{liShen2020SAVMAC,
  title   = {Error Analysis of the {SAV-MAC} Scheme for the {Navier--Stokes} Equations},
  author  = {Li, Xiaoli and Shen, Jie},
  journal = {SIAM~J.~Numer.~Anal.},
  volume  = {58},
  number  = {5},
  pages   = {2465--2491},
  year    = {2020}
}

@article{linNiYangDong2020,
  title   = {An Energy-Stable Scheme for Incompressible {Navier--Stokes} Equations with Periodically Updated Coefficient Matrix},
  author  = {Lin, Lianlei and Ni, Naxian and Yang, Zhiguo and Dong, Suchuan},
  journal = {J.~Comput.~Phys.},
  volume  = {418},
  pages   = {109624},
  year    = {2020}
}

@article{elAmrani2024,
  title   = {Error Estimates for a Viscosity-Splitting Scheme in Time Applied to {Non-Newtonian} Fluid Flows},
  author  = {El-Amrani, Mofdi and Obbadi, Anouar and Seaid, Mohammed and Yakoubi, Driss},
  journal = {Comput.~Methods Appl.~Mech.~Eng.},
  volume  = {419},
  pages   = {116639},
  year    = {2024}
}



\end{document}